\documentclass[aap]{imsart}

\RequirePackage{amsthm,amsmath,amsfonts,amssymb}
\RequirePackage[authoryear]{natbib}
\RequirePackage[colorlinks,citecolor=blue,urlcolor=blue]{hyperref}
\RequirePackage{graphicx}

\usepackage{mathrsfs}
\usepackage{dsfont}
\usepackage{mathtools}
\usepackage[dvipsnames]{xcolor}

\startlocaldefs
\mathtoolsset{showonlyrefs} 

\theoremstyle{plain}
\newtheorem{thm}{Theorem}[section]
\newtheorem{prop}[thm]{Proposition}
\newtheorem{lem}[thm]{Lemma}

\theoremstyle{definition}

\newtheorem{assu}[thm]{Assumption}

\newtheorem{rmk}[thm]{Remark}

\newcommand{\R}{\mathbb{R}}
\newcommand{\N}{\mathbb{N}}
\newcommand{\Z}{\mathbb{Z}}
\renewcommand{\d}{d} 
\renewcommand{\L}{\mathcal{L}}
\newcommand{\U}{\mathcal{U}}
\newcommand{\X}{\mathfrak{X}}
\newcommand{\W}{\mathfrak{W}}

\newcommand{\T}{\mathbb{T}}
\renewcommand{\S}{\mathbb{S}}

\renewcommand{\P}{\mathbb{P}}
\newcommand{\E}{\mathbb{E}}
\newcommand{\PP}[1]{\mathbb{P}\left\{#1\right\}}

\newcommand{\1}{\mathds{1}}

\newcommand{\mce}{\mathcal{E}}
\newcommand{\mcq}{\mathcal{Q}}
\newcommand{\mcb}{\mathcal{B}}
\newcommand{\mci}{\mathcal{I}}

\newcommand{\mcs}{\mathcal{S}}
\newcommand{\mcp}{\mathcal{P}}
\newcommand{\mck}{\mathcal{K}}

\newcommand{\cerny}[1]{\v Cern\'y}

\endlocaldefs

\begin{document}
\begin{frontmatter}
\title{The maximal displacement of radially symmetric branching random walk in $\R^{\lowercase{d}}$}
\runtitle{The maximal displacement of a multidimensional BRW}

\begin{aug}
\author[A, C]{\fnms{Viktor}~\snm{Bezborodov}\ead[label=e1]{viktor.bezborodov@pwr.edu.pl}},

\author[B]{\fnms{Nina}~\snm{Gantert}\ead[label=e2]{gantert@ma.tum.de}}
\address[A]{University of Göttingen,  Institute for Mathematical Stochastics}

\address[B]{Technical University of Munich\printead[presep={,\ }]{e2}}

\address[C]{Wrocław University of Science and Technology, Faculty of Information and Communication Technology\printead[presep={,\ }]{e1}}
\end{aug}

\begin{abstract}
	We consider discrete-time branching random walks with a radially symmetric distribution. Independently of each other individuals generate offspring whose relative locations are given by a copy of a radially symmetric point process $\L$. The number of particles at time $t$ form a supercritical  Galton-Watson process. 
We investigate the maximal distance to the origin of such branching random walks. Conditioned on survival, we show that, under some assumptions on $\L$, it grows in the same way as for branching Brownian motion or a broad class of one-dimensional branching random walks:
the first term is linear in time and the second logarithmic. The constants in front of these terms are explicit and depend only on the mean measure of  $\L$ and dimension. Our main tool in the proof is a ballot theorem with moving barrier which may be of independent interest. 
\end{abstract}

\begin{keyword}[class=MSC]
\kwd[Primary ]{
	60J80	}
\kwd{60K35
}
\kwd[; secondary ]{60G50}
\kwd[, ]{82C22}
\end{keyword}

\begin{keyword}
\kwd{branching random walk}
\kwd{ballot theorem}
\kwd{maximal displacement}
\end{keyword}

\end{frontmatter}
\tableofcontents


\section{Introduction}
A (binary) branching Brownian motion is a system of particles where particles branch independently into two daughter particles, and all particles move independently as Brownian motions. The study of its extrema in dimension $d=1$ is a classical topic going back to
\cite{KolmogorovPetrovskiiPiskunov1937, McKean1975}. In the seminal paper \cite{Bramson1978}, it was proved that
the maximum $M_t$ of a one-dimensional branching Brownian motion satisfies
$$
M_t = m_t +O(1) \text{ where } m_t = \sqrt{2}t - \frac{3}{2\sqrt{2}}\ln t\, .
$$
Much more is know, namely that the law of $M_t - m_t$ converges to a random shift of a Gumbel law, see \cite{LalleySellke1987}.
Branching random walks are systems of particles that reproduce and move according to independent copies of a point process $\L$, centered at the location of the mother particle.
Often, the case of i.i.d.~increments is studied, where the particles generate offspring according to a branching process and the offspring particles take i.i.d.~increments from their mother particle.

Branching random walks are used in mathematical biology as a model for the evolution of a population in space. Many natural populations live in a spatially extended habitat, with a range that is much larger than the typical distance that any
	individual may travel during its lifetime. Branching random walks on finite graphs are also used to describe the evolution of a population in terms of its genotypes.
	We refer to the book \cite{BW21}
	and in particular Chapters 2, 14 and 16 therein.

The study of the maximal distance to the origin $R_t$
of  branching random walk has a long history.
The linear growth of the first term 
was established
in 1970s \cite{Hammersley74, Big76} (see also \cite{Big97} for a more general result):
under certain integrability assumptions 
there exists $\gamma > 0$
satisfying 
$$\lim\limits _{t \to \infty} \frac{R_t}{t} = \gamma.$$
Over the next few decades it was established
that the fluctuations
of the maximum depend on
the finer structure of the 
mean measure $\mu$ of $\L$ (see \eqref{mu-def} for the definition of the mean measure).
The focus was mainly on the one-dimensional 
case $d = 1$.
If the support of $\mu$
is bounded from one side
by a unit mass atom,
then it can happen that
the second term in the representation of the maximum is of order $\frac{\ln \ln t}{\ln 2}$. This was proven in
\cite{BramsonBRW}
under the assumption of i.i.d.~increments, i.e. if
\begin{equation}\label{indepdisplacements}
	\L \overset{d}{=} \sum\limits _{i=1} ^N \delta _{Y_i},
\end{equation}
where the $\{Y_i\}$ are i.i.d.~random variables and independent of  the $\Z_+$-valued random variable $N$.
Throughout we adopt the convention that $\N = \{1,2,3, ...\}$
and $\Z _+ = \{0,1,2,3, ...\}$.
For instance, this is the case if $\mu(\{0\}) = 1$, 
$\mu(\{-1\}) = \mu(\R)-1 > 0$.
Related results and
further extensions can be found in \cite{Dur79, DH91}.
In contrast, 
for a class of branching random walks 
that can be reduced with a linear transformation to the so called 
`critical case', that is,
there exist $\beta >0$, $\gamma \in \R $ such that
\begin{equation}
 \mu (\R) > 1, \ \ \ 
	\int\limits _\R e^{-\beta (x-\gamma )} \mu (dx) = 1, 
	\ \ \ 
	\int\limits _\R (x-\gamma) e^{-\beta (x-\gamma )} \mu (dx) = 0,
\end{equation}
the maximum $M_t$ is of order $\widetilde{m}_t = \gamma t - \frac{3}{2\beta} \ln t$. 
A similar result was established in	\cite{McD95}; a preliminary discussion 
appears already in \cite{Dur79}. 
It took several years until further progress was made \cite{ABR09} and \cite{HS09}, and the convergence in law of $M_t - \widetilde{m}_t$ was proved for branching random walks in \cite{aidekon2013convergence}.
Extrema of branching Brownian motions and branching random walks in dimension $d=1$ are by now well-understood.
There are fewer results in dimension $d\geq 2$. Let $R_t^{(d)}$ be the maximal distance of a branching Brownian motion to the origin. It was shown by Mallein, \cite{Mal15} that 
\begin{equation}\label{malleintight}
	(R_t^{(d)} - m_t^{(d)})_{t \geq 0} \text{ is tight, where } m_t^{(d)} =  \sqrt{2}t + \frac{d-4}{2\sqrt{2}}\ln t\, .
\end{equation}
Later, Kim, Lubetzky and Zeitouni proved that $R_t^{(d)} - m_t^{(d)}$ converges in law to a random shift of a Gumbel law, see \cite{KimLubetzkyZeitouni2023}. Building on \cite{KimLubetzkyZeitouni2023} and \cite{StasinskiBerestyckiMallein2021}, 
\cite{BerestyckiKimLubetzkyMalleinZeitouni2023} investigate the corresponding extremal point process and prove convergence to a randomly shifted decorated Poisson point process.\\
There are also recent results for branching random walks in random environment. For the one-dimensional branching random walk in random environment the maximum scales as ${\gamma t + \beta B_1 \sqrt{t}}$, where $\gamma$ and $\beta$ are constants and $B_1$ is a standard normal random variable.
In fact a stronger result holds with the entire rescaled trajectory behaving like
a Brownian motion \cite[Theorem 2.1]{CD20}. Interestingly when accounting for the mean spread  conditionally on the environment
the correction is only logarithmic \cite[Theorem 2.4]{CD20}; the $\sqrt{t}$ terms comes from the randomness in the environment. For the maximum of a BRW in a time-inhomogeneous random environment see \cite{Kriechbaum}.\\
  
In this paper we make a first step in generalizing \eqref{malleintight} to branching random walks. Our main result, Theorem~\ref{thm fartherst point},
gives a similar statement for branching random walks with radially symmetric laws.
As  is common in the branching random
walk literature, a first moment estimate is used to provide an upper bound on the maximum
and a  truncated second moment method gives a lower bound.  The main new tool is a ballot theorem with moving barrier which may be of independent interest, see Theorem~\ref{thm beinhalten = include}.\\ 
  
The paper is organized as follows. 
In Section~\ref{sec 2} our assumptions and results 
are given. In Section~\ref{sec many-to-few} many-to-one and many-to-two lemmas are formulated. 
Section~\ref{sec max disp} contains the proof
of our main result, Theorem \ref{thm fartherst point}. Finally, Section~\ref{sec ballot} is devoted to the proof of our ballot theorem with a moving barrier, Theorem \ref{thm beinhalten = include}.

\section{Branching random walk}\label{sec 2}
We consider a branching random walk  (BRW) in $\R ^\d$, $\d \geq 2$. The process starts with a single particle located at the origin at time $0$.
All particles at time $t \in \Z_+$ are replaced 
at time $t+1$
by independent copies of
a finite  point process $\L$ translated to their position. 
The genealogical structure of the process
is therefore given by a  Galton--Watson tree $\T$. For an individual node $v\in \T$,
its generation is denoted by $|v|$ and its spatial position by $X (v)$.
Thus, a BRW is defined as a pair $(\T, \{X(v) \mid v \in \T \})$. \\
  
Denote by $\T _t$ the nodes of $\T$ in generation $t \in \N$,
i.e. the particles alive at time $t$.
In this paper we study the maximal distance 
of the BRW $(\T, \{X(v) \mid v \in \T \})$ to the origin,
defined as 
\begin{equation*}
	R _t = \max\{ |X(v)|: v \in \T _t \}.
\end{equation*}
Let $\theta $ be a direction in $\R ^\d$, that is,
$\theta \in \S^{\d - 1}$,  where $\S^{\d - 1}$ is the unit sphere in $\R ^\d$.
Let 
$${\X ^{(\theta)} (v) = \langle X(v) , \theta \rangle}$$
be the position
of the projection of $X(v)$
onto the line $\{ a \theta : a \in \R \}$.
The  pair
\begin{equation}\label{Umfang = extent, scope}
	(\T, \{\X ^{(\theta)} (v) \mid v \in \T \})
\end{equation}
constitutes then a one-dimensional BRW 
with the same genealogical tree $\T$. 
The point process $\L$ of the first generation 
particles can be written as a measure: $\L = \sum\limits_{|v| = 1} \delta _{X(v)}$.
An introduction to point processes can be found in \cite[Page 2 and elsewhere]{Kallenberg_Rand_meas} or \cite[Chapter 2]{LastPenroseBook}.
Note that
$$
\L(\R^d) = \sum\limits_ {|v| = 1}1.
$$
The  process in \eqref{Umfang = extent, scope}
can be seen as a BRW 
with the genealogical tree $\T$
and the particles
replaced by independent copies of the process  $\L ^{(\theta)}$
whose point locations are  the projections of the points of $\L$ onto $\{ a \theta : a \in \R \}$:
$$
\L ^{(\theta)} = \sum\limits _{|v| = 1} \delta _{\langle X(v) , \theta \rangle }.
$$
Denote by $0_\d$ the origin of $\R ^\d$.
\begin{assu}\label{sphersymm}
	The process $\L$ is spherically symmetric, i.e. all the $\L ^{(\theta)}$ have the same law. Further, $\E \left[\L(\R^d)\right]
	> 1$, $\E \left[\L(\R^d) ^2 \right] < \infty$, and $\E \big[\L (\R^\d \setminus \{0_\d\})\big]
	> 0$.
\end{assu}
  
Denote $m = \mu (\R ^\d)$, where $\mu$ is the mean measure of $\L$, given by
\begin{equation}\label{mu-def}
	\mu(A) = \E [ \L(A) ]  = \E \big[\sum\limits _{|v | = 1}\1 \{X(v) \in A\}\big], \ \ \ A \in \mathscr{B}(\R ^\d).
\end{equation} 
Then, $m= \E \left[\L(\R^d)\right]$. Assumption \ref{sphersymm} gives in particular that the Galton-Watson tree $\T$ is supercritical, i.e.
\begin{equation}\label{supercrit}
	m  \in  (1, \infty), \text{ implying that } \P\{\mathcal{S}\} > 0,
\end{equation}
where $\mathcal{S}$ is the event that 
the BRW survives or, in other words, that the tree $\T$ is infinite.
Define for $u \in \R$
\begin{align*}
	\Phi (u) & = \E \Big[ \sum\limits_{|v| = 1} e^{ u \X_v  ^{(\theta)}} \Big] \  \ \ \ \   = \E  \Big[\int\limits_{\R^\d} e^{ u x} \L ^{(\theta)} (dx)\Big],
	\\
	\big( \Phi  \big)' (u) & = \E \Big[ \sum\limits_{|v| = 1} \X_v  ^{(\theta)} e^{ u \X_v  ^{(\theta)}} \Big]  =\E  \Big[\int\limits_{\R^\d}x e^{ u x} \L ^{(\theta)} (dx)\Big],
\end{align*}
and  
\begin{equation}
	\Psi (u) = \ln \Phi (u), \qquad  \Psi ' (u) = \frac{ \Phi  ' (u)}{ \Phi (u)} .
\end{equation}
Note that due to Assumption \ref{sphersymm},
the functions
$\Phi (u)$, $\Phi '(u)$,
$\Psi (u)$, and $\Psi '(u)$
do not depend on the direction $\theta \in \S^{\d-1}$.  
The following equation plays an important role 
in the analysis of the behavior of the extremal points:
\begin{equation}\label{saggy ne soggy}
	u   \Psi ' (u ) - \Psi (u ) = 0.
\end{equation}
We will always work under 
\begin{assu}\label{lambda-existence}
	There exists $\lambda  \in (0,\infty)$
	with $\Phi (\lambda) < \infty$
	solving  equation \eqref{saggy ne soggy}.  Furthermore,
	there exists $\varepsilon > 0$
	such that
	$\Phi (\lambda +\varepsilon) < \infty$ 
	and for $\theta \in \S ^{\d-1}$
	\begin{equation}\label{crosshair}
		\E \Big[ \L (\R ^\d)  \int\limits_{\R^\d}
		e^{ (\lambda + \varepsilon) \langle x, \theta \rangle} \L  (dx)\Big] < \infty.
	\end{equation}
\end{assu}
  
Note that \eqref{crosshair} implies $\E \left[\L(\R^d) ^2 \right] < \infty$,
and that in view of Assumption \ref{sphersymm} the expectation in \eqref{crosshair}
does not depend on $\theta \in \S ^{\d-1}$.
Recall that $R_t$ is the maximal distance
of the particles in the branching random walk in generation $t \in \N$
to the origin.
We can now present our main result.
\begin{thm}\label{thm fartherst point}
	Let Assumptions \ref{sphersymm} and \ref{lambda-existence}  be satisfied.
	Let
	\begin{equation}\label{rdef}
		r_t: = \frac{\Psi (\lambda)}{\lambda}t
		+ 
		\frac{\d-4}{2\lambda} \ln t\, .
	\end{equation}
	Then,
	\begin{equation}\label{mainstatement}
		\sup\limits_{t \geq 0} \P\big\{\big|R_t - r_t| \geq y 
		\big| \mathcal{S}  \big\} \to 0  \  \mbox { as } \ y \to \infty\, .
	\end{equation}
	In other words, conditioned on survival, the laws of $(R_t - r_t)$ are tight.
\end{thm}
  
Theorem \ref{thm fartherst point} is proven in Section \ref{sec max disp}. The main  building blocks are a ballot theorem with a moving barrier (Theorem~\ref{thm beinhalten = include}),  many-to-one and many-to-two lemmas given in Section~\ref{sec many-to-few}, and 
random walk related results such as a	
local limit theorem, large deviation  bounds, and exponential tilting.	A brief description of the main steps of the proof is given on Page~\pageref{roadmap}.

\textbf{Examples.} Let us describe some examples which satisfy our assumptions. Note that apart from the inequalities $\E \left[\L(\R^d) ^2 \right] < \infty$ and \eqref{crosshair} and the rotational symmetry, both Assumptions \ref{sphersymm}
and \ref{lambda-existence} can be viewed as conditions
on the mean measure $\mu$ of $\L$. 
We note  that checking  
Assumption~\ref{lambda-existence} is facilitated by the formula
for the distribution of the projection of a random point on the unit sphere given in Lemma \ref{gully = ravine}. Fix $\kappa > 0$.
\begin{itemize}
	\item 
	An important class of branching random walks 
	are those with i.i.d.~increments, see \eqref{indepdisplacements}.
	If $Y_1$ has a continuous rotationally invariant distribution 
	with $\E[e^{ \kappa |Y_1| ^{a}}] < \infty$ for some $a>1$,
	$\E[N^2] < \infty$, and
	$\E[N ] > 1$, then 
	Assumptions \ref{sphersymm}
	and \ref{lambda-existence} are satisfied.
	\item 
    In dimensions $\d = 2,3$ 
	 under the assumptions	$\E \left[\L(\R^d)\right]
		> 1$, $\E \left[\L(\R^d) ^2 \right] < \infty$, and
	$$ \mu \big( \{x \in \R ^\d: |x| \geq r \}\big) =  Ce^{-\kappa r}$$
	for some $C> 0$ equation \eqref{saggy ne soggy} has a unique solution $\lambda >0$.  Condition \eqref{crosshair} needs to be verified separately;
	note that it is satisfied in the case of i.i.d.~increments described above. The case $\d \geq 4$ requires further analysis; the existence of $\lambda$ solving \eqref{saggy ne soggy} may
	depend on other parameters like $m$. 
	\item In the case of a compactly supported $\mu$ a solution to \eqref{saggy ne soggy} always exists;  the other conditions on $\mu$ are satisfied as well
	provided that $\mu(\R ^d) > 1$ and $\mu(\R^\d \setminus \{0_\d\}) > 0$.
	\item 
	Let us also mention here 
	Branching Brownian motion at integer times: each 
	particle moves independently like a Brownian motion
	and at rate $1$ splits independently into a random 
	number of particles following a distribution 
	$\nu$ on $\Z_+ $, $\nu (1) = 0$. 
	To ensure  \eqref{crosshair} it
	suffices to assume that $$\sum\limits _{n \in \Z_+}n^2 \nu (n) < \infty,$$
	see \cite[Corollary 6.1]{Branch1}.
	The inequality $\mu(\R ^d) > 1$ is equivalent to 
	$$\sum\limits _{n \in \Z_+}n \nu (n) >1.$$
\end{itemize}

Let us introduce some important notation used throughout the paper.
For $u \in \T _t$, 
$t \in \N$,
we denote by $X_s(u)$, $0\leq s < t$,
the spatial location of the ancestor of $u$ in generation $s$,
and by $X_t(u)$ the location of $u$.
Let $\{Q_t, t \in \Z_+\}$
be the random walk associated with $\L$,
that is, the random walk with 
$Q_0 = 0$ and 
\begin{equation*}
	\P \{Q_1 \in A \} = \frac 1m \E \big[\sum\limits _{|v|=1}\1 \{X(v) \in A\}\big] 
	= \frac 1m 
	\E [ \L(A) ], \ \ \ A \in \mathscr{B}(\R ^\d),
\end{equation*} 
and let $\{\X_t, t \in \Z_+\}$ \label{X def} be
a one-dimensional projection of $Q_t$
(since $\L$ is radially symmetric 
the choice of the direction does not affect the distribution of $\{\X_t, t \in \Z_+\}$). 
In other words, the distribution of $Q_1$
is the normalized mean measure of $\L$, namely $\frac{1}{\mu(\R^\d)} \mu$.
Note that we have
\begin{equation*}
	\Phi  (u)  = m \E \Big[e^{ u \X_1 }\Big],
	\ \ \ 
	\Phi  ' (u)  = m \E \Big[\X_1 e^{ u \X_1 }\Big], 
	\ \ \ 
	\Psi (u) = \ln \Phi  (u), 
	\ \ \ 
	\Psi ' (u) = \frac{\Phi  ' (u)}{\Phi   (u)},
\end{equation*}
and that the functions $\Phi$ and $\Psi $
are even.

In the proof of  Theorem \ref{thm fartherst point} 
we need the following
ballot theorem with moving barrier.
For other ballot theorems 
for a general random walk 
see e.g.
\cite[Lemma 2.1]{BDZ16}, \cite{Ballot_unpublished},
\cite{BR08}. The paper \cite{Mallein_interfaces}
contains
a few similar and related results
in Lemmas 3.6-3.9; see also \cite[Lemma 2.3]{BDZ16}
and \cite[Section A.2]{Shi15}.
\begin{thm}\label{thm beinhalten = include}
	Let $(S_n, n \in \Z_+)$   be a centered non-lattice
	one-dimensional random walk with $\E \big[S_1^2 \big]> 0$
	and for some $\varepsilon > 0$, $\E \big[|S_1|^{3+\varepsilon} \big] < \infty$.
	For a sequence of functions $f_n: \{ 0,1,...,n\} \to \R$
	define the `time reversals'
	$g_n(k) = f_n(n - k) - f_n(n)$
	and the running maximums
	$\bar f_n (k) = \max\limits _{0 \leq i \leq k}
	|f_n (i)|$ and 
	$\bar g_n (k) = \max\limits _{0 \leq i \leq k}
	|g_n (i)|$.
	Let $f:\Z_+ \to [0,\infty)$ be increasing and such that
	$$f(k) \geq \sup \{\bar f_n(k) 
	\vee \bar g_n(k): n \geq k \}. $$
	Assume that
	\begin{equation}\label{fib = unimportant lie}
		\sum\limits _{m \in \N } \frac{f(m)}{m^{3/2}} < \infty.
	\end{equation}
	Then there is some $C  > 0$ and $n_f \in \N$ such that for  
	$n\geq n_f$,
	$n \in \N$,  and for all $a, b \in  [0, n^{1/2}]$
	\begin{align}
		\frac{ (a+ 1) (b + 1) }{Cn^{3/2} }
	&	\leq
		\P\{S_k \geq  f _n(k)-a, k \leq n,  S _n \in [f_n(n) -a+ b , f_n(n)-a + b + 1 ]  \}
		\\
		\label{verschleiern = disguise}
	 &	\leq \frac{C (a+ 1) (b + 1)}{n^{3/2} }.
	\end{align}
	The constant $C$ in \eqref{verschleiern = disguise}
	depends on $f$ and the distribution of $S_1$.
\end{thm}
  
The applications of 
Theorem \ref{thm beinhalten = include}
are not restricted to branching random walks. Let us give some remarks.
\begin{rmk}
	Theorem \ref{thm beinhalten = include} is formulated
	for a non-lattice random walk $(S_n, n \in \Z_+)$.
	It can be seen in the proof that 
	the statement of the theorem holds
	for lattice random walks too, 
	provided that the length
	of the `terminal' interval 
	at least matches
	the lattice step size $\delta$;
	that is, the interval
	$[f_n(n) -a+ b , f_n(n)-a + b + 1 ]$
	in \eqref{verschleiern = disguise} 
	is replaced by 
	$[f_n(n) -a+ b , f_n(n)-a + b +  \delta ]$
	if $\delta > 1$.
\end{rmk}
\begin{rmk}
	The constant $C$
	in Theorem \ref{thm beinhalten = include}
	can always be taken so large that
	the second inequality in \eqref{verschleiern = disguise}
	holds for all $n \in \N$.
\end{rmk}
\begin{rmk}\label{ueppig = lush, lavish}
	We will also need to apply Theorem \ref{thm beinhalten = include} to a modification of the random walk 
	$(S_n, n \in \Z_+)$. Let $t \in \N$,
	$\widetilde S _0  = \widehat S _0 = 0$
	and for $s \geq 1$, 
	$$\widetilde S _s = \W + S_{s-1}, \ \ \   
	\widehat S _s = S_{s \wedge (t-1)} + \W \1 \{s=t\},
	\ \ \ s = 0, 1, ..., t,$$
	where $\W$ is a random variable with  $\E \big[|\W|^2\big] < \infty$
	independent of $(S_n, n \in \Z_+)$. The process
	$(\widetilde S_n, n \in \Z_+)$ 
	can be seen as a random walk with the same 
	step distribution as $(S_n, n \in \Z_+)$
	but with the first step following a different distribution. 
	Then \eqref{verschleiern = disguise} remains valid
	for $(\widetilde S_n, n \in \Z_+)$ and $(\widehat S_n, n \in \Z_+)$: 
	under the assumptions of Theorem \ref{thm beinhalten = include}
	there exists $C  > 0$ such that for  
	sufficiently large
	$t=n \in \N$  and for all $a, b \in  [0, n^{1/2}]$
	\begin{align*}
		\frac{ (a+ 1) (b + 1) }{Cn^{3/2} }
	&	\leq
		\P\{\widetilde S_k \geq  f _n(k)-a, k \leq n, \widetilde S _n \in [f_n(n) -a+ b , f_n(n)-a + b + 1 ]  \}
		\\
	&	\leq \frac{C (a+ 1) (b + 1)}{n^{3/2} }
	\end{align*}
	and 
	\begin{align}%
		\frac{ (a+ 1) (b + 1) }{Cn^{3/2} }
	&	\leq 
		\P\{\widehat S_k \geq  f _n(k)-a, k \leq n, \widetilde S _n \in [f_n(n) -a+ b , f_n(n)-a + b + 1 ]  \}
		\\
	&	\leq  \frac{C (a+ 1) (b + 1)}{n^{3/2} }. 
	\label{verschleiern = disguise modified}
	\end{align}
	This can be seen by conditioning on the first step
	of $(\widetilde S_n, n \in \Z_+)$ or the $t$-th step of
	$(\widehat S_n, n \in \Z_+)$, i.e.\ on $\W$.
	Note that for $n < t$ the probabilities 
	in the middle coincide in \eqref{verschleiern = disguise} and \eqref{verschleiern = disguise modified}.
\end{rmk}
  
Theorem \ref{thm beinhalten = include}
is proven in Section \ref{sec ballot}.
In the proof of Theorem \ref{thm beinhalten = include} we rely on another theorem 
which is a minor extension of certain parts of Theorem 3.2 in \cite{PP95}.
In a more general set-up asymptotic tails for the crossing moment are given in \cite{DAW18}.
The difference to \cite[Theorem 3.2]{PP95}
is that the right hand sides of the bounds explicitly involve $a$.
The statement that  \eqref{a hill to die on} and \eqref{being a trooper} hold uniformly in $a \in [0,n^{1/2}]$
	does not seem to be covered in the literature.
\begin{thm}\label{thm homing device = guidance system}
	Let $f: \N \to [0,+\infty) $ be an increasing function and $(S_n, n \in \Z_+)$ be a random walk as in Theorem \ref{thm beinhalten = include}. 
	Assume that 
	\eqref{fib = unimportant lie} holds.
	Then there exists $C  > 0$ 
	and $n_f \in \N$
	such that 
	\begin{equation}\label{a hill to die on}
		\P\{ S_k  \geq f(k) - a \text{ for }
		n_f \leq k\leq n  \} \geq \frac {(a+1)}{C\sqrt{n}},
		\ \ \ a \in [0,n^{1/2}], \, \, n \geq n_f,
	\end{equation}
	and 
	\begin{equation}\label{being a trooper}
		\P\{ S_k  \geq -f(k) - a \text{ for }
		1 \leq k\leq n  \} \leq \frac {C(a+1)}{\sqrt{n}},
		\ \ \ a \in [0,n^{1/2}], \, \, n \geq 1.
	\end{equation}
\end{thm}
  
The proof of Theorem \ref{thm homing device = guidance system}
can be found in Section~\ref{sec ballot}.
The beginning of that section also contains a very brief description of the proof ideas.\\
  
We conclude this section with a remark on notation. Throughout the paper $C$ and sometimes $c$ refer  to 
positive
constants whose exact value is not important
and can change from line to
line and even within the same line. The constants  cannot depend on time or on the realization of the BRW
but can depend 
on the distribution 
of the involved processes, for example on $\L$ or $S_1$.

\section{Many-to-one and many-to-two lemmas}\label{sec many-to-few}
The following lemma
relates the 
BRW $(\T, \{X(v) \mid v \in \T \})$
and the random walk
$\{Q_t, t \in \Z_+\}$.
More specifically
a sum over 
the particles in generation $n$
can be expressed in terms of 
a single path of the random walk.
For a particle $v \in \T$
with $|v| = n$ denote by $v_i$
its ancestor in generation $i$,
so that $v_0$ is the root of $\T$,
and $v_n = v$. 
\begin{lem}[Many-to-one]\label{lem manyto1}
	Let $n \in \N$ and let $f: (\R ^\d)^n \to \R$
	be a bounded measurable function. It holds that
	\begin{equation}
		\E\Big[\sum\limits _{|v| = n}
		f(X(v_0), X (v_1), ..., X(v_n) )\Big]
		= m^n\E \big[f(Q_0, Q_1, ..., Q_n)\big].
	\end{equation}
\end{lem}
  
The lemma is obvious in the case of i.i.d.~increments. It can be proved in the same way as Theorem 1.1. in \cite{Shi15}, taking $a=0$ there.\\
  
Next we present a many-to-two lemma 
which helps us to deal with the second moment.
A similar computation appears in the proof 
of (1.6) in \cite{BRW_Min_SimpleProof}.
The many-to-two lemma will be used
for deriving bounds
on the second moment
of random variables of the form 
$$
\# \{  |u| = t: X(u) \geq g(t)-1, \  \forall s \leq t: X(u_s) \leq g(s)\},
$$
where $t \in \N$ and $g:\{0,1,...,t\} \to \R_+$
is a certain function.\\
  
Set $m_2 = \E \big[\L(\R^d)^2\big] - \E\big[\L(\R^d)\big]  = 
\E \big[\L(\R^d)^2\big] - m $. Note that $m_2 >0$ since $\E\big[\L(\R^d)\big] > 1$.
Let $\Delta = (\Delta_1, \Delta_2)$ 
be a $(\R^\d)^2$-valued exchangeable pair of random variables
(i.e. $(\Delta_1, \Delta_2)  \overset{d}{=} (\Delta_2, \Delta_1)$)
with distribution
\begin{equation}\label{morose, sulky, sullen}
	\P \{  (\Delta_1, \Delta_2) \in A \}
	= m_2 ^{-1}\E \Big[ \sum\limits _{\substack{|u| = |v| = 1 \\ u \ne v}} \1 
	\{(X(u), X(v)) \in A  \}  \Big] 
\end{equation}
where $A \in \mathscr{B}(\R ^{2\d})$
is symmetric, that is, for $x, y \in \R ^\d$, 
$(x,y ) \in A$ if and only  if $(y,x) \in A$.\\
  
We note here that for $a > 0$ whenever $\Phi(a) < \infty$
we also have
$$\E\big[e^{a | \langle Q_1, \mathbf{e}_1 \rangle|}\big]
 \leq
\E\big[e^{-a  \langle Q_1, \mathbf{e}_1 \rangle} + 
e^{a  \langle Q_1, \mathbf{e}_1 \rangle} \big]
= 
 2 m^{-1} \Phi (a) < \infty$$
and that whenever $\E \big[ \L (\R ^\d)  \int_{\R^\d}
e^{ a \langle x, \theta \rangle} \L  (dx)\big] < \infty $
we have 
\begin{align}
	m_2 \E \big[e^{a \Delta _1}\big] = \E \Big[ \sum\limits _{\substack{|u| = |v| = 1 \\ u \ne v}} e^{a X(v)}  \Big] 
	= 
	\E \Big[ (\L(\R ^\d) - 1) \sum\limits _{ |v| = 1 } e^{a X(v)}  \Big] < \infty.
\end{align}
Let  $\{Q ' _t, t \in \Z_+\}$
and  $\{Q '' _t, t \in \Z_+\}$
be copies 
of  $\{Q_t, t \in \Z_+\}$
and let the walks 
$\{Q  _t, t \in \Z_+\}$,
$\{Q ' _t, t \in \Z_+\}$,
$\{Q '' _t, t \in \Z_+\}$,
the BRW $(\T, \{X(v) \mid v \in \T \})$,
and the pair $\Delta$ be independent.
Define for $i,k  \in \Z_+$
\begin{equation}\label{die Salbe = ointment}
	Q^{\langle k \rangle}_i = \begin{cases}
		Q_i,  & \text{ if }  i \leq k,
		\\
		Q_k + \Delta _1 + Q ' _{i-k-1}, & \text{ if }  i > k
	\end{cases}
\end{equation}
and similarly
\begin{equation} \label{selbstsuchtig = selfish}
	Q^{[k]}_i = \begin{cases}
		Q_i,  & \text{ if }  i \leq k,
		\\
		Q_k + \Delta _2 + Q '' _{i-k-1}, & \text{ if }  i > k.
	\end{cases}
\end{equation}
Note that in law the walk 
defined by \eqref{die Salbe = ointment}
can be written as 
\begin{equation}\label{zu meinen Gunsten}
	\{ Q^{\langle k \rangle}_i , i \in \Z_+ \}\overset{d}{=} \{  Q_{i - \1 \{i > k  \}} +
	\Delta _1  \1 \{i > k  \} , i \in \Z_+ \},
\end{equation}
and a similar equality holds for $ \{ Q^{[k]}_i , i \in \Z_+ \}$.
Thus $ \{ Q^{\langle k \rangle}_i , i \in \Z_+ \}$
is a random walk with the same
step 
distribution as  $ \{ Q_i , i \in \Z_+ \}$
but with a single step replaced
by $\Delta _1$.
\begin{lem}[Many-to-two] \label{lem many-to-two}
	Let $n \in \N$ and let $f: (\R ^\d)^{2n} \to \R$
	be a bounded measurable function. 
	It holds that
	\begin{align}
		\E 	\Big[\sum\limits _{v,u \in \T: |v| = |u| = n}
		f(X(v_0), & \, X (v_1), ..., X(v_n), X(u_0), X (u_1), ..., X(u_n) )\Big] \notag 
		\\
		= \ &  \label{Skrei}
		m^n\E \Big[f(Q_0, Q_1, ..., Q_n, Q_0 , Q_1, ..., Q_n)\Big]
		\\   & + m_2 m^{2n-2}
		\sum\limits _{k = 0} ^{n-1} m^{-k} \E\Big[ f(Q_0 ^{\langle k \rangle}, 
        Q_1 ^{\langle k \rangle}, ..., Q_n ^{\langle k \rangle}, Q^{[k]} _0  ,  Q^{[k]} _1 , ..., Q^{[k]} _n )\Big].
		\notag
	\end{align}
\end{lem}

   The proof of Lemma~\ref{lem many-to-two} 
can be found in the appendix on Page~\pageref{proof of many-to-two}.

\section{Maximal distance to the origin: Proof of Theorem \ref{thm fartherst point}}
\label{sec max disp}
To streamline the proof we split 
the statement of Theorem \ref{thm fartherst point}
into two parts which are handled separately.
Recall that
$
r_t = \frac{\Psi (\lambda)}{\lambda}t
+ 
\frac{\d-4}{2\lambda} \ln t
$.
\begin{prop}\label{prop fartherst point upper bound}
	Under the assumptions of Theorem \ref{thm fartherst point} we have
	\begin{equation}\label{mainstatement upper}
		\sup\limits_{t \geq 0} \P\big\{R_t   \geq r_t + y 
		\big| \mathcal{S}  \big\} \to 0  \  \mbox { as } \ y \to \infty\, .
	\end{equation}
\end{prop}
\begin{prop}\label{prop fartherst point lower bound}
	Under the assumptions of Theorem \ref{thm fartherst point} we have
	\begin{equation}\label{mainstatement lower}
		\sup\limits_{t \geq 0} \P\big\{R_t   \leq r_t - y 
		\big| \mathcal{S}  \big\} \to 0  \  \mbox { as } \ y \to \infty\, .
	\end{equation}
\end{prop}
  
We note that Propositions \ref{prop fartherst point upper bound} and \ref{prop fartherst point lower bound}
together
give Theorem \ref{thm fartherst point}.
The proof of Proposition \ref{prop fartherst point upper bound} can be found on Page \pageref{Proof prop fatherst point upper bound}. It is preceded by the auxiliary Lemma  \ref{anschalten = switch on}. Here we rely on the many-to-one lemma, 
an observation \cite[Lemma 2.3]{Mal15} of how many small caps are needed to cover the unit sphere, and
bounds in the large deviation regime.\\
  
We have to work much harder in order to complete the proof of
Proposition \ref{prop fartherst point lower bound}.
A roadmap of the proof can be found on Page \pageref{roadmap}.
The auxiliary statements  are given by Lemmas \ref{hoodwink = deceive, trick} to \ref{umbrage; take an umbrage new}
and Proposition \ref{volcano plume}. 
The final part of the proof  of Proposition \ref{prop fartherst point lower bound}
is located on 
Page \pageref{Proof prop fartherst point lower bound}.\\
  
We now introduce some notation and collect some preliminary material 
used throughout this section.
Set $\X_t = \langle Q_t, \mathbf{e}_1 \rangle $, 
where $\mathbf{e}_1 = (1,0,...,0) \in \S ^{\d-1}$.
Define then a random walk 
$\{\X^\perp_t, t \in \Z_+\}$
by 
$\X^\perp_t = Q_t - \X_t \mathbf{e}_1$,
that is, $\X^\perp_t$ is the orthogonal
projection of $Q_t$ on the hyperplane orthogonal
to $\mathbf{e}_1$. Note that $|Q_t|^2 = \X^2_t + |\X^\perp _t|^2$.

Recall that $\lambda > 0$ is the solution
to \eqref{saggy ne soggy}.
For $s,t \in \N$, $s \leq t$, $s \in \Z_+$
and $y \geq 0$ set 
\begin{equation}\label{fdef}
	f_s^{t,y} = \frac{\Psi (\lambda)}{\lambda}s+ \frac{\d-1}{2\lambda} \ln (s+1)
	- \frac{3}{2\lambda} \ln \frac{t+1}{t-s+1} +
	\frac{3}{\lambda} M + \frac 1\lambda y,
\end{equation} 
where $M > 1 \vee  \lambda$ is sufficiently large so that
\begin{equation}\label{die Ampel}
	{\Psi (\lambda)}s+ \frac{\d-1}{2} \ln (s+1)
	- \frac{3}{2} \ln \frac{t+1}{t-s+1} +
	M \geq 0,  \ \ \ s, t\in \Z_+, s\leq t,
\end{equation} 
and
\begin{equation*}
	\P \Big\{  |Q _1| \leq \frac{M}{\lambda}  \Big\} \geq \frac 12.
\end{equation*}
We also note here that
by \eqref{die Ampel}
for $s \in \Z_+, t \in \N$, $s \leq t$, 
$y \geq 0$,
\begin{equation}\label{der Platz}
	f_s^{t,y} \geq  \frac{2}{\lambda} M.
\end{equation}
It will be important later
that
the family of maps $\{f_s^{t,0} - \frac{\Psi (\lambda)}{\lambda}s, 0\leq s \leq t \} _{t \in \N}$ satisfies the assumptions 
of Theorem \ref{thm beinhalten = include}
for example  with $f(s) = \frac{\d+2}{2\lambda} \ln (s+1) + \frac{3}{\lambda}M$.
Indeed,
\begin{equation*}
	\Big|f_s^{t,0} - \frac{\Psi (\lambda)}{\lambda}s - \frac{3}{\lambda}M \Big|
	=  \Big|\frac{\d-1}{2\lambda} \ln (s+1)
	- \frac{3}{2\lambda} \ln \Big( 1 + \frac{s}{t-s+1}\Big) \Big|
	\leq \frac{\d+2}{2\lambda} \ln (s+1)
\end{equation*}
and 
\begin{align*}
	\Big|f_{t-s}^{t,0} - \frac{\Psi (\lambda)}{\lambda}(t-s)
	-  & \, f_t^{t,0} - \frac{\Psi (\lambda)}{\lambda}t  \Big|
	\\
	= & \, \Big|
	\frac{\d-1}{2\lambda} \ln (t-s+1)
	-
	\frac{3}{2\lambda} \ln \frac{t+1}{s+1}
	- \frac{\d-1}{2\lambda} \ln (t+1)
	+ \frac{3}{2\lambda} \ln (t+1)
	\Big|
	\\
	\leq  & \, \frac{\d-1}{2\lambda} \Big| \ln \Big( 1 + \frac{s}{t-s+1}\Big) \Big|
	+  \frac{3}{2\lambda} \ln (s+1)
	\leq  \frac{\d+2}{2\lambda} \ln (s+1).
\end{align*}
By Jensen's inequality 
$$
\Phi( u ) =  m \E \big[e^{ u \X_1 }\big] \geq  m e^{ \E [u \X_1] } = m,
$$
and hence $\Psi(\lambda ) \geq  \ln m = \Psi(0) $.

For $b \in \R$
with $\Phi(b)< \infty$
denote by $(\X^{(b)}_s, s\in \Z_+ )$,
$\X^{(b)}_0 = 0$,
the random walk 
with distribution 
\begin{equation}\label{Verwahrung = custody, safekeeping}
	\P \{\X^{(b)}_1 \in A \} = 
	\frac{\E \big[e^{b \X_1} \1 \{\X_1 \in A\}\big]   }{
		\E \big[e^{b \X_1}\big] },
	\ \ \ 
	A \in \mathscr{B}(\R ).
\end{equation}
The change of measure  \eqref{Verwahrung = custody, safekeeping}
is well-known
and is discussed in more detail 
in \cite[Section 2.3.1]{srw}.
Note that 
\begin{equation}\label{erheblich = ziemlich viel}
	\E \big[\X^{(b)}_1\big] = \frac{\E \big[\X_1 e^{b \X_1}\big] }{
		\E  \big[e^{b \X_1}\big] }
	= \frac{\Phi' (b)}{\Phi (b)} = \Psi ' (b),
\end{equation}
\begin{equation}\label{Verwahrung = custody, safekeeping2}
	\P \{(\X_1, ..., \X_n) \in A \}
	= \E \big[e^{b \X_1}\big]^n \E \left[ \frac{
		\1 \{(\X^{(b)}_1, ..., \X^{(b)}_n) \in A\} }{e^{b \X^{(b)}_n}}\right],
	\ \ \ 
	A \in \mathscr{B}(\R ^n),
\end{equation}
or equivalently for a measurable function $f: \R ^n \to \R$
\begin{equation*}
	\E [f(\X_1, ..., \X_n) ]
	= \E \big[e^{b \X_1}\big]^n \E  \left[\frac{
		f(\X^{(b)}_1, ..., \X^{(b)}_n)}{e^{b \X^{(b)}_n}}\right].
\end{equation*}
Here and throughout we adopt the convention that 
$\E [A]^b = \big(\E [A] \big)^b$.\\
  
We also note that 
$\E \big[e^{b \X_1}\big] = m^{-1} \Phi  (b) = 
m^{-1}
\E \big[\sum\limits _{|v| = 1} e^{b \langle X(v), e_1 \rangle } \big]$,
and that
by \eqref{saggy ne soggy} and \eqref{erheblich = ziemlich viel} $$
\E \big[\X_1^{(\lambda)}\big] = \Psi'(\lambda) = 
\frac{\Psi(\lambda)}{\lambda},
$$
and that $\X_1^{(\lambda)}$ has an exponential moment
by Assumption \ref{lambda-existence} and \eqref{Verwahrung = custody, safekeeping}.
In particular, for $\varepsilon > 0$ from Assumption \ref{lambda-existence}
\[
\E \Big[ e ^{\varepsilon | \X^{(\lambda)}_1 | } \Big] = \frac{\E \big[ e ^{\lambda \X_1 +  \varepsilon | \X_1 | }  \big] }{\E \big[e ^{\lambda \X_1   }\big]} 
\leq \frac{\E \big[ e ^{(\lambda  +  \varepsilon) \X_1} + e ^{(\lambda  - \varepsilon) \X_1  }  \big] }{\E \big[e ^{\lambda \X_1   }\big]} 
\leq \frac{2 \E \big[ e ^{(\lambda + \varepsilon)  \X_1  }  \big] }{\E \big[e ^{\lambda \X_1   }\big]} < \infty
\]
where we used that the map $ \R_+ \ni q \mapsto \E \big[e ^{q \X_1  }\big]  = \E \big[ \big(e ^{\X_1  } \big)^q\big]$ is increasing
by Hölder's inequality
 because the distribution of $\X_1$
is symmetric and non-degenerate: 
for $0 < q_1 < q_2$
$$
\E \Big[ e ^{q_1 \X_1  } \Big] < 
\bigg(\E \Big[ e ^{ \frac{q_2}{q_1} q_1 \X_1  } \Big] \bigg) ^{\frac{q_1}{q_2}} \cdot 1
= \bigg(\E \Big[ e ^{  q_2 \X_1  } \Big]\bigg) ^{\frac{q_1}{q_2}} 
\leq  \E \Big[ e ^{  q_2 \X_1  } \Big].
$$

Finally 
	let us 
	note that the random walk $\{\X_t, t \in \Z_+\}$
	satisfies the non-lattice condition in
	Theorem \ref{thm beinhalten = include},
	since apart from a possible atom at zero $\X_1$  has a continuous distribution as
	 an orthogonal projection of a spherically-symmetric random walk:
	 the distribution of $\X_1$ can be found in \eqref{take a rain check=decline}
	 in the appendix.
	This also applies to $\{\X^{(b)}_t, t \in \Z_+\}$ whenever the latter is well-defined. 

  \subsection{Proof of Proposition \ref{prop fartherst point upper bound}}
The following lemma contains the bulk of the arguments
leading to the proof of Proposition \ref{prop fartherst point upper bound}.
\begin{lem}\label{anschalten = switch on}
	There exists $C> 0$ such that 
	for any $t \in \N$ and $y \in [0 ,\sqrt{t}]$.
	\begin{equation}
		\P \left\{ \exists s\leq t, u \in \T _s:
		|X_s(u)| \geq f_s^{t,y} \right\} \leq 
		C (1+y) e^{-y}.
	\end{equation}
\end{lem}
\begin{proof} We have 
\begin{align}
	\notag
	\P \{  \exists  s\leq t, & u \in \T _s:
	|X_s(u)| \geq f_s^{t,y} \} 
	\\ \notag
	&
	=
	\P \Big( 
	\bigcup\limits _{k = 1} ^{t} 
	\Big\{ \exists  u \in \T _k: |X_s(u)|  <  f_s^{t,y}, s \leq k-1, 
	|X_k(u)| \geq f_k^{t,y} \Big\}   \Big) 
	\notag
	\\ \notag
	& \leq \P \Big( \bigcup\limits _{k = 1} ^{t} 
	\Big\{\exists  u \in \T _k: |X_s(u)| <  f_s^{t,y}, s \leq k-1, 
	|X_{k-1}(u) | + | X_{k}(u) - X_{k-1}(u) | \geq f_k^{t,y} \Big\}    \Big) 
	\notag
	\\
	& \leq \sum\limits 
	_{k = 1} ^{t}
	\P \Big\{ \exists  u \in \T _k: 
	|X_s(u)| <  f_s^{t,y}, s \leq k-1, 
	|X_{k-1}(u) | +  | X_{k}(u) - X_{k-1}(u) | \geq f_k^{t,y} \Big\} 
	\notag
	\\
	& \leq \sum\limits 
	_{k = 1} ^{t}
	\E\Big[ \sum\limits _{u \in \T _k} \1 \Big\{ 
	|X_s(u)| <  f_s^{t,y}, s \leq k-1, 
	|X_{k-1}(u) | +  | X_{k}(u) - X_{k-1}(u) | \geq f_k^{t,y} \Big\} \Big]
	\notag
\end{align}
\vspace{-0.25cm}
\begin{align}
	\mkern-35mu = \sum\limits 
	_{k = 1} ^{t}
	m^{k}
	\P \big\{|Q_s| <  f_s^{t,y}, s \leq k-1, 
	|Q_{k-1}| + |Q_k - Q_{k-1}| \geq f_k^{t,y} \big\}
	\label{encore = additional performance (bis)}
\end{align}
where we used the many-to-one formula (Lemma \ref{lem manyto1}) for the last equality. 
For $r \in [0, f_k^{t,y}] $
set $\phi_k (r) = \P \big\{|Q_s| <  f_s^{t,y}, s \leq k-1, 
|Q_{k-1}| + r \geq f_k^{t,y} \big\} $,
so that 
\begin{multline}\label{felt is cloth}
	\P \big\{|Q_s| <  f_s^{t,y}, s \leq k-1, 
	|Q_{k-1}| + |Q_k - Q_{k-1}| \geq f_k^{t,y} \big\}
	\\
	= 
	\E\big[\phi _k (|Q_k - Q_{k-1}|)\big]  = 
	\E\big[\phi _k (|Q_1|)\big] 
\end{multline}
Lemma 2.3 in \cite{Mal15} 
gives, for each  $R > 1$, the existence of 
a finite subset $\U (R)$ of the unit sphere 
satisfying 
$$ \sup\limits _{R>1} \frac{\#\U (R)}{R^{\frac{\d-1}{2} }} < \infty
\ \ \  
\text{ and }
\  \ \
\S ^{\d-1} \subseteq \bigcup \limits _{\theta \in \U (R)}
\Big\{  \vartheta \in \S ^{\d-1}: \langle \theta, \vartheta \rangle   \geq 1 - \frac 1 R \Big\}.$$
Therefore
\begin{multline*}
	\{|Q_s| <  f_s^{t,y}, s \leq k-1, 
	|Q_{k-1}| + r \geq f_k^{t,y} \}
	\\
	\subseteq 
	\mkern-9mu
	\bigcup\limits _{\theta \in \U (f_k^{t,y})}
	\mkern-10mu
	\{\langle Q_s, \theta \rangle <  f_s^{t,y}, s \leq k-1, 
	\langle Q_{k-1}, \theta \rangle  \geq f_k^{t,y} - r - 1 \},
\end{multline*}
and
hence 
\begin{align}
	\phi _k (r)
	& \leq 
	\P \Big( \bigcup\limits _{\theta \in \U (f_k^{t,y})}
	\{\langle Q_s, \theta \rangle <  f_s^{t,y}, s \leq k-1, 
	\langle Q_{k-1}, \theta \rangle  \geq f_k^{t,y} - r -1 \}
	\Big) 
	\notag
	\\
	&  \leq \# \U (f_k^{t,y})
	\P \{\X_s <  f_s^{t,y}, s \leq k-1, 
	\X_{k-1}  \geq f_k^{t,y} - r -1 \}
	\notag
	\\
	& \leq  C (1+k)^{\frac{\d-1}{2} }
	\P \{\X_s <  f_s^{t,y}, s \leq k-1, 
	\X_{k-1}  \geq f_k^{t,y} - r -1  \}.
	\label{fudge is a sweet thing}
\end{align}
Consider the random walk $Y_s = -\X ^{(\lambda)} _s + \frac{\Psi (\lambda)}{\lambda}s$. Note that $\E [Y_1] = 0$.
Denote $$G ^t _s =   \frac{\Psi (\lambda)}{\lambda}s - f_s^{t,0}
= \frac{\Psi (\lambda)}{\lambda}s - f_s^{t,y} + \frac y \lambda .$$
Barring a possible atom at zero, the distribution of $Y_s$ 
	is continuous and  therefore
 non-lattice, $1 \leq s \leq t$.
Applying Theorem \ref{thm beinhalten = include}
we get 
for $n \in \N $ and $r > 0$
\begin{align}
	\P \{ 
	\X ^{(\lambda)} _s \leq  \, & f_s^{t,y}, s \leq n, 
	\X ^{(\lambda)} _{n}  \geq f_n^{t,y} - r -1 \}
	\\
	=  & \,
	\P \big\{ 
	Y _s \geq    \frac{\Psi (\lambda)}{\lambda}s - f_s^{t,y}, s \leq n, 
	Y _{n}  \leq  \frac{\Psi (\lambda)}{\lambda}n
	- f_n^{t,y}  + r +1 \big\}
	\\ \label{limpid = clear, transparent} 
	= & \,  \P \big\{ 
	Y _s \geq   G ^t _s - \frac y \lambda, s \leq n, 
	Y _{n}  \leq G ^t _n - \frac y \lambda + r +1 \big\}
	\\
	\leq  & \, 
	\sum\limits _{\ell = 0 } ^{ \lceil r + 1 \rceil }
	\P \big\{ 
	Y _s \geq  G ^t _s - \frac y \lambda, s \leq n, 
	Y _{n}   \in [  G ^t _n - \frac y \lambda +  \ell , G ^t _n - \frac y \lambda +  \ell + 1 ] \big\}
	\\
	\leq  & \, C  \sum\limits _{\ell = 0 } ^{ \lceil r + 1 \rceil } (1+y)(1+\ell) (n+1) ^{-3/2}
	\\
	\leq  & \,
	C (r+2)^2 (1+y) (n+1) ^{-3/2} .
\end{align} 
With this in mind and since 
$$
e^{\lambda f_k^{t,y}  } =   e^{\Psi (\lambda) k}
\left( k+1\right)^{\frac{\d - 1}{2 }} \left( t-k+1\right)^{\frac{3}{2 }} 
e^{3M+y} (t+1)^{-3/2} 
$$
the change of measure \eqref{Verwahrung = custody, safekeeping} with $b = \lambda$  yields
\begin{align}
	\P \{ 
	\X_s <  f_s^{t,y}, s & \leq  k-1, 
	\X_{k-1}  \geq f_k^{t,y} - r  -1 \} 
	\notag 
	\\ \notag
	& = \E\big[e^{\lambda \X_1}\big]^{k-1}
	\E\left[\frac{\1 \{ 
		\X ^{(\lambda)} _s <  f_s^{t,y}, s \leq k-1, 
		\X ^{(\lambda)} _{k-1}  \geq f_k^{t,y} - r -1  \}}
	{e^{\lambda \X^{(\lambda)}_{k-1}}} \right]
	\\  \notag
	& \leq
	\left( \frac{\Phi  (\lambda)}{m} \right)^{k-1}
	\frac{\P \{ 
		\X ^{(\lambda)} _s <  f_s^{t,y}, s \leq k-1, 
		\X ^{(\lambda)} _{k-1}  \geq f_k^{t,y} - r -1 \} }
	{e^{\lambda(f_k^{t,y} - r-1 )}}
	\\ \notag
	& \leq
	\left( \frac{\Phi  (\lambda)}{m} \right)^{k-1}
	\frac{
		C (y+1) \left( t+1\right)^{\frac{3}{2 }}
		(r+2)^2 e^{\lambda r} }
	{(k+1)^{3/2}  e ^{\Psi (\lambda) k}
		\left( k+1\right)^{\frac{\d - 1}{2 }}   
		\left( t-k+1\right)^{\frac{3}{2 }} 
		e^y 
	}
	\\
	&
	= \frac{
		C (y+1)  \left( t+1\right)^{\frac{3}{2 }}
	}
	{m ^{k-1} (k+1)^{3/2} \Phi   (\lambda) 
		\left( k+1\right)^{\frac{\d - 1}{2 }}   
		\left( t-k+1\right)^{\frac{3}{2 }} 
		e^y  } \times (r+2)^2 e^{\lambda r}.
	\label{in Anbetracht dessen = in view of}
\end{align}
By Assumption \ref{lambda-existence},
for some $\varepsilon > 0$
$$\E\big[ (\X_1 + 2)^2 e^{(\lambda + \varepsilon) \X_1}\big] < \infty$$
which combined with the symmetry of the law of $Q_1$ implies  
$$ \E\big[ (|Q_1| + 2)^2 e^{\lambda |Q_1|}\big] < \infty $$
and by \eqref{fudge is a sweet thing}
and \eqref{in Anbetracht dessen = in view of}
\begin{align}
	\notag
	\E\big[\phi _k (|Q_1|)\big] 
	& \leq  C (1+k)^{\frac{\d-1}{2} } 
	\frac{
		(y+1)  \left( t+1\right)^{\frac{3}{2 }}
	}
	{m ^{k-1} (k+1)^{3/2} 
		\left( k+1\right)^{\frac{\d - 1}{2 }}   
		\left( t-k+1\right)^{\frac{3}{2 }} 
		e^y  } 
	\E\big[ (|Q_1| + 2)^2 e^{\lambda |Q_1|}\big]
	\\ 
	& \leq C 
	\frac{
		(y+1)  \left( t+1\right)^{\frac{3}{2 }}
	}
	{m ^{k-1} (k+1)^{3/2}  
		\left( t-k+1\right)^{\frac{3}{2 }} 
		e^y  }. \label{zittig = shaky}
\end{align}
Combining
\eqref{encore = additional performance (bis)},
\eqref{felt is cloth},
and \eqref{zittig = shaky}
we get
\begin{align}
	\P \{ \exists s\leq t, & u \in \T _s:
	|X(u)| \geq f_s^{t,y} \} 
	\notag 
	\\
	& \leq 
	\sum\limits 
	_{k = 0} ^{t-1}
	m^{k}
	\E\big[\phi _k (|Q_1  |)\big] 
	\notag 
	\\
	& \leq  
	C\sum\limits 
	_{k = 0} ^{t-1}
	m^{k}
	\frac{
		(y+1) \left( t+1\right)^{\frac{3}{2 }} }
	{m ^{k-1} (k+1)^{3/2}   
		\left( t-k+1\right)^{\frac{3}{2 }} 
		e^y  }
	\notag 
	\\
	& \leq  
	C (y+1) e^{-y} \sum\limits 
	_{k = 0} ^{t-1}
	\frac{
		\left( t+1\right)^{\frac{3}{2 }} }
	{ (k+1) ^{\frac 32}   
		\left( t-k+1\right)^{\frac{3}{2 }} 
	}. \label{nachhaltig = sustainable}
\end{align}
It remains to note that the last sum 
in \eqref{nachhaltig = sustainable}
is uniformly bounded in $t$:
for some $C_2>0$
\begin{align}
	\sum\limits 
	_{k = 0} ^{t-1}
  &	\frac{
		\left( t+1\right)^{\frac{3}{2 }} }
	{ (k+1) ^{\frac 32}   
		\left( t-k+1\right)^{\frac{3}{2 }} 
	} 
	\\
	\label{riff-raff}
&	\leq 2 \sum\limits 
	_{k \in \Z_+: k\leq t/2 + 1} 
	\frac{
		\left( t+1\right)^{\frac{3}{2 }} }
	{ (k+1) ^{\frac 32}   
		\left( t-k+1\right)^{\frac{3}{2 }} 
	} \leq C_2
	\sum\limits 
	_{k = 0} ^{\infty}
	\frac{1 }
	{ (k+1) ^{\frac 32}   
	}.
\end{align}
\end{proof}

Recall \eqref{fdef} and note that  
$$
f_t^{t,y} = \frac{\Psi (\lambda)}{\lambda}t+ \frac{\d-4}{2\lambda} \ln (t+1) +
\frac{3}{\lambda} M + \frac 1\lambda y = r_t + \frac{3}{\lambda} M + \frac 1\lambda y.
$$

\begin{proof}[Proof of Proposition \ref{prop fartherst point upper bound}.]
	\label{Proof prop fatherst point upper bound}
	Fix $\delta > 0$.
	Let $y_0 > 3M$ be large enough so that for $y \geq y_0$
	$$
	C (1 + y - 3M) e^{-(y - 3M)} 
	\leq  \delta,
	$$
	where $C >0$ is the constant from Lemma \ref{anschalten = switch on}.
	 For $y \geq y_0$, $t _0 \geq y^2$
	 we have
	by Lemma \ref{anschalten = switch on} 
	      \begin{align*}
	      	 \sup\limits _{t \geq t_0} \P \Big\{ R_t - r_t \geq \frac 1\lambda y \Big\}
	      	 =
	      	  \sup\limits _{t \geq t_0} \P \Big\{ R_t  \geq f_t^{t,y - 3M}  \Big\}
	      	 \leq C (1 + y - 3M) e^{-(y - 3M)}  \leq \delta. 
	      \end{align*}      
	      For a large  $y_1 > 0$ we have by tightness
	      $$
	      \sup\limits _{t \geq t_0} \P \Big\{ R_t - r_t \geq \frac 1\lambda y_1 \Big\}
	     \leq  \delta;
	      $$
	      hence for $y \geq y_0 \vee y_1$
	      $$
	      \sup\limits _{t \in \N} \P \Big\{ R_t - r_t \geq \frac 1\lambda y \Big\} 
	      \leq \delta.
	      $$
	      The proof is concluded by observing that $\P (\mathcal{S}) > 0$ and 
	      $\delta >0$ is arbitrary.
\end{proof}

    \subsection{Proof of Proposition \ref{prop fartherst point lower bound}: a roadmap}
From \label{roadmap} now on and until the end of the section we are working toward establishing Proposition \ref{prop fartherst point lower bound}. We now give a short roadmap of the remaining part of the proof. 
Recall that $f^{t,y} _s$ was introduced in \eqref{fdef}.
For $\theta \in \S^{\d - 1}$  
define
\begin{equation}\label{der Teil}
	\mathcal{A}^{t,y}_\theta =  \{u \in \T _t: 
	| X_s(u)| \leq f^{t,y} _s, s \leq t, \langle X_t(u), \theta\rangle \geq f^{t,y} _t - 1 \}.
\end{equation}
Our aim is to apply the second moment method to 
$$\mathcal{N}_t = 
\sum \limits _{\theta \in \mathscr{K}_t}
\#\mathcal{A}^{t,y}_\theta$$
with   $\mathscr{K}_t$ being a finite subset of $\S ^{\d-1}$ satisfying certain conditions; in particular $\# \mathscr{K}_t \xrightarrow{t \to \infty}  \infty$. 
More specifically, we will use the inequality 
\begin{equation}
	\P \{ \mathcal{N}_t \geq 1   \}
	\geq \frac{\E\big[\mathcal{N}_t \big] ^2}{\E \big[\mathcal{N}_t ^2\big]}
\end{equation}
in the proof of Proposition \ref{volcano plume}.
For this inequality to be helpful  we need an upper bound 
on the second moment $ \E  \big[ \mathcal{N}_t ^2  \big] $
and a lower bound on the first moment $ \E  \big[ \mathcal{N}_t  \big]$. We start with an upper bound on $ \E  \big[ \# \mathcal{A}^{t,y}_\theta  \big]$
in Lemma \ref{hoodwink = deceive, trick}. 
Lemma \ref{fledgling} gives then an upper bound on $ \E  \big[ (\# \mathcal{A}^{t,y}_\theta  \big) ^2]$.
A lower bound on $ \E  \big[ \# \mathcal{A}^{t,y}_\theta  \big]$ is obtained
Lemma \ref{hoodwink = deceive, trick2}.
Lemma \ref{hoodwink = deceive, trick2} relies on  an auxiliary statement, Lemma \ref{colic}, where it is shown
that roughly speaking the function $f^{t,y}_{\cdot}$
dominates in a certain sense the sum of two functions $g_{\cdot}$ and $h_{\cdot}$, with $g$ behaving similarly to  $f^{t,y}_{\cdot}$ as far as Theorem \ref{thm beinhalten = include} is concerned and $h$ having certain desirable properties.
This immediately gives a lower bound on 
$ \E  \big[ \mathcal{N}_t  \big]$ since 
\[
\E  \big[ \mathcal{N}_t  \big] \geq \# \mathscr{K}_t
\E  \big[ \# \mathcal{A}^{t,y}_\theta  \big].
\]
We have to work much harder for the upper bound on 
$ \E  \big[ \mathcal{N}_t ^2 \big] $.
In Subsection \ref{subsec challenging} we come to probably the most challenging part of the entire proof: the upper bound on 
$\E  [\#  \mathcal{A}^{t,y}_\theta \# \mathcal{A}^{t,y}_{\theta'}]$ for $\theta \ne \theta '$.
Note that since $\# \mathscr{K}_t \xrightarrow{t \to \infty}  \infty$, we  have to treat the case 
when $\theta$ and $\theta '$ may be very close (the distance between them can be of order $O(t^{-1/2})$);
the bound has to depend on the distance between $\theta$ and $\theta '$.  An upper bound on 
$\E  [\#  \mathcal{A}^{t,y}_\theta \# \mathcal{A}^{t,y}_{\theta'}]$ is given in
Lemma \ref{umbrage; take an umbrage new}.
In the preceding Lemma \ref{pauper = poor person}
a necessary bound on a certain 
expectation containing an exponential expression 
is obtained;  Lemma \ref{pauper = poor person}
is arguably the most technical part of the proof. 
Proposition \ref{volcano plume} 
gives then  an exponential bound
used in the proof of Proposition \ref{prop fartherst point lower bound}.

\subsection{Upper bound on $\E  [(\#  \mathcal{A}^{t,y}_\theta)^2]$}
The following lemma is used in the proof of Lemma \ref{fledgling}, where an upper bound on 
$\E  [(\#  \mathcal{A}^{t,y}_\theta)^2]$ is derived.

\begin{lem}\label{hoodwink = deceive, trick}
	There exists $C>1$ such that 
	for $t \in \N$ 
	and $y \in [0, t^{1/2}]$
	\begin{equation} \label{reticence1}
		\E  \big[ \# \mathcal{A}^{t,y}_\theta  \big]
		\leq C (1+y) e^{-y} t^{-\frac{\d-1}{2}}.
	\end{equation}
\end{lem}
	\begin{proof}
		Recall that $ m \E [e^{\lambda \X _1}] = \Phi(\lambda)$.
		Using the many-to-one lemma (Lemma \ref{lem manyto1}) and 
		\eqref{Verwahrung = custody, safekeeping2}
		with $b = \lambda$
		we get 
		\begin{align}
			\E  \big[ \# \mathcal{A}^{t,y}_\theta  \big]
			&  = 
			m^t
			\P\{ 
			| \X_s| \leq f^{t,y} _s, s \leq t,  \X_t  \geq f^{t,y} _t - 1 \}
			\\
			& \leq  
			\Phi ^t (\lambda)
			e^{-\lambda f^{t,y} _t}
			\P\{ 
			| \X ^{(\lambda)}_s| \leq f^{t,y} _s, s \leq t,  \X^{(\lambda)}_t  \geq f^{t,y} _t - 1 \}
			\\
			& = (t+1) ^{-\frac{\d-4}{2}} e^{-y -3M}
			\P\{ 
			| \X ^{(\lambda)}_s| \leq f^{t,y} _s, s \leq t,  \X^{(\lambda)}_t  \geq f^{t,y} _t - 1 \}.
		\end{align}
		Applying Theorem \ref{thm beinhalten = include}
		to the random walk $Y_s = -\X ^{(\lambda)} _s + \frac{\Psi (\lambda)}{\lambda}s$
		just like in \eqref{limpid = clear, transparent}
		we get
		\[
		\P\{ 
		| \X ^{(\lambda)}_s| \leq f^{t,y} _s, s \leq t,  \X^{(\lambda)}_t  \geq f^{t,y} _t - 1 \}
		\leq C (1+y) e^{-y }(t+1) ^{-3/2}.
		\]
		The inequality in  \eqref{reticence1} follows.
	\end{proof}

Now we give an upper bound on $\E  [(\#  \mathcal{A}^{t,y}_\theta) ^2]$;
 an upper bound on 
$\E  [\#  \mathcal{A}^{t,y}_\theta \# \mathcal{A}^{t,y}_{\theta'}]$
with $\theta \ne \theta'$ 
is derived later
in Lemma \ref{umbrage; take an umbrage new}.
\begin{lem}\label{fledgling}
	There exists $C>0$
	such that for $t \in \N$, $y \in [0, t^{1/2}]$,
	and $\theta \in \S ^{\d-1}$
	\begin{equation*}
		\E  [(\#  \mathcal{A}^{t,y}_\theta)^2]
		\leq  C (1+y) e^{-y}  t^{-\frac{\d-1}{2}}.
	\end{equation*}
\end{lem}

\begin{proof}
	By the many-to-two Lemma \ref{lem many-to-two}
	\begin{align} \label{in the groove}
		\ \ \ \ \ \ \ 
		\E   [(\#  \mathcal{A}^{t,y}_\theta)^2]  & = \E \big[ \# \mathcal{A}^{t,y}_\theta  \big] 
		+  m_2 m^{2t-2} 
		\sum\limits _{s=0} ^{t-1} m ^{-s} 
		\P\Big[|Q ^{\langle s \rangle} _k|  \leq f^{t,y} _k, \langle Q ^{\langle s \rangle} _t, \theta \rangle
		\geq f^{t,y} _t - 1,  \\ & |Q^{[s]}_k|  \leq f^{t,y} _k, \langle Q^{[s]}_t, \theta \rangle
		\geq f^{t,y} _t - 1, k = 0,...,t  \Big], 
		\notag
	\end{align}
	where $\{Q ^{\langle s \rangle}_k, k = 0,1,...,t \}$ and $\{Q^{[s]}_k, k = 0,1,...,t \}$ are random walks defined
	in \eqref{die Salbe = ointment} and \eqref{selbstsuchtig = selfish}.\\
	
	Set 
	$\X^{\langle s \rangle}_k =  \langle  Q^{\langle s \rangle}_k, \theta \rangle$ and
	$\X^{[s]}_k =  \langle  Q^{[s]}_k, \theta \rangle$, and
	define for $s = 0,...,t-1$
	\begin{equation} \label{phi def}
		\phi_s(x) =  \P \{\X  _{k}  + x \leq f^{t,y} _{s+k + 1}, k \leq t-s-1, \X _{t-s-1}  + x \geq   f^{t,y} _t - 1  \}.
	\end{equation}
	By the Markov property 
	\begin{multline}\label{der Verkehr}
		\P\Big[|Q ^{\langle s \rangle}_k| \leq f^{t,y} _k, \langle Q ^{\langle s \rangle} _t, \theta \rangle
		\geq f^{t,y} _t - 1, |Q^{[s]}_k| \leq f^{t,y} _k, \langle Q^{[s]}_t, \theta \rangle
		\geq f^{t,y} _t - 1, k = 0,...,t  \Big]
		\\
		\leq 
		\P\Big[  \X ^{\langle s \rangle} _k \leq f^{t,y} _k, \X ^{\langle s \rangle} _t\geq f^{t,y} _t - 1,
		\X ^{[s]}_k \leq f^{t,y}  _k, \X ^{[s]}_t \geq f^{t,y} _t - 1, k=0,...,t  \Big]
		\\
		=  \E \left[ \phi_s   (\X_s + \W _1) \phi_s   (\X_s + \W _2)\1 \{\X_k \leq f^{t,y} _k, k \leq s  \}       \right].
	\end{multline}
	Recall that $\Psi(\lambda) = \ln \Phi (\lambda)$.
	Since the maps 
		$$\{0,1,...,t\} \ni s \mapsto f^{t,y} _{s}$$ 
		are Lipschitz uniformly in $t$ and $y$, we have by  Theorem \ref{thm beinhalten = include}  for $x \leq f^{t,y} _{s}$
	\begin{align}
		\P \{ \X ^{(\lambda)} _{k} &  \leq f^{t,y} _{s+k+1} - x, k \leq t-s-1, \X_{t-s-1} ^{(\lambda)}  \geq   f^{t,y} _t -x - 1  \}
		\\
		&	= \P \Big\{ \X ^{(\lambda)} _{k} - \frac{\Psi (\lambda)}{\lambda}k  \leq f^{t,y} _{s+k+1} - \frac{\Psi ( \lambda)}{\lambda}k - x, k \leq t-s-1,
		\\ & 
		\mkern 47mu \X_{t-s-1} ^{(\lambda)} - \frac{\Psi (\lambda)}{\lambda}(t-s-1) \geq   f^{t,y} _t - \frac{\Psi (\lambda)}{\lambda}(t-s-1) -x - 1  \Big\}
		\\ \label{confabulate}
		& \leq C \frac{1+ f^{t,y} _{s+1}   - x }{(t-s+1)^{3/2}}
		\leq C \frac{1+ f^{t,y} _{s}   - x }{(t-s+1)^{3/2}}.
	\end{align}
	Applying the transformation in \eqref{Verwahrung = custody, safekeeping},
	\eqref{Verwahrung = custody, safekeeping2} with 
	$b = \lambda $ and \eqref{confabulate}
	we get for $x \leq f^{t,y} _{s}$
	\begin{align}
		\phi_s(x ) & =    \P \{\X  _{k}  + x \leq f^{t,y} _{s+k+1}, k \leq t-s-1, \X _{t-s-1}  + x \geq   f^{t,y} _t - 1  \} \notag
		\\
		&
		= \E  [e^{\lambda \X_1}]^{t-s-1}
		\E \bigg[  \frac{ \1 \{ \X ^{(\lambda)} _{k}  \leq f^{t,y} _{s+k+1} - x, k \leq t-s-1, \X_{t-s-1} ^{(\lambda)}  \geq   f^{t,y} _t -x - 1\} }{e^{\lambda \X_{t-s-1} ^{(\lambda)}  }}
		\bigg] \notag
		\\
		&
		\leq C \left( \frac{\Phi (\lambda)}{m} \right)^{t-s}
		\frac{\P \{ \X ^{(\lambda)} _{k}  \leq f^{t,y} _{s+k+1} - x, k \leq t-s-1, \X_{t-s-1} ^{(\lambda)}  \geq   f^{t,y} _t -x - 1  \}}{\exp \{ \lambda f^{t,y} _t - \lambda x - \lambda   \}} \notag
		\\
		&
		\leq C \left( \frac{\Phi (\lambda)}{m} \right)^{t-s}
		\frac{\big(1 + f^{t,y} _{s} - x\big) e^{\lambda x} }{(t-s+1)^{3/2}}
		\times
		\frac{(t+1)^{{3}/{2}}}{(t+1)^{{(\d-1)}/{2}}e^{\Psi(\lambda)t} e^y}
		\notag
		\\
		&
		\leq
		C m^{-(t-s)} \Phi(\lambda) ^{-s}
		e^{-y}
		\frac{\big(1 + f^{t,y} _{s} - x\big) e^{\lambda x} }{(t-s+1)^{3/2}(t+1)^{{(\d-4)}/{2}}}. \label{disparate}
	\end{align}
	Therefore by \eqref{in the groove} and \eqref{der Verkehr}
	\begin{align}
		\E   [(\#  \mathcal{A}^{t,y}_\theta)^2]    \leq 
		\E \big[ \# \mathcal{A}^{t,y}_\theta \big] + 
		C & \sum\limits _{s=1} ^t
		\bigg\{  \, 
		\frac{m ^{2t-s} m^{-2(t-s)} \Phi(\lambda) ^{-2s}
			e^{-2y}}{(t-s+1)^{3}(t+1)^{{(\d-4)}}} \notag
		\\ \notag
		&	\times 
		\E\Big[\big(1  + f^{t,y} _{s+1}  - \X_s - \W _1 \big)
		\big(1  + f^{t,y} _{s+1} - \X_s - \W _2 \big)  
		e^{\lambda (2\X_s + \W_1 + \W_2 )} 
		\\ 
		&	\times 
		\1 \{\X_k \leq f^{t,y} _k, k \leq s,
		\X_k + \W _1 \vee \W_2 \leq f^{t,y} _{k+1}   \}   \Big]
		\bigg\}\, .
		\label{fete = prazdnik}
	\end{align}
	Denote by $\mathcal{E}$ the expectation on the right hand side of 
	\eqref{fete = prazdnik}.
	Applying the inequality \\
	
	${(u+a)(u+b) \leq \big(u+\frac 12 a + \frac 12 b\big)^2}$, $u,a,b \in \R$
	and using that a.s.
	$$
	\1 \{\X_k \leq f^{t,y} _k, k \leq s,
	\X_k + \W _1 \vee \W_2 \leq f^{t,y} _{k+1}  \} 
	\leq \1 \big\{\X_k \leq f^{t,y} _k, k \leq s,
	\X_k + \frac 12 \W _1  + \frac 12 \W_2 \leq f^{t,y} _{k+1}  \big\} 
	$$
	we get
	\begin{equation*}
		\mathcal{E} \leq \E \bigg[ \Big(1  + f^{t,y} _{s+1} - \X_s  - \frac{\W _1 + \W_2}{2} \Big) ^2
		\1 \big\{\X_k \leq f^{t,y} _k, k \leq s,
		\X_k + \frac 12 \W _1  + \frac 12 \W_2 \leq f^{t,y} _{k+1}  \big\}
		\bigg]
	\end{equation*}
	Recall that the random walk  $(\X_s, s\in \Z_+ )$
	and the pair $(\W_1, \W_2)$ are independent.
	Next we apply the transformation in 
	\eqref{Verwahrung = custody, safekeeping} to the  
	walk $(\X)$ and  the random variable 
	$ \frac 12 \W _1  + \frac 12 \W_2 $
	to get
	\begin{align}
		\mathcal{E} \leq \, &	
		\E \big[e^{\lambda \X_1}\big]^s \E \big[e^{\frac  \lambda 2 (\W _1  + \W_2) } \big]
		\E \Bigg[   \frac{
			\Big(1  + f^{t,y} _{s+1} - \X_s ^{(\lambda)} -  \frac{\W ^{(\lambda)}  _1 + \W ^{(\lambda)}  _2}{2} \Big) ^2
			e^{\lambda (2\X_s ^{(\lambda)} + \W ^{(\lambda)}  _1 + \W ^{(\lambda)}  _2 )}}{e^{\lambda \X^{(\lambda)}_s + \frac  \lambda 2 (\W ^{(\lambda)}  _1  + \W ^{(\lambda)}  _2)}  }
		\\ \notag 
		&  \times	 	\1 \big\{\X_k ^{(\lambda)} \leq f^{t,y} _k, k \leq s,
		\X_s ^{(\lambda)} + \frac 12 \W^{(\lambda)}  _1  + \frac 12 \W ^{(\lambda)} _2 \leq f^{t,y} _{s+1}  \big\} \Bigg]
		\\ \notag
		\leq \, &
		C \left( \frac{\Phi (\lambda)}{m} \right)^{s}
		\E \bigg[  \Big(1  + f^{t,y} _{s+1} - \X_s ^{(\lambda)} - \frac{\W ^{(\lambda)}  _1 + \W ^{(\lambda)} _2}{2} \Big) ^2 e^{\lambda (\X_s ^{(\lambda)} + \frac 12\W ^{(\lambda)}  _1 + \frac 12 \W ^{(\lambda)}  _2 )}
		\\ \notag
		&  \times \1 \big\{\X_k \leq f^{t,y} _k, k \leq s,
		\X_s + \frac 12 \W ^{(\lambda)}  _1  + \frac 12 \W ^{(\lambda)}  _2 \leq f^{t,y} _{s+1}  \big\} \bigg].
	\end{align}
	Next decomposing over the values of $ \X_s + \frac 12 \W _1  + \frac 12 \W_2$ in intervals of the form ${(f^{t,y} _{s+1} -n - 1, f^{t,y} _{s+1} - n ]}$, $n \in \Z_+$,
	and applying Theorem \ref{thm beinhalten = include}
	(taking
	into account Remark \ref{ueppig = lush, lavish} as well)
	we obtain
	\begin{align}
		\notag
		\mathcal{E} \leq \, &
		C \left( \frac{\Phi (\lambda)}{m} \right)^{s} \sum\limits _{n = 0} ^\infty \bigg( 
		(n +2  )^2   e^{\lambda f^{t,y} _{s} - \lambda n }
		\\ \notag
		&
		\times 
		\P \Big\{\X_k \leq f^{t,y} _k, k \leq s,
		\X_s ^{(\lambda)} + \frac 12 \W ^{(\lambda)}  _1  + \frac 12 \W ^{(\lambda)}  _2 \in (f^{t,y} _{s+1} -n - 1, f^{t,y} _{s+1} - n ] \Big\} \bigg)
		\\ \notag
		\leq \, &
		C \left( \frac{\Phi (\lambda)}{m} \right)^{s} \sum\limits _{n = 0} ^\infty 
		(n+2)^2 e^{- \lambda n} e^{\lambda f^{t,y} _{s}}
		\frac{(1+y)n }{(s+1)^{3/2}} 
		\\ 
		\label{in der Grafik}
		\leq \, &
		\frac{C (1+y)}{(s+1)^{3/2}}
		\left( \frac{\Phi (\lambda)}{m} \right)^{s} 
		e^{\lambda f^{t,y} _{s}} 
		\\ \notag
		= \, &  \frac{C (1+y)}{(s+1)^{3/2}}
		\left( \frac{\Phi (\lambda)}{m} \right)^{s}
		\left(\Phi (\lambda)\right) ^s
		e^y
		\frac{(s+1)^{(\d-1)/2}(t-s+1)^{3/2}}{(t+1)^{3/2}},
	\end{align}
	and consequently by
	Lemma \ref{hoodwink = deceive, trick}
	we get from \eqref{fete = prazdnik}
	\begin{align*}
		\E  [(\#  \mathcal{A}^{t,y}_\theta)^2]\leq & \,
		C (1+y) e^{-y} t^{-\frac{\d-1}{2}} +
		C (1+y) e^{-y}
		\sum\limits _{s=1} ^t 
		\frac{(s+1)^{(\d-1)/2}(t-s+1)^{3/2}}{(s+1)^{3/2}(t+1)^{3/2}(t-s+1)^{3}(t+1)^{{(\d-4)}}}
		\\
		\leq & \, C (1+y) e^{-y} t^{-\frac{\d-1}{2}} + \frac{C (1+y) e^{-y}}{(t+1)^{{(\d-1)/2}}}
		\sum\limits _{s=1} ^t 
		\left( \frac{s+1}{t+1} \right) ^{(\d-1)/2}
		\frac{ (t+1)^{3/2} }{(s+1)^{3/2}(t-s+1)^{3/2}}
		\\
		\leq & \, C (1+y) e^{-y} t^{-\frac{\d-1}{2}} + \frac{C (1+y) e^{-y}}{(t+1)^{{(\d-1)/2}}}
		\sum\limits _{s=1} ^t 
		\frac{ (t+1)^{3/2} }{s^{3/2}(t-s+1)^{3/2}}.
	\end{align*}
	It remains to note that the last series is 
	uniformly bounded in $t$
	by \eqref{riff-raff}.
\end{proof}

  \subsection{Lower bound on the first moment $\E  [\#  \mathcal{A}^{t,y}_\theta]$}  
The next lemma plays an auxiliary role.
It shows that the function 
$ f  $ dominates in a certain sense
two   functions, $g$ and $h$, with desirable properties.
These two functions are used
in the proof of Lemma \ref{hoodwink = deceive, trick2}.
The function $g$ does not deviate too much from the linear
map $s \mapsto \frac{\Psi(\lambda)}{\lambda}s$,
which allows to apply
Theorem \ref{thm beinhalten = include}
to $g$ and the random walk $\big(-\X ^{(\lambda)} _t + \frac{\Psi (\lambda)}{\lambda}t, t \in \Z_+\big)$
similarly to the transformations
in \eqref{limpid = clear, transparent}.
The function $h$ grows faster than the square root.
Below the operations of minimum
and maximum, $\wedge$ and $\vee$ respectively,
succeed multiplication in the operation order,
that is, for $a,b,c \in \R$, $ab \vee c = (ab) \vee c$.
\begin{lem}\label{colic}
	For $s,t \in \Z_+$, $y \in [1,t ^{1/2}]$, $s \leq t$ define
	$$g^{t,y}_s = f^{t,y}_s - 
	\varepsilon\big[s \wedge (t-s) \big] ^{0.4}
	\1 \Big\{s > \frac{ M^{5/3} }{\varepsilon ^{5/3} \lambda ^{5/3}} \Big\}
	- \frac{M}{\lambda}
	\1 \Big\{s \leq  \frac{ M^{5/3} }{\varepsilon ^{5/3} \lambda ^{5/3}} \Big\} - \frac12 \1 \{s=t\}.
	$$
	and 
	$$
	h^{t }_s = \varepsilon   [s^{0.6} \wedge t ^{1/2}] \vee \frac M \lambda.
	$$
	Provided that $\varepsilon > 0$
	is sufficiently small it holds that for  large $t \in \N$
	for all $s \in \{0,1,...,t\}$
	\begin{equation}\label{tattle on}
		(g^{t,y} _s)^2 + (h^{t} _s)^2 \leq (f^{t,y} _s)^2.
	\end{equation}
\end{lem}
  
\begin{proof} 
Recall that $f^{t,y} _s \geq 2\frac{M}{\lambda}$ by \eqref{der Platz}
and hence for some $c > 0$
\begin{equation*}
	f^{t,y} _s \geq c s \vee 2\frac{M}{\lambda}.
\end{equation*}
The inequality in \eqref{tattle on}
can be checked separately for $${s \leq \varepsilon ^{-5/3} \lambda ^{-5/3} M ^{5/3}},$$
and then under the assumption $s > \varepsilon ^{-5/3} \lambda ^{-5/3} M ^{5/3}$
in the cases
$s \leq t^{5/6}$, $t^{5/6} \leq s \leq \frac t2 $,
$\frac t2 \leq s \leq t - 1 $, $s = t$.
If $s \leq \varepsilon ^{-5/3} \lambda ^{-5/3} M ^{5/3} $
then 
\begin{multline*}
	(f^{t,y} _s)^2
	-
	(g^{t,y} _s)^2 -  (h^{t }_s)^2
	=
	(f^{t,y} _s)^2
	-
	\big(f^{t,y}_s - 
	\frac{M}{\lambda}\big)^2 -  \frac{M^2}{\lambda ^2} 
	= 2  f^{t,y} _s \frac{M}{\lambda} - 
	\frac{M^2}{\lambda ^2}-
	\frac{M^2}{\lambda ^2}
	> 0
\end{multline*} 
since $f^{t,y} _s \geq 2\frac{M}{\lambda}$ by \eqref{der Platz}.\\
  
In the rest of the proof  assume $s \geq \varepsilon ^{-5/3} \lambda ^{-5/3} M ^{5/3}$.
If $s \leq t^{5/6}$, then 
\begin{equation*}
	(f^{t,y} _s)^2
	-
	(g^{t,y} _s)^2 -  (h^{t }_s)^2
	=
	(f^{t,y} _s)^2
	-
	(f^{t,y}_s - 
	\varepsilon s  ^{0.4})^2 - \varepsilon ^2 s^{1.2} 
	=2 \varepsilon s  ^{0.4}  f^{t,y}_s - 
	\varepsilon ^2 s  ^{0.8}
	- \varepsilon ^2 s^{1.2}  \geq 0
\end{equation*}
provided $\varepsilon$ is sufficiently small.
If $t^{5/6} \leq s \leq \frac t2 $, then 
\begin{equation*}
	(f^{t,y} _s)^2
	-
	(g^{t,y} _s)^2 -  (h^{t }_s)^2
	=
	(f^{t,y} _s)^2
	-
	(f^{t,y}_s - 
	\varepsilon s  ^{0.4})^2 - \varepsilon ^2 t 
	= 2 \varepsilon s  ^{0.4}  f^{t,y}_s  - 
	\varepsilon ^2 s  ^{0.8}
	- \varepsilon ^2 t \geq 0
\end{equation*}
for a small $\varepsilon$.
If $\frac t2 \leq s \leq t-1 $, then 
\begin{multline*}
	(f^{t,y} _s)^2
	-
	(g^{t,y} _s)^2 -  (h^{t }_s)^2
	=
	(f^{t,y} _s)^2
	-
	(f^{t,y}_s - 
	\varepsilon (t-s)  ^{0.4})^2 - \varepsilon ^2 t 
	\\
	= 2 \varepsilon (t-s)  ^{0.4}  f^{t,y}_s - 
	\varepsilon ^2 (t-s)  ^{0.8}
	- \varepsilon ^2 t \geq 0
\end{multline*}
provided $\varepsilon$ is   small.
Finally, for $s = t$ we have 
\begin{equation*}
	(f^{t,y} _s)^2
	-
	(g^{t,y} _s)^2 -  (h^{t }_s)^2
	=
	(f^{t,y} _s)^2
	-
	(f^{t,y}_s - 
	1/2)^2 - \varepsilon ^2 t 
	=   f^{t,y}_t + 1/4
	- \varepsilon ^2 t \geq 0
\end{equation*}
provided $\varepsilon$ is sufficiently small.
\end{proof}

The next lemma provides a lower bound on the first moment.
\begin{lem}\label{hoodwink = deceive, trick2}
	There exists $C>1$ such that 
	for large $t \in \N$ and $y \in [1, t^{1/2}]$
	\begin{equation} \label{reticence2}
		\E \big[ \# \mathcal{A}^{t,y}_\theta  \big]
		\geq C^{-1} y e^{-y} t^{-\frac{\d-1}{2}}.
	\end{equation}
\end{lem}
  
\begin{proof}
Consider
the functions
$$g^{t,y}_s = f^{t,y}_s - 
\varepsilon\big[s \wedge (t-s) \big] ^{0.4}
\1 \{s > \varepsilon ^{-1} M^{5/3} \}
- \frac{M}{\lambda}
\1 \Big\{s \leq  \frac{ M^{5/3} }{\varepsilon ^{5/3} \lambda ^{5/3}} \Big\} - \frac12 \1 \{s=t\} 
$$
and 
$$
h^{t }_s = \varepsilon   [s^{0.6} \wedge t ^{1/2}] \vee \frac M \lambda 
$$
from Lemma \ref{colic}
with a small $\varepsilon > 0$,
so that \eqref{tattle on} is satisfied for $s,t \in \Z_+$, $y \in [1,t ^{1/2}]$, $s \leq t$.
For a
fixed  $\theta \in \S ^{\d-1} $ 
and
$u \in \T _t$ denote $\X ^u _s = \langle X_s(u), \theta \rangle$ the projection 
of the spatial positions of the ancestral line of $u$
on the line $\{ a\theta: a \in \R\} $, 
and let $\X^{u,\perp} _s  = X_s(u) - \X ^u _s \theta $
be the projection on the orthogonal complement of 
$\{ a\theta: a \in \R\} $. 
Note that for a small $\varepsilon > 0$, $g^{t,y} _s \geq 0$ and
by \eqref{tattle on} for large $t \in \N$
\begin{multline} 
	\{u \in \T _t: 
	\X^u _s \leq g^{t,y} _s, s \leq t, \X^u _t \geq f^{t,y} _t - 1 \}
	\cap \{u \in \T _t: |\X^{u,\perp} _s| \leq h^{t} _s   \}
	\\ \subset 
	\{u \in \T _t: 
	|X_s(u)|  \leq f^{t,y} _s, s \leq t, \X^u _t \geq f^{t,y} _t - 1 \}
	= \mathcal{A}^{t,y}_\theta
\end{multline}
since $g^{t,y} _t  = f^{t,y} _t - \frac 12$
and for every $s \in \{0,1,...,t\}$
\begin{equation}\label{in der Sonne sitzen}
	(f^{t,y} _s)^2
	-
	(g^{t,y} _s)^2 -  ( h^t _s  )^2
	\geq 0.
\end{equation}
Recall the random walks $\{\X_t, t \in \Z_+\}$
and $\{\X^\perp_t, t \in \Z_+\}$
introduced on Page \pageref{X def}.
By the many-to-one lemma (Lemma \ref{lem manyto1})
\begin{align}
	\E \big[ \# \mathcal{A}^{t,y}_\theta \big] \geq 
	\ &
	\E \big[ \# \{u \in \T _t: 
	\text{for } s \leq t, 
	\X^u _s \leq g^{t,y} _s,  \X^u _t \geq f^{t,y} _t - 1  \text{ and } |\X^{u,\perp} _s| \leq h^{t} _s   \} \big]
	\notag 
	\\ \notag
	=  \ &  m ^t
	\P  \Big\{ 
	\text{for } s \leq t, 
	\X _s \leq g^{t,y} _s,  \X _t \geq g^{t,y} _t - \frac 12  \text{ and } |\X^{\perp} _s| \leq h^{t} _s   \Big\}
	\\ \label{der Tresor = safe}
	=  \ &  m ^t
	\P (B) \P \big\{ |\X^{\perp} _s| \leq h^{t} _s, s \leq t \big| B \big\},
\end{align}
where 
\begin{equation}\label{lindern = alleviate, ease, mitigate}
B := \Big\{ 
\X _s \leq g^{t,y} _s, s \leq t,   \X _t \geq g^{t,y} _t - \frac 12 \Big\}.	 
\end{equation}
Similarly to \eqref{lindern = alleviate, ease, mitigate} let us introduce
	\begin{equation}\label{lindern lambda}
		B^{(\lambda)} := \Big\{ 
		\X ^{(\lambda)} _s \leq g^{t,y} _s, s \leq t,   \X ^{(\lambda)} _t \geq g^{t,y} _t - \frac 12 \Big\}.	 
\end{equation}

For $b \in \R$ with $\Phi(b) < \infty$
we have for $\text{Law}(\X _t)$-almost all $a \in \R $
the equality of conditional distributions
\begin{equation}\label{jdm etw um die Ohren hauen}
	\text{Law}\big( \{\X _s - \X _{s-1}, 1\leq s\leq t \} \big| \X _t = a \big) 
	=
	\text{Law}\big( \{\X ^{(b)} _s - \X ^{(b)} _{s-1},  1\leq s\leq t \} \big| \X ^{(b)} _t = a \big)
\end{equation}
and hence for any Borel $\mathfrak{C}\subset \R ^t $
	the events 
	$$
\mathcal{C} = 
	\Big\{  \{\X _s - \X _{s-1}, 1\leq s\leq t \} \in \mathfrak{C} \Big\} 
	 \ \ \ \text{and} \ \ \
	 \mathcal{C} ^{(\lambda)} = 
	\Big\{ \{\X ^{(\lambda)} _s - \X ^{(\lambda)} _{s-1}, 1\leq s\leq t \} \in \mathfrak{C}
	\Big\} 
	$$
	satisfy
	\begin{equation}\label{erbärmlich = pathetic, abject}
		e^{ g^{t,y} _t - 1/2  }	\leq \frac{\P (\mathcal{C} ^{(\lambda)} \cap B^{(\lambda)})}
		{\P (\mathcal{C} \cap B  )} \leq e ^ { g^{t,y} _t }.
	\end{equation}
Conditioning 
on $\X ^{(\lambda)} _s - \X ^{(\lambda)} _{s-1}$
we get   
by Lemma \ref{deign}
\label{pageref Lemma deign}
for   $s$ with $1 \leq s \leq t$
\begin{equation}\label{die Flosse = fin}
	 \P \big( B^{(\lambda)} \big|  \X ^{(\lambda)} _s - \X ^{(\lambda)} _{s-1} \big)
	 \leq C \Big(  |\X ^{(\lambda)} _s - \X ^{(\lambda)} _{s-1}|^4+1 \Big)
	 \P (B^{(\lambda)} ).
\end{equation}
Consider the map
\begin{equation}
  \R \ni a \mapsto  h(a) = \E \Big[ e^{\frac \varepsilon 2 |Q_1|}  \Big| \X _1 = a \Big]
  \in [1, \infty)
\end{equation}
defined for $\text{Law}(\X _1)$-almost all $a\in \R$;
note that $ \E \big[ h
	\big(\langle Q_1 , \theta \rangle\big)
	 \big] = \E \big[ e^{\frac \varepsilon 2 |Q_1|}  \big] $ for $\theta \in \S ^\d$.  
We have by \eqref{die Flosse = fin}
\begin{align}
	  \E \big[ h( \X ^{(\lambda)} _s - \X ^{(\lambda)} _{s-1}) \1 _{B^{(\lambda)}}\big]
	 & = \E\bigg[ h( \X ^{(\lambda)} _s - \X ^{(\lambda)} _{s-1})  \P \big[ B^{(\lambda)} 
	 \big|   \X ^{(\lambda)} _s - \X ^{(\lambda)} _{s-1} \Big]
	 \bigg]
	 \\
	& \leq 
	C \E\bigg[ h( \X ^{(\lambda)} _s - \X ^{(\lambda)} _{s-1}) 
	 \Big(  |\X ^{(\lambda)} _s - \X ^{(\lambda)} _{s-1}|^4+1 \Big)
	 \P (B^{(\lambda)} )
	 \bigg]
	 \\
	& = 
	\label{ledrig = leathery}
	C \P (B^{(\lambda)} )\E \Big[ 
	h( \X ^{(\lambda)} _s - \X ^{(\lambda)} _{s-1}) 
	 \Big(  |\X ^{(\lambda)} _s - \X ^{(\lambda)} _{s-1}|^4+1 \Big)
	\Big].
\end{align}
Next we obtain by Assumption \ref{lambda-existence} for some $\overline C  > 0$
\begin{align}
	\E \Big[ 
	h( \X ^{(\lambda)} _s - \X ^{(\lambda)} _{s-1}) 
	\Big( & |\X ^{(\lambda)} _s - \X ^{(\lambda)} _{s-1}|^4+1 \Big)
	\Big]
	\\
	& \leq 
	C	\E \Big[ 
	h(\X  _s - \X  _{s-1})
	\Big(  |\X  _s - \X  _{s-1}|^4+1 \Big) e^{\lambda (\X  _s - \X  _{s-1}) }
	\Big]
		\\
	& = 
		C	\E \Big[ 
	 e^{\frac \varepsilon 2 |Q_s  - Q_{s-1} |}
	\Big(  |\X  _s - \X  _{s-1}|^4+1 \Big) e^{\lambda (\X  _s - \X  _{s-1}) }
	\Big]
	 \overset{ \eqref{crosshair}}{\leq} \overline C.
\end{align}
Given $\X  _s - \X  _{s-1}$,
$B$ and $|Q_s  - Q_{s-1} |$ 
are conditionally independent.
Hence by \eqref{ledrig = leathery}
\begin{equation}
	\E \big[ e^{\frac \varepsilon 2 |Q_s ^{(\lambda)} - Q_{s-1} ^{(\lambda)}|} \big| B^{(\lambda)} \big]
	=
	\frac{ \E \big[ h( \X ^{(\lambda)} _s - \X ^{(\lambda)} _{s-1}) \1 _{B^{(\lambda)}}\big]}{\P (B^{(\lambda)} )}
	 \leq \overline C.
\end{equation}
Therefore by \eqref{erbärmlich = pathetic, abject} 
	for all $t, s \in \N$ with   $1 \leq s \leq t$
	for some $C>0$
	\begin{equation}
		\E \big[ e^{\frac \varepsilon 2 |Q_s - Q_{s-1} |}  \big| B \big] \leq C
	\end{equation}
	and hence also 
	\[
	\E \bigg[ \exp\Big\{ \frac{\varepsilon}{2} \big|\X^{\perp} _{s} -\X^{\perp} _{s-1} \big| \Big\} \bigg|  B \bigg] \leq C <  \infty,
	\ \ \ s= 1,...,t.
	\]
A martingale invariance principle 
(e.g. \cite[Chapter 4]{HH80}) implies 
\begin{equation}\label{erroeten = rot werden}
	\inf\limits _{t \in \N}\P\Big\{ |\X^{\perp} _s| \leq h^{t} _s,  s\leq t \Big|  B \Big\} > 0. 
\end{equation}
  
This gives a lower bound on 
the conditional probability  in \eqref{der Tresor = safe}.
To obtain a  lower bound on $ \P (B) =   \P  \Big\{ 
\X _s \leq g^{t,y} _s, s \leq t,   \X _t \geq g^{t,y} _t - \frac 12 \Big\} $
we
apply the transformation \eqref{Verwahrung = custody, safekeeping2}
with $b = \lambda $,
\begin{align}
	\P (B)
	& = \E\big[ e^{\lambda \X _1}  \big]^t\E \Bigg[\frac{\1\Big\{ 
		\X ^{(\lambda)} _s \leq g^{t,y} _s, s \leq t,   \X ^{(\lambda)} _t \geq g^{t,y} _t - \frac 12 \Big\}}{e^{\lambda \X ^{(\lambda)} _t}}  \Bigg]
	\\ \label{chasm}
	& \geq C \left( \frac{\Phi (\lambda)}{m} \right)^{t}
	e^{-\lambda g^{t,y} _t}
	\P \Big\{ 
	\X ^{(\lambda)} _s \leq g^{t,y} _s, s \leq t,   \X ^{(\lambda)} _t \geq g^{t,y} _t - \frac 12 \Big\}.
\end{align}
Next
we proceed similarly to \eqref{limpid = clear, transparent}.
Using the random walk $Y_s = -\X ^{(\lambda)} _s + \frac{\Psi (\lambda)}{\lambda}s$ 
and 
$$G ^t _s =   \frac{\Psi (\lambda)}{\lambda}s - g_s^{t,0}
= \frac{\Psi (\lambda)}{\lambda}s - g_s^{t,y} + \frac y \lambda $$
we obtain by Theorem \ref{thm beinhalten = include}
for some $c > 0$
\begin{align} 
	\notag
	\P \Big\{ 
	\X ^{(\lambda)} _s \leq g^{t,y} _s, s \leq t,   \X ^{(\lambda)} _t \geq g^{t,y} _t - \frac 12 \Big\}
	& =  \P \Big\{ 
	Y_s \geq  G^{t} _s - \frac y \lambda, s \leq t,   Y _t \leq  G^{t} _t - \frac y \lambda + \frac 12 \Big\}
	\\
	& =
	\P \Big\{ 
	2Y_s \geq  2G^{t} _s - 2\frac y \lambda, s \leq t,   Y _t \leq  2G^{t} _t - 2\frac y \lambda + 1 \Big\}
	\\
	& \geq  \frac{c(1+y)}{(t+1)^{3/2}}
\end{align}
and hence by \eqref{chasm}
\begin{equation}\label{Struemer = forward}
	\P (B) \geq c  (y+1) m^{-t} e^{-y}  t ^{\frac{ 1 -\d }{2 } },
\end{equation}
where we took into account that
$ g^{t,y} _t =  f^{t,y} _t - \frac 12$ and 
$$
e^{-\lambda f^{t,y} _t} \geq C \big(\Phi (\lambda) \big)^{-t}
(t+1)^{-\d/2 + 2} e^{-y}.
$$

Finally,
combining \eqref{der Tresor = safe}, \eqref{erroeten = rot werden},
and \eqref{Struemer = forward}, we get 
for large $t$
\begin{equation}
	\E \big[ \# \mathcal{A}^{t,y}_\theta \big] \geq 
	c  (y+1)  e^{-y}  t ^{\frac{ 1 -\d }{2 } }.
\end{equation}
\end{proof}

  \subsection{Upper bound on $\E  [\#  \mathcal{A}^{t,y}_\theta \# \mathcal{A}^{t,y}_{\theta'}]$ for $\theta \ne \theta'$}
  \label{subsec challenging}
In this subsection   we encounter arguably the most technical part of the paper.
Recall our order convention:
$ab \wedge c = (ab) \wedge c$ for $a,b,c \in \R$.
\begin{lem} \label{rationell = efficient}
	For $a, b  \geq 0$ and $\alpha \in \left[ 0, \frac \pi 2 \right]$
	\begin{equation}
		\inf\limits _{x \in [0,1]} \Big[ a\Big( 1 - \cos \big(\alpha -  \arcsin ( x  ) \big) \Big) + b x \Big] \geq 
		\pi ^{-1} \big(a \alpha ^2 \wedge b \alpha \big).
	\end{equation}
\end{lem}
  
\begin{proof}
For $\beta \in \left[ - \frac \pi 2, \frac \pi 2 \right]$
the inequality $1 - \cos \beta \geq \frac{\beta ^2}{\pi}$ holds. We have
\begin{equation*}
	1 - \cos \big(\alpha -  \arcsin ( x  ) \big)
	\geq \frac{(\alpha -  \arcsin ( x  ) )^2}{\pi}
\end{equation*}
and hence 
\begin{align}
	\label{beeintraechtigen = affect, impair}
	\inf\limits _{x \in [0,1]} \Big[ a\Big( 1 - \cos \big(\alpha -  \arcsin ( x  ) \big) \Big) + b x \Big]
	\geq & \,
	\inf\limits _{x \in [0,1]} \Big[ a \frac{(\alpha -  \arcsin ( x  ) )^2}{\pi} + b x \Big] 
	\\  \notag 
	= & \, \inf\limits _{x \in [0,\sin \alpha  ]} \Big[ a \frac{(\alpha -  \arcsin ( x  ) )^2}{\pi} + b x \Big] 
	\\ \notag
	\geq  & \, \inf\limits _{x \in [0,\frac 2 \pi \alpha  ]} \Big[ a \frac{(\alpha -  \arcsin ( x  ) )^2}{\pi} + b x \Big] 
	\\ \notag
	\geq  & \, \inf\limits _{x \in [0,\frac 2 \pi \alpha  ]} \Big[ a \frac{(\alpha -  \frac \pi 2 x )^2}{\pi} + b x \Big]
	\\ \notag
	=  & \, \pi ^{-1} \inf\limits _{y \in [0,  \alpha  ]} \Big[ a {(\alpha -  y )^2} + 2b  y \Big],
\end{align}
where in the last two inequalities  we used that 
$\arcsin ( x  ) \leq \frac \pi 2 x $ 
for $x \in \big[0, 1 \big]$
and $\frac 2  \pi \alpha \leq  \sin \alpha  $ for
$\alpha \in \big[0, \frac \pi 2 \big]$.
The analysis of the parabola $y \mapsto a {(\alpha -  y )^2} + 2b  y$
shows that in the case $\alpha > \frac ba$
the minimum is achieved at $\alpha - \frac ba$ and
\begin{equation*}
	\inf\limits _{y \in [0,  \alpha  ]} \Big[ a {(\alpha -  y )^2} + 2b  y \Big] = \frac{b^2}{a} + 2 b \alpha - 2\frac{b^2}{a}   = b\Big( 2\alpha - \frac ba \Big) \geq b \alpha,
\end{equation*}
whereas in the case $\alpha \leq \frac ba$
the minimum is achieved at $0$: 
\begin{equation*}
	\inf\limits _{y \in [0,  \alpha  ]} \Big[ a {(\alpha -  y )^2} + 2b  y \Big] = a \alpha ^2.
\end{equation*}
Combined with \eqref{beeintraechtigen = affect, impair}
this gives the statement of the lemma. 
\end{proof}

Let us introduce the rate functions related to the walks $\{\X_t, t \in \Z_+\}$
and $\{Q_t, t \in \Z_+\}$,
\begin{equation}\label{I introduced}
	I(x)
	=
	\sup\limits _{y \in \R ^\d} \big[ \langle y, x \rangle  - \ln \E \big[e^{\langle y, Q_1 \rangle}\big] \big],
	\ \ \ x \in \R ^\d,
\end{equation}
and
$I_1$ is defined on $\R _+$ by
$$
I_1(b)
=
\sup\limits _{a \in \R } \big[ ab - \Psi(a) + \ln m \big]
= \sup\limits _{a \in \R } \big[ ab - \ln \E\big[ e^{a \X_1 }\big] \big].
$$
Note that $I_1(b) = I(b \theta)$
for a direction $\theta  \in \S ^{\d-1}$
and $b \geq 0$ and $I(x) = I_1(|x|)$
for $x \in \R ^\d$.
We claim that the rate function $I_1$ is convex and satisfies 
\begin{equation}\label{I1prop}
	I_1\Big(\frac{\Psi(\lambda)}{\lambda}\Big)  = \ln m, \ \ \ 
	I_1'\Big(\frac{\Psi(\lambda)}{\lambda}\Big) = \lambda.
\end{equation}
To see that \eqref{I1prop} holds true, note that the function $\R_+ \ni q \mapsto \sup\limits _{a > 0} \big[a q - \Psi(a) \big]$ is the Legendre transform
of the function $\Psi$. Its derivative
is given by the inverse of the derivative $\frac{\Phi'}{\Phi}$ of $\Psi = \ln \Phi$.
Consequently the derivative at $q = \frac{\Psi(\lambda)}{\lambda}$,
$
\kappa := \frac{d
}{dq}\sup\limits _{a > 0} \big[a q - \Psi(a) \big]
\bigg\vert _{q = \frac{\Psi(\lambda)}{\lambda}}
$, satisfies
$\frac{\Phi '(\kappa) }{\Phi (\kappa)} = \frac{\Psi(\lambda)}{\lambda}$.
Recall that $\lambda > 0$ is a unique solution to \eqref{saggy ne soggy}; in particular,
	$$
	\frac{\Phi '(\lambda) }{\Phi (\lambda)} = \Psi'(\lambda) =  \frac{\Psi(\lambda)}{\lambda}.
	$$
	Since the map $\R_+ \ni \kappa \mapsto \frac{\Phi '(\kappa) }{\Phi (\kappa)} = \Psi'(\kappa)$ is increasing
	we have $\kappa = \lambda$,
that is, $I_1'\Big( \frac{\Psi(\lambda)}{\lambda} \Big) = \lambda$. \label{I of lambda}
Furthermore $\sup\limits _{a > 0} \left[a\frac{\Psi(\lambda)}{\lambda}
- \Psi(a) \right] = 0$ because the supremum is 
achieved for $a = \lambda$, and hence $I_1\Big( \frac{\Psi(\lambda)}{\lambda} \Big) =  \ln m$.

For a radius $r> 0$ we define
$$\mcb _r = \{ (x, y_1, y_2) \in (\R^d)^3: |x| \leq r, |x+y_i| \leq r,
i = 1,2 \}.$$
We introduce more notation related to $(x, y_1, y_2) \in \mcb _r$.
Denote $ \mathfrak{m} = |x + y _1| \vee |x + y _2| $ and $2 \mathfrak{d} = |y_1 - y_2|  $,
and denote by $\gamma$
the angle between $\mathbf{e}_1$
and $x + y_1 $, $- \pi < \gamma \leq \pi$. 
Recall that $\varepsilon > 0$ is introduced in Assumption \ref{lambda-existence}.

\begin{lem}\label{in Gang kommen}
	Let $r = r(s) = \frac{\Psi (\lambda)}{\lambda}s + o (s)$ 
	as $s \to \infty$.	
	For $\alpha \in [0, \frac{\pi}{2}]$
	let
	$\theta _1 = (\cos \alpha, \sin \alpha)$
	and
	$\theta _2 = (\cos \alpha, -\sin \alpha)$. 
	Define a function $h =  h_{\alpha}$ from $ \mcb _r$ to $ \R$ by
	\begin{equation*}
		h_{\alpha}(x, y_1, y_2) = \lambda \langle x + y_1, \theta _1 \rangle
		+ \lambda \langle x + y_2, \theta _2 \rangle
		- I\Big(\frac xs \Big) s - (\lambda +\varepsilon) (|y_1|\vee |y_2|).
	\end{equation*}
	Then for large $s$
	\begin{itemize}
		\item[$(i)$]  
		$h_{\text{max}} = h_{\text{max}}(r) := \sup\limits _{\alpha \in [0, \frac{\pi}{2}]}
		\sup \limits _{(x, y_1, y_2) \in \mcb _r}  h_{\alpha}(x, y_1, y_2) = h_0( r \mathbf{e}_1, 0_2, 0_2 ) $,
		where $\mathbf{e}_1 = (1,0)$, $0_2 = (0,0)$.
		\item[$(ii)$] 
		If $\alpha \leq \frac \pi 6$
		then 
		for $(x, y_1, y_2) \in \mcb _r $ 
		with $\gamma \in [0,  \frac \pi 4]$
		\begin{equation*}
			h_{\text{max}} -  h_{\alpha}(x, y_1, y_2)
			\geq    \Big( \lambda -  \frac{\varepsilon}{10}  \Big)
			(r - \mathfrak{m} )  
			+
			2 \lambda \mathfrak{m} \Big[ 1 - \cos \Big(\alpha -  \arcsin \big( \frac{\mathfrak{d}}{ \mathfrak{m}}  \big) \Big)   \Big]
			+ 
			0.99 \varepsilon \mathfrak{d}.
		\end{equation*}
		\item[$(iii)$] For some $c> 0$ 
		for
		$(x, y_1, y_2) \in \mcb _r $
		\begin{equation}\label{sich verplappern = spill the beans}
			h_{\text{max}} -  h_{\alpha}(x, y_1, y_2)
			\geq    \big( \lambda -  \frac{\varepsilon}{10}  \big)
			(r - \mathfrak{m} ) 
			+ c  \alpha ^2.
		\end{equation}
	\end{itemize}
\end{lem}
\begin{rmk}
	Figure \ref{f_for_proof} is an illustration to Lemma \ref{in Gang kommen}.	
	The `positive' part of $h$, $\lambda \langle x + y_1, \theta _1 \rangle
	+ \lambda \langle x + y_2, \theta _2 \rangle$,
	can be seen 
	as the gain, 
	whereas the `negative' part, $I\big(\frac xs \big) s + (\lambda +\varepsilon) (|y_1|\vee |y_2|)$,
	can be thought of as the cost. 
	A higher-dimensional version of the function $h$ will appear
	below
	in the exponent of a certain expression 
	related to the expectation $\E  [\#  \mathcal{A}^{t,y}_\theta \# \mathcal{A}^{t,y}_{\theta'}]$, $\theta \ne \theta '$, $\theta, \theta ' \in \S ^{\d-1}$.
\end{rmk}

\begin{proof}[Proof of Lemma \ref{in Gang kommen}]

Without loss of generality we assume that $\varepsilon \leq \lambda$.
Recall that 
the function $I_1 : \R_+ \to \R _+$
introduced on Page \pageref{I introduced}
is an increasing convex
function which satisfies \eqref{I1prop}.
Therefore
for large $s$ for all $x \leq r $
\begin{equation}\label{pusten = puff}
	I_1\Big(\frac xs \Big) \leq \ln m + \frac{\varepsilon}{100}, \ \ \ 
	I_1'\Big(\frac xs\Big) \leq \lambda +  \frac{\varepsilon}{100}.
\end{equation}
and
for   $0\leq a \leq \frac{\Psi(\lambda)}{\lambda} $
\begin{equation}\label{Errungenschaft = achievement}
	I_1 (a  ) \leq \ln m, \ \ \ 
	I_1' (a  ) \leq \lambda.
\end{equation}
Note that $|y_1| \vee |y_2| \geq \mathfrak{d}$
and hence by \eqref{pusten = puff} 
\begin{multline}\label{das Loch}
	I\Big(\frac xs \Big) s + (\lambda +\varepsilon) (|y_1|\vee |y_2|) \geq 
	I_1\Big(\frac {\mathfrak{m} - (|y_1|\vee |y_2|)}{s} \Big) s + (\lambda +\varepsilon) (|y_1|\vee |y_2|)
	\\
	\geq 
	I_1\Big(\frac{\mathfrak{m} - \mathfrak{d}}{s} \Big)s + (\lambda +\varepsilon) \mathfrak{d}
	\geq I_1\Big(\frac{\mathfrak{m}}{s} \Big) s + 
	0.99 \varepsilon \mathfrak{d}.
\end{multline}

\begin{figure}[t!]
	\vspace{-1.5cm}
	\centering
	\includegraphics[scale=0.51]{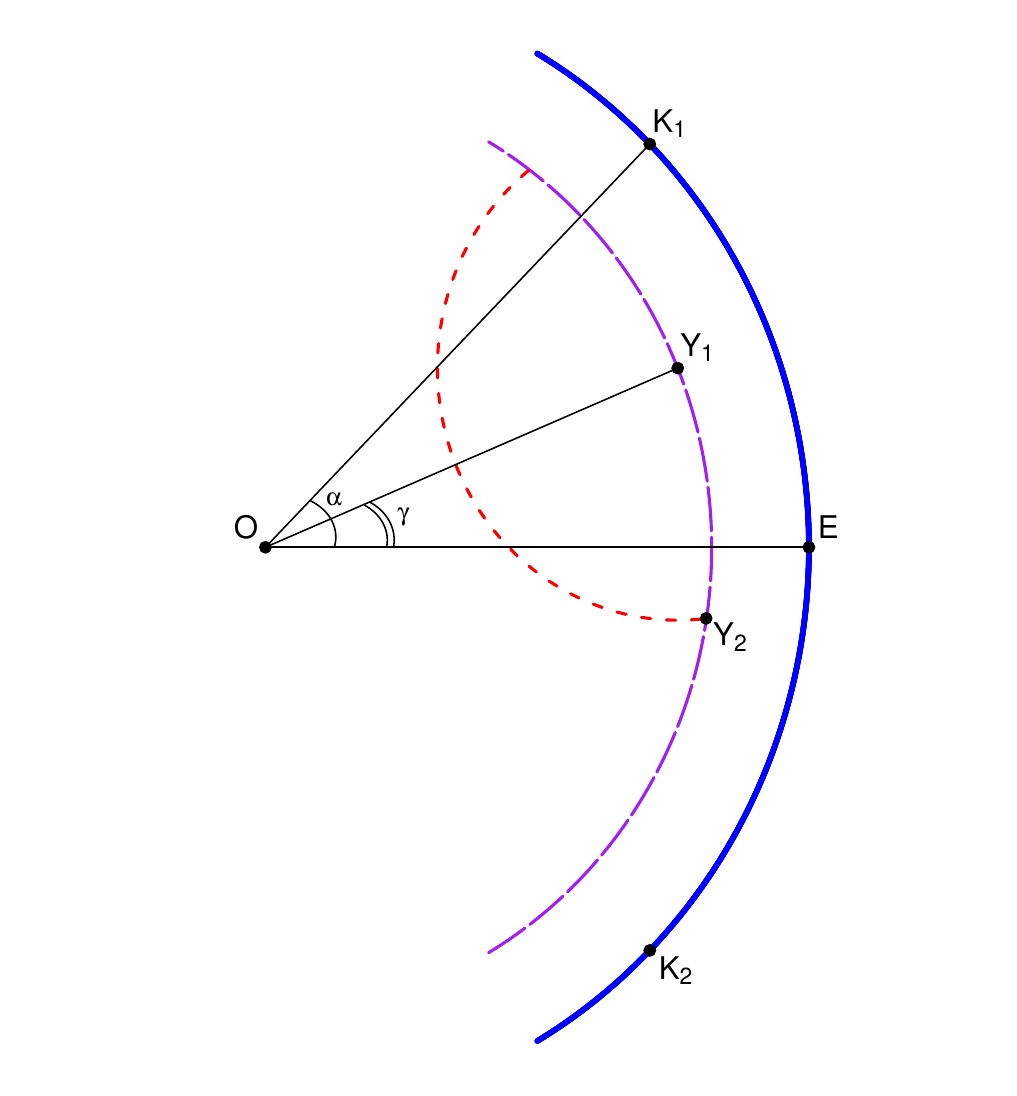}
	\caption{ \small An illustration to the proof of Lemma \ref{in Gang kommen}. The point coordinates are ${\rm{K}}_1 = (r \cos \alpha, r \sin \alpha )$, ${\rm{K}}_2 = (r \cos \alpha, - r \sin \alpha )$, ${\rm{E}} =(r, 0)$,
		${\rm{O}} = (0,0)$,
		and ${\rm{Y}}_1 = (\mathfrak{m} \cos \gamma, \mathfrak{m} \sin \gamma) = x + y_1$. The thick solid blue  arc is a part of the circle of radius $r$ around ${\rm{O}}$, the long-dashed purple  arc belongs to the circle
		of radius $\mathfrak{m}$ around ${\rm{O}}$. The dashed  red arc contains points $z_2 $  satisfying $|z_2 - x - y_1| = 2\mathfrak{d}$
		and $|z_2| \leq \mathfrak{m} $. We see that the maximum 
		of $\langle z_2, \theta _2 \rangle $ on the dashed red arc is achieved 
		for $z_2 = {\rm{Y}}_2$ (at least whenever $\alpha \leq \frac \pi 6$ and $\gamma \in [0, \frac \pi 4]$). Since the distance between ${\rm{Y}}_1$ and ${\rm{Y}}_2$ is $|y_2 - y_1| = 2\mathfrak{d}$,  the angle between ${\rm{OY}}_1$ and ${\rm{OY}}_2$ is $\angle {\rm{Y}}_1{\rm{O}}{\rm{Y}}_2 = 2 \arcsin \big( \frac{\mathfrak{d}}{ \mathfrak{m}}  \big)$, 
		and $ \angle {\rm{Y}}_2 {\rm{O}} {\rm{K}}_2 = \alpha + \gamma - 2 \arcsin \big( \frac{\mathfrak{d}}{ \mathfrak{m}}  \big) $.
	}
	\label{f_for_proof}
	\vspace{1cm}
\end{figure}
  
Therefore  for large $s$
\begin{multline*}
	h_0( r \mathbf{e}_1, 0_2, 0_2 ) -  h_{\alpha}(x, y_1, y_2) \geq 
	2 \lambda r - I_1\Big(\frac rs\Big)s
	- 2 \lambda \mathfrak{m} - I_1\Big(\frac{\mathfrak{m}}{s} \Big) s
	\\ \geq 2 \lambda (r - \mathfrak{m} ) - \Big( \lambda +  \frac{\varepsilon}{100} \Big) (r - \mathfrak{m} ) \geq 0.
\end{multline*}
and $(i)$ is proven.\\
  
To show $(ii)$ recall that $\gamma $
the angle between $x + y_1 $ and $\mathbf{e}_1$,
and note that 
\[
\langle x + y_1, \theta_1 \rangle = |x + y_1| \cos (\alpha - \gamma).
\]
For $x, y_1, y_2$ with
$ \mathfrak{m} = |x + y _1| \vee |x + y _2| $ and $ 2 \mathfrak{d} = |y_1 - y_2|  $ we assume without loss of generality
that $ \mathfrak{m} = |x + y _1|$.
The locus of points $z_2 \in \R ^2$
satisfying $  |z _2| \leq \mathfrak{m}  $ and $ |x+ y_1 - z_2| = 2 \mathfrak{d}  $ 
is given by the dashed red arc on Figure \ref{f_for_proof}.
In particular we have
\begin{align}
	\langle x + y_1, \theta _1 \rangle
	+  \langle x + y_2, \theta _2 \rangle \leq
	& \,  \mathfrak{m} \Big[\cos (\alpha - \gamma) + 
	\cos \Big(\alpha - 2 \arcsin \Big( \frac{\mathfrak{d}}{ \mathfrak{m}}  \Big) + \gamma \Big) \Big]
	\notag 
	\\ \notag 
	= & \, 2 \mathfrak{m} \cos \Big(\alpha -  \arcsin \Big( \frac{\mathfrak{d}}{ \mathfrak{m}}  \Big) \Big)
	\cos \Big(  \arcsin \Big( \frac{\mathfrak{d}}{ \mathfrak{m}}  \Big) - \gamma \Big)
	\\ \label{der Trost = comfort}
	\leq & \,
	2 \mathfrak{m} \cos \Big(\alpha -  \arcsin \Big( \frac{\mathfrak{d}}{ \mathfrak{m}}  \Big) \Big),
\end{align}
with the last inequality turning into equality 
for $\gamma =  \arcsin \Big( \frac{\mathfrak{d}}{ \mathfrak{m}}  \Big)$.
By
\eqref{pusten = puff}, 
\eqref{das Loch}, and
\eqref{der Trost = comfort}
\begin{align*}
	h_{\text{max}} -  h_{\alpha}(x, y_1, y_2) 
	\geq & \, 
	2\lambda r - 
	2 \lambda \mathfrak{m} \cos \Big(\alpha -  \arcsin \Big( \frac{\mathfrak{d}}{ \mathfrak{m}}  \Big) \Big)
	-I_1 \Big(\frac rs\Big) s +  I_1 \Big(\frac{\mathfrak{m}}{s} \Big) s + 0.99 \varepsilon \mathfrak{d}
	\\ 
	\geq & \, 
	2\lambda r - 
	2 \lambda \mathfrak{m} \cos \Big(\alpha -  \arcsin \Big( \frac{\mathfrak{d}}{ \mathfrak{m}}  \Big) \Big)
	- \Big( \lambda +  \frac{\varepsilon}{100}  \Big)
	(r - \mathfrak{m} ) + 
	0.99 \varepsilon \mathfrak{d}
	\\ 
	= & \, 
	\Big( \lambda -  \frac{\varepsilon}{100}  \Big)
	(r - \mathfrak{m} )  
	+
	2 \lambda \mathfrak{m} \Big[ 1 - \cos \Big(\alpha -  \arcsin \Big( \frac{\mathfrak{d}}{ \mathfrak{m}}  \Big) \Big)   \Big]
	+ 
	0.99 \varepsilon \mathfrak{d},
\end{align*}
and $(ii)$ is proven.
Applying Lemma \ref{rationell = efficient} with $a = 2 \lambda \mathfrak{m}$,
$x = \frac{\mathfrak{d}}{\mathfrak{m}}$ (note that $\mathfrak{m} \geq \mathfrak{d}$), and $b = 0.99 \varepsilon \mathfrak{m}$ to the cases covered by $(ii)$ we find 
\begin{multline}
	h_{\text{max}} -  h_{\alpha}(x, y_1, y_2) 
	\geq 
	\Big( \lambda -  \frac{\varepsilon}{100}  \Big)
	(r - \mathfrak{m} ) 
	+ \pi ^{-1} \Big[ 2 \lambda \mathfrak{m} \alpha ^2
	\wedge 0.99 \varepsilon \mathfrak{m} \alpha    \Big]
	\\
	\geq 
	\Big( \lambda -  \frac{\varepsilon}{100}  \Big)
	(r - \mathfrak{m} ) 
	+ 0.99 \varepsilon  \mathfrak{m} \pi ^{-1}  [   \alpha ^2
	\wedge    \alpha   ].
\end{multline}
Since  $\alpha \geq \frac 2 \pi \alpha  ^2$
for $\alpha \in [0, \frac \pi 2]$, \eqref{sich verplappern = spill the beans} follows. The remaining case
when either $\alpha \geq \frac \pi 6$
or $\gamma \geq \frac \pi 4$, 
can be analyzed by a simpler version of the same arguments.
For instance if $\alpha \geq \frac \pi 6$, then 
in the case $2\mathfrak{d} \leq \mathfrak{m} \sin \mkern-2mu \frac {\pi}{12}$
either the (absolute value of the) angle between $x + y_1$  and $\theta _1$
or the angle between  $x + y_2$ and $\theta _2$ is going to be 
at least $\frac {\pi}{12}$, and hence  
$$
\langle x + y_1, \theta _1 \rangle
+ \langle x + y_2, \theta _2 \rangle
\leq \mathfrak{m} + \mathfrak{m} \cos \frac {\pi}{12}
= \Big( 1 + \cos \frac {\pi}{12} \Big) \mathfrak{m}  
$$
gives a bound on the gain;
in the case $2\mathfrak{d} \geq \mathfrak{m} \sin \frac {\pi}{12}$
the bound on the cost \eqref{das Loch} gives
the desired inequality.
Finally, the case $\gamma \leq 0$ then follows by symmetry 
considerations. 
\end{proof}

In the next lemma we establish an inequality
later
used to get an upper bound on
$\E  [\#  \mathcal{A}^{t,y}_{\theta_1 } \# \mathcal{A}^{t,y}_{\theta_2}]$  for $\theta _1 \ne \theta _2$, $\theta _1, \theta _2 \in \S ^{\d-1}$.

\begin{lem}\label{pauper = poor person}
	Let 
	$\theta _1, \theta _2 \in \S ^{\d - 1}$ 
	and let $2 \alpha \in [0, \pi]$
	be the angle between $\theta _1$ and $\theta_2$, 
	$\cos (2\alpha) = \langle \theta_1, \theta _2 \rangle$.	
	Set
	\begin{align*}
		\mce = \mce_\alpha =  \, & \E \bigg[
		\big(1 + r   - \langle Q_s + \Delta _1, \theta _1 \rangle \big)
		\big(1 + r  - \langle Q_s + \Delta _2, \theta _2 \rangle \big) e^{\lambda \langle Q_s+ \Delta _1, \theta _1 \rangle
			+	 \lambda \langle Q_s+ \Delta _2, \theta _2 \rangle}
		\\
		& \times
		\1 \{|Q_k| \leq f^{t,y} _k, k = 0,...,s, |Q_s+\Delta _1|\vee |Q_s+\Delta _2| \leq f^{t,y} _{s+1}  \}    
		\bigg].
	\end{align*}
	Then for some $C, c> 0$
	\begin{equation}\label{pull someone up for something}
		\mce \leq  \frac{C e^{2 \lambda f^{t,y} _{s} -c \alpha ^2s -y }}{m^s} \times  \frac{  (t+1) ^{3/2} (1+y)}{(s+1)^{(\d + 2)/2} (t-s+1) ^{3/2}} .
	\end{equation}
\end{lem}
\begin{rmk}
	The proof of Lemma \ref{pauper = poor person} is quite lengthy. Let us point out that assuming the finiteness of the $(2 \lambda + \varepsilon)$-th moment:
	$$\E \Big[ \L (\R ^\d)  \int\limits_{\R^\d}
	e^{ (2\lambda + \varepsilon) \langle x, \theta \rangle} \L  (dx)\Big] < \infty$$ 
	would significantly simplify the proof. Indeed under this assumption $\E \big[ e^{\lambda \langle  \Delta _1, \theta _1 \rangle
		+	 \lambda  \langle  \Delta _2, \theta _2 \rangle}\big] < \infty $ implying that 
	\begin{align*}
		\mce \leq C   \, & \E \bigg[
		\big(1 + r   - \langle Q_s + \Delta _1, \theta _1 \rangle \big)
		\big(1 + r  - \langle Q_s + \Delta _2, \theta _2 \rangle \big) e^{\lambda \langle Q_s, \theta _1 \rangle
			+	 \lambda \langle Q_s, \theta _2 \rangle}
		\\
		& \times
		\1 \{|Q_k| \leq f^{t,y} _k, k = 0,...,s \}    
		\bigg].
	\end{align*}
	and the following analysis can be based on the values of $\lambda \langle Q_s, \theta _1 + \theta _2 \rangle	$. That is, considering  a projection on the line defined by $\frac{\theta _1 + \theta _2}{|\theta _1 + \theta _2|}$ would be sufficient. In the actual proof we work with a two-dimensional projection containing $\theta _1$  and  $\theta _2$.		
\end{rmk}
  
\begin{proof}[Proof of Lemma \ref{pauper = poor person}]
To simplify notation in the analysis below we start by observing that we can
treat a slightly simpler expression by
replacing $f^{t,y} _{s+1}$
by $f^{t,y} _s$: for some $C> 0$
\begin{align*}
	\mce \leq C   \, & \E \bigg[
	\big(1 + r   - \langle Q_s + \Delta _1, \theta _1 \rangle \big)
	\big(1 + r  - \langle Q_s + \Delta _2, \theta _2 \rangle \big) e^{\lambda \langle Q_s+ \Delta _1, \theta _1 \rangle
		+	 \lambda \langle Q_s+ \Delta _2, \theta _2 \rangle}
	\\
	& \times
	\1 \{|Q_k| \leq f^{t,y} _k, k = 0,...,s, |Q_s+\Delta _1|\vee |Q_s+\Delta _2| \leq f^{t,y} _{s}  \}    
	\bigg].
\end{align*}
Denote $\mck = \Big( - \frac 12, \frac 12 \Big] ^\d$
and 
$$\mci  = \mci _s  = \{ (x, y_1, y_2) \in (\Z^\d ) ^3: |x| \leq f^{t,y} _s, |x+y_i| \leq f^{t,y} _s,
i = 1,2 \}.$$
Denote also 
$$
{\mathrm{\mathbf{p}}}(x) = \P \{|Q_k| \leq f^{t,y} _k, k = 0,...,s 
| Q_s \in x + \mck \}.
$$
  
Recall that $0_\d$ is the origin of $\R ^\d$.
In this proof we will use that for any polynomial
${g :  (\R ^\d ) ^3 \to \R}$ and any $0 < c_1 < c_2$
for large $s$
\begin{equation}\label{Gliederung = outline, structure}
	\sum\limits _{(x,y_1, y_2) \in \mci _s} g(x, y_1, y_2)
	e^{-c_2 s} \leq e^{-c_1 s}.
\end{equation}
Since $Q_s$ and the pair $(\Delta_1, \Delta_2)$
are independent
we have for $x, y_1, y_2 \in \Z ^\d$
\begin{align}\label{der Haufen = heap}
	\P \{Q_s \in x +\mck, \Delta _i \in 
	y_i + \mck, i=1,2   \} =  
	\P \{Q_s \in x +\mck  \}
	\P \{ \Delta _i \in 
	y_i + \mck, i=1,2   \}.
\end{align}
We will need the following key bounds
on the probabilities on the right hand side of \eqref{der Haufen = heap}.
By Assumption \ref{lambda-existence}
\begin{align} \label{funkeln = sparkle}
	\P \{ \Delta _i \in 
	y_i + \mck, i=1,2   \} 
	\leq   C  \exp\big\{ - (\lambda + \varepsilon)(|y_1|\vee|y_2|)\big\},
\end{align}
whereas by a multidimensional local limit theorem
\cite[Theorem 1]{BR65}, or \cite[(A.1)]{Hu91}
\begin{align}
	\P \{Q_s \in x +\mck  \} 
	\leq  C s^{-\d /2} \exp\Big\{-I_1\Big(\frac{|x|}{s}\Big) s \Big\}.
\end{align}
First we consider the easier case $\alpha \geq  \frac{\pi}{360}$.
Dropping the restriction $|Q_k| \leq f^{t,y} _k, k = 0,...,s-1$ on the entire trajectory we obtain
\begin{align*}
	\mce_\alpha \leq   \, & \E \bigg[
	\big(1 + f^{t,y} _s  - \langle Q_s + \Delta _1, \theta _1 \rangle \big)
	\big(1 + f^{t,y} _s  - \langle Q_s + \Delta _2, \theta _2 \rangle \big) e^{\lambda \langle Q_s+ \Delta _1, \theta _1 \rangle
		+	 \lambda \langle Q_s+ \Delta _2, \theta _2 \rangle}
	\\
	& \times
	\1 \{|Q_s| \leq f^{t,y} _s, |Q_s+\Delta _1|\vee |Q_s+\Delta _2| \leq f^{t,y} _s  \}    
	\bigg].
\end{align*}
and hence
\begin{align}  
	\mce \leq & \,  C s^{-\d /2} \sum\limits _{(x,y_1, y_2) \in \mci }
	\bigg[ 
	\exp  \Big\{ \lambda \langle x + y_1, \theta _1 \rangle +
	\lambda \langle x + y_2, \theta _2 \rangle -I_1\Big(\frac{|x|}{s}\Big) s - (\lambda + \varepsilon)(|y_1|\vee|y_2|) \Big\}
	\\
	\notag
	& \hspace{4cm}
	\times 
	\big(1 + f^{t,y} _s  - \langle x + y _1, \theta _1 \rangle \big)
	\big(1 + f^{t,y} _s  - \langle x + y _2, \theta _2 \rangle \big)
	\bigg]
	\\
	\notag
	= & \, 
	C s^{-\d /2}
	\sum\limits _{(x,y_1, y_2) \in \mci }
	e^{\widetilde h(x,y_1, y_2)}
	\big(1 +f^{t,y} _s - \langle x + y _1, \theta _1 \rangle \big)
	\big(1 + f^{t,y} _s - \langle x + y _2, \theta _2 \rangle \big),
\end{align}
where $\widetilde{h}$ is the function defined by
$$
\widetilde h(x,y_1, y_2)  = 
\lambda \langle x + y_1, \theta _1 \rangle +
\lambda \langle x + y_2, \theta _2 \rangle -I_1\Big(\frac{|x|}{s}\Big) s - (\lambda + \varepsilon)(|y_1|\vee|y_2|).
$$
This function is closely related 
to $h$
from Lemma \ref{in Gang kommen} which was defined in dimension two.
Without loss of generality assume 
$\theta _1 = (\cos \alpha, \sin \alpha, 0, ...,0)$ and 
$\theta _2 = (\cos \alpha, -\sin \alpha, 0, ...,0)$.
We have
\[
\widetilde h(x,y_1, y_2) \leq h(\mcp  x, \mcp y_1, \mcp y_2),
\]
where ${\mcp}$ is the orthogonal projection of $\R ^\d$
onto the plane
$$
\mathfrak{L}_2 = \{ z \in \R ^ \d: z=  (u_1, u_2,0,...,0) \}
$$
spanning $\theta _1$ and $\theta _2$. Note also that  when $x, y_1, y_2$ belong to this plane, then 
\[
\widetilde h(x,y_1, y_2) = h(\mcp  x, \mcp y_1, \mcp y_2).
\]
Lemma \ref{in Gang kommen}, $(iii)$,
with $r = f^{t,y} _s$ 
and \eqref{Gliederung = outline, structure}
now imply \eqref{pull someone up for something}
for $\alpha \geq \frac{\pi}{360}$.\\
  
Next we consider  the more challenging case $\alpha  \in [0, \frac{\pi}{360}]$.
Set
$r = r_s = f^{t,y} _{s  }$.
We start off with a few observations.
Consider a copy of the random walk $(\X_t, t\geq 0)$
defined by  $\mathfrak{Y} _t= \langle Q_t, \frac{x}{|x|} \rangle $, $ 0< |x| \leq r$. 
Proceeding similarly to 
\eqref{limpid = clear, transparent}
we apply  the transformation 
\eqref{Verwahrung = custody, safekeeping} 
with $b = \lambda$ and Theorem \ref{thm beinhalten = include}
to find that for some $C > 0$
for large $s$
\begin{equation}\label{Laien = lay ppl}
	{\mathrm{\mathbf{p}}}(x) \leq \frac{C (r + 1 - |x|) (1+y)}{s}.
\end{equation}
Decomposing over the values of $Q_s$ and $\Delta _i$
we have
\begin{align} \label{straff = firm, tight}
	\mce \leq C \sum\limits _{(x,y_1, y_2) \in \mci _s}
	& \P \{Q_s \in x +\mck, \Delta _i \in 
	y_i + \mck, i=1,2   \}{\mathrm{\mathbf{p}}}(x)
	e^{\lambda \langle x + y_1, \theta _1 \rangle +
		\lambda \langle x + y_2, \theta _2 \rangle }
	\\
	\notag
	& 
	\times
	\big(1 + r  - \langle x + y _1, \theta _1 \rangle \big)
	\big(1 + r  - \langle x + y _2, \theta _2 \rangle \big)
\end{align}
The next step is to reduce the problem to  an essentially  two-dimensional case.
Recall that we assumed
$\theta _1 = (\cos \alpha, \sin \alpha, 0, ...,0)$ and 
$\theta _2 = (\cos \alpha, -\sin \alpha, 0, ...,0)$.
Note that  $ \frac{\theta _1 + \theta _2}{|\theta _1 + \theta _2|} = \mathbf{e}_1$. 
Using the inequality 
$$a + \frac{b_1+...+b_n}{4n(a \vee 1)} \leq \sqrt{a^2 + b_1 ^2 + ... + b_n^2}$$
which holds for $a \geq 0$, $n \in \N$, $b_i \geq 0$,
we obtain for
$z_1, z_2 \in \mathfrak{L}_2$ 
in the case $\d > 2$
\begin{align*}
	\sum\limits  _{y_1, y_2  \in \Z ^\d:   \mcp y_i = z_i} & 
	\mkern-14mu
	\exp  \big\{  - (\lambda + \varepsilon)(|y_1|\vee|y_2|) \big\}
	\\
	& \mkern-67mu = 
	\sum\limits _{\substack{ k_3, ..., k_\d, \ell _3, ...,\ell _\d \in \Z}}
	\mkern-14mu
	\exp  \bigg\{  - (\lambda + \varepsilon)\Big( \sqrt{|z_1|^2 + k_3 ^2 + ... + k_{\d} ^2} \vee \sqrt{|z_2|^2 + \ell_3 ^2 + ... + \ell_{\d} ^2} \Big ) \bigg\}
	\\
	& \mkern-67mu \leq 
	\sum\limits _{\substack{ k_3, ..., k_\d, \ell _3, ...,\ell _\d \in \Z}}
	\mkern-14mu
	\exp  \bigg\{  - (\lambda + \varepsilon)\bigg( \Big[  |z_1| + \frac{|k_3| + ... + |k_{\d}|}{4(\d -2)(|z_1|\vee 1)} \Big] \vee \Big[ |z_2| + \frac{|\ell_3| + ... + |\ell_{\d}|}{4(\d -2)(|z_2|\vee 1)} 
	\Big] \bigg ) \bigg\}
	\\
	& \mkern-67mu \leq
	C (|z_1|\vee|z_2| + 1 )^{2\d - 2}
	\exp  \big\{  - (\lambda + \varepsilon)(|z_1|\vee|z_2|) \big\}.
\end{align*}
\nopagebreak
(when $\d = 2$, the inequality between the first and the last expression in the above chain is still satisfied.)
Combining this with \eqref{straff = firm, tight} and
\eqref{funkeln = sparkle} 
we get
\begin{align} 
	\notag
	\mce \leq & \, C \sum\limits _{(x,y_1, y_2) \in \mci _s}
	\Big[ \P \{Q_s \in x +\mck, \Delta _i \in 
	y_i + \mck, i=1,2   \}{\mathrm{\mathbf{p}}}(x)
	e^{\lambda \langle x + y_1, \theta _1 \rangle +
		\lambda \langle x + y_2, \theta _2 \rangle
		- (\lambda + \varepsilon)(|y_1|\vee|y_2|) }
	\\
	\notag
	& 
	\mkern 78mu
	\times
	\big(1 + r  - \langle x + y _1, \theta _1 \rangle \big)
	\big(1 + r  - \langle x + y _2, \theta _2 \rangle \big) \Big]
	\\ 
	\leq & \, C \sum\limits _{\substack{(x,y_1, y_2) \in \mci _s: 
			\\
			y_1, y_2 \in \mathfrak{L}_2}}
	\Big[
	\P \{Q_s \in x +\mck \}{\mathrm{\mathbf{p}}}(x)
	e^{\lambda \langle x + y_1, \theta _1 \rangle +
		\lambda \langle x + y_2, \theta _2 \rangle
		- (\lambda + \varepsilon)(|y_1|\vee|y_2|) }
	\\ &
	\mkern 78mu
	\times
	\big(1 + r  - \langle x + y _1, \theta _1 \rangle \big)
	\big(1 + r  - \langle x + y _2, \theta _2 \rangle \big)
	(|y_1|\vee|y_2| + 1 )^{2\d-2}\Big]
\end{align}

Denote $\mck _2 = \Big( - \frac 12, \frac 12 \Big] ^2 \times \{  \underbrace{(0,...,0)}_{\d-2} \}
= \mcp \mck$ and let $ \mci _s ^{(2)} = \mci _s \cap \mathfrak{L}_2$.
Since for $|x_1| \leq |x_2|$ the inequality ${\mathrm{\mathbf{p}}}(x _1) \geq {\mathrm{\mathbf{p}}}(x_2)$ holds true, we have
\begin{align}
	\mce
	\leq  & \, C \sum\limits _{\substack{(x,y_1, y_2) \in \mci _s: 
			\\
			y_1, y_2 \in \mathfrak{L}_2}}
	\Big[
	\P \{\mcp Q_s \in \mcp x +\mck  \}{\mathrm{\mathbf{p}}}( \mcp x)
	e^{\lambda \langle x + y_1, \theta _1 \rangle +
		\lambda \langle x + y_2, \theta _2 \rangle 
		- (\lambda + \varepsilon)(|y_1|\vee|y_2|)}
	\\ &
	\mkern 78mu
	\times
	\big(1 + r  - \langle x + y _1, \theta _1 \rangle \big)
	\big(1 + r  - \langle x + y _2, \theta _2 \rangle \big)
	(|y_1|\vee|y_2| + 1 )^{2\d-2}\Big]
	\\ 
	=  & \, C \sum\limits _{(x,y_1, y_2) \in \mci _s ^{(2)}}
	\Big[
	\P \{\mcp Q_s \in  x +\mck _2 \}{\mathrm{\mathbf{p}}}( x)
	e^{\lambda \langle x + y_1, \theta _1 \rangle +
		\lambda \langle x + y_2, \theta _2 \rangle
		- (\lambda + \varepsilon)(|y_1|\vee|y_2|) }
	\\ &
	\mkern 78mu
	\times
	\big(1 + r  - \langle x + y _1, \theta _1 \rangle \big)
	\big(1 + r  - \langle x + y _2, \theta _2 \rangle \big)
	(|y_1|\vee|y_2| + 1 )^{2\d-2}\Big]
\end{align}
With this we have reduced the $\d$-dimensional case
to essentially a two-dimensional case.
Over the following couple of pages we focus on deriving a useful bound on $ \widetilde h(x,y_1, y_2)$.
Let $\mcq _s$ be the  projection 
of $Q_s$ onto $\mathfrak{L}_2$:
$\mcq _s = \mcp Q_s$.
Then by the above
\begin{align}
	\mce \leq 
	& \, C \sum\limits _{(x,y_1, y_2) \in \mci _s ^{(2)}}
	\Big[
	\P \{\mcq _s \in  x +\mck _2 \}{\mathrm{\mathbf{p}}}( x)
	e^{\lambda \langle x + y_1, \theta _1 \rangle +
		\lambda \langle x + y_2, \theta _2 \rangle - (\lambda + \varepsilon)(|y_1|\vee|y_2|) }
	\\ &
	\mkern 78mu
	\times
	\big(1 + r  - \langle x + y _1, \theta _1 \rangle \big)
	\big(1 + r  - \langle x + y _2, \theta _2 \rangle \big)
	(|y_1|\vee|y_2| + 1 )^{2\d-2}\Big].
\end{align}
By Assumption \ref{lambda-existence}
and a multidimensional local limit theorem
\cite[Theorem 1]{BR65}, or \cite[(A.1)]{Hu91}, 
\begin{equation}\label{schraeg = skewed, oblique}
	\P \{\mcq _s \in x +\mck _2  \}   
	\leq  C s^{-1} \exp\Big\{-I_1\Big(\frac{|x|}{s}\Big) s\Big\}.
\end{equation}
Therefore
\begin{align} 
	\mce \leq & \,
	C s^{-1} \sum\limits _{(x,y_1, y_2) \in \mci _s ^{(2)}}
	\mkern-12mu
	\Big[ e^{\widetilde h(x,y_1, y_2)}
	{\mathrm{\mathbf{p}}}(x)
	\\
	&
	\times
	\big(1 + r  - \langle x + y _1, \theta _1 \rangle \big)
	\big(1 + r  - \langle x + y _2, \theta _2 \rangle \big)
	(|y_1|\vee|y_2| + 1 )^{2\d-2} \Big].
\end{align}
Denote by $\gamma _i $ the angle 
between $\mathbf{e}_1 $ and $x +y_i$, 
$- \pi < \gamma _i \leq \pi$.
We adopt the orientation under which the angle between $\mathbf{e}_1$ and $\theta _1$
equals $\alpha$, whereas the angle between $\mathbf{e}_1$ and $\theta _2$
equals $-\alpha$.
Set 
\begin{multline}
	\mci _{360}
	= \Big\{ (x, y_1, y_2) \in \mci _s ^{(2)}: -\frac{\pi}{360}  \leq \gamma _2 \leq \gamma _1 \leq \frac{\pi}{360},
	\\
	|x| \leq f^{t,y} _{s}, |x + y_i| \leq f^{t,y} _{s},
	i = 1,2, |x + y_1| = \mathfrak{m}, |y_1| \geq |y_2| \Big\}.
\end{multline}
Rotation and symmetry considerations imply that  for large $s$ 
\begin{align} 
	\label{ausgiebig = extensive, thorough}
	\mce \leq & \,
	C s^{-1} \sum\limits _{(x,y_1, y_2) \in \mci _{360}}
	\mkern-12mu
	\Big[ e^{ \widetilde h(x,y_1, y_2)}
	{\mathrm{\mathbf{p}}}(x)
	\\
	&
	\times
	\big(1 + r  - \langle x + y _1, \theta _1 \rangle \big)
	\big(1 + r  - \langle x + y _2, \theta _2 \rangle \big)
	(|y_1|\vee|y_2| + 1 )^{2\d-2} \Big].
\end{align}
Recall that $r = r_s = f^{t,y} _{s}$,
$ \mathfrak{m} = |x + y _1| \vee |x + y _2| $ and $2 \mathfrak{d} = |y_1 - y_2|  $.
Define also $\partial = \mathfrak{m} - 
|x + y _1| \wedge |x + y _2|$. Note that
for $(x,y_1, y_2) \in \mci _{360}$
\[
\mathfrak{m} = |x + y _1| 
\ \ \ 
\text{and} 
\ \ \ 
\partial = \mathfrak{m} - |x + y _2|.
\] 
Note that
$$
\langle x + y_1,\theta _1 \rangle =  |x + y_i| \cos ( \alpha - \gamma _1), \ \ \ 
\langle x + y_2,\theta _2 \rangle =  |x + y_2| \cos ( -\alpha - \gamma _2),
$$
and hence 
\begin{equation}\label{das Kreuz}
	\langle x + y_1, \theta _1 \rangle +
	\langle x + y_2, \theta _2 \rangle \leq 
	\mathfrak{m} (\cos ( \alpha - \gamma _1) + \cos ( \alpha + \gamma _2)) - \partial \cos ( \alpha - \gamma _2).
\end{equation}
Let us also note that since
$- \frac{\pi}{360} \leq \gamma _2, \alpha \leq  \frac{\pi}{360}$,
\[
\partial \cos ( \alpha - \gamma _2) \geq 0.9 \partial.
\]
The above inequalities give  a bound on the `gain' part of $\widetilde{h}$.
Next we proceed to bound the `cost' part. 
Since $|y_1|\vee|y_2| = |y_1|$
for $(x,y_1, y_2) \in \mci _{360}$, similarly to  \eqref{das Loch}
we have
for large $s$
\begin{equation}\label{das Loch2}
	I\Big(\frac xs \Big) s + (\lambda +\varepsilon) (|y_1|\vee |y_2|) \geq 
	I_1\Big(\frac {\mathfrak{m} - |y_1|}{s} \Big) s + (\lambda +\varepsilon) |y_1|
	\geq I_1\Big(\frac{\mathfrak{m}}{s} \Big) s + 
	0.99 \varepsilon |y_1|.
\end{equation}
Moreover, for $(x,y_1, y_2) \in \mci _{360}$
we have 
$$
2 \mathfrak{d} = |y_1 - y_2| \leq 2 |y_1|
$$
and since 
$- \frac{\pi}{360} \leq \gamma _2 \leq \gamma _1 \leq \frac{\pi}{360}$ (see  Figure \ref{f_for_proof2} for an illustration)
\begin{equation}\label{gaeren = ferment, fester}
	|y_1| \geq  \mathfrak{d} \geq
	\frac {\mathfrak{m}}{2} \sin (\gamma _1 - \gamma _2) \geq 0.48\mathfrak{m} (\gamma _1 - \gamma _2).
\end{equation}
\begin{figure}[t!]
	\vspace{-4.8cm}
	\centering
	\includegraphics[scale=0.75]{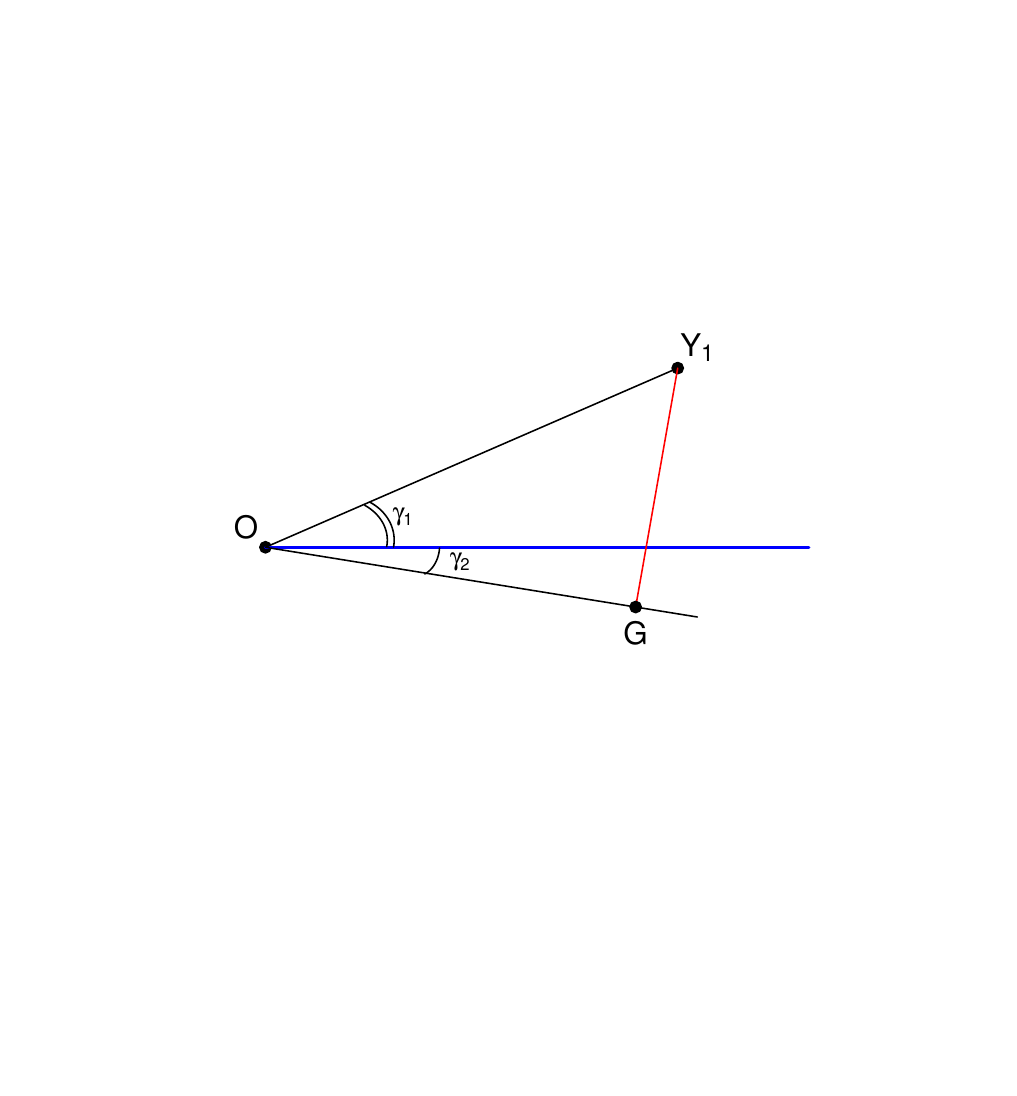}
	\vspace{-5.8cm}
	\caption{ \small
		An illustration to \eqref{gaeren = ferment, fester}. 
		The point coordinates are 
		${\rm{O}} = (0,0)$,
		and ${\rm{Y}}_1 = (\mathfrak{m} \cos \gamma, \mathfrak{m} \sin \gamma) = x + y_1$. We have $|x+y_1| = \mathfrak{m}$. In this picture $\gamma _1 > 0$, $\gamma _2 < 0$. The point ${\rm{G}} $ is then chosen so that   $\angle {\rm{O}} {\rm{G}} {\rm{Y}}_1 = \frac \pi 2$. No matter where on the line ${\rm{OG}} $ the point ${\rm{Y}}_2 = x + y_2$ lies,
		$2 \mathfrak{d} =  |{\rm{Y_1Y_2}}| \geq |{\rm{Y_1G}}|  = \mathfrak{m}  \sin (\gamma _1 - \gamma _2).$
	}
	\label{f_for_proof2}
	\vspace{0.21cm}
\end{figure}
Combining these bounds with \eqref{das Kreuz} we find 
\begin{equation}
	\widetilde h(x,y_1, y_2) \leq \lambda \mathfrak{m} (\cos ( \alpha - \gamma _1) + \lambda  \mathfrak{m} \cos (- \alpha - \gamma _2))
	- 0.9 \partial 
	- 
	I_1\Big(\frac{\mathfrak{m}}{s} \Big) s -
	0.5 \varepsilon |y_1|
	- 0.22 
	\mathfrak{m}\varepsilon  (\gamma _1 - \gamma _2).
\end{equation}
Elementary inequalities 
give
\begin{equation*}
	\lambda \cos ( \alpha - \gamma _1) + \lambda \cos ( -\alpha - \gamma _2) 
	\leq 2 \lambda  - 0.3 \lambda( \alpha - \gamma _1) ^2 
	- 0.3\lambda ( \alpha + \gamma _2)^2.
\end{equation*}
Therefore 
\begin{multline}
	\widetilde h(x,y_1, y_2)
	 \leq  2\lambda \mathfrak{m} 
	- 0.3 \lambda \mathfrak{m}( \alpha - \gamma _1) ^2 
	- 0.3 \lambda \mathfrak{m}( \alpha + \gamma _2)^2
	- 0.9 \partial 
	\\
	- 
	I_1\Big(\frac{\mathfrak{m}}{s} \Big) s -
	0.5 \varepsilon |y_1|
	- 0.22 
	\mathfrak{m}\varepsilon  (\gamma _1 - \gamma _2).
\end{multline}
Furthermore, for some $c>0$
\[
0.3 \lambda( \alpha - \gamma _1) ^2
+ 0.3 \lambda( \alpha + \gamma _2)^2
+ 0.1 
\varepsilon  (\gamma _1 - \gamma _2)
\geq c \alpha ^2 + c (\gamma_1 + \gamma_2)^2
\]
and hence 
\begin{equation}
	\widetilde h(x,y_1, y_2) \leq  2\lambda \mathfrak{m} 
	- 0.9 \partial 
	- c \mathfrak{m} \alpha ^2 - c \mathfrak{m} (\gamma_1 + \gamma_2)^2
	-
	I_1\Big(\frac{\mathfrak{m}}{s} \Big) s -
	0.5 \varepsilon |y_1|
	- 0.12
	\mathfrak{m}\varepsilon  (\gamma _1 - \gamma _2).
\end{equation}
Since $r = \frac{\Psi(\lambda)}{\lambda} s 
+ o(s) $ and $I_1'\Big(\frac{\Psi(\lambda)}{\lambda}\Big) = \lambda$
and because $I_1$ is convex
for large $s$ we have  ${I_1'\big(\frac{\mathfrak{m}}{s} \big) \leq \lambda +  \frac{\varepsilon} {100}}$
whenever $0\leq \mathfrak{m} \leq r$.
Consequently 
for large $s$
\[
- 
I_1\Big(\frac{\mathfrak{m}}{s} \Big)s
+ I_1\Big(\frac{r}{s} \Big)s \leq 
\Big(\lambda +  \frac{\varepsilon} {100} \Big)(r - \mathfrak{m})
\]
and 
\[
2\lambda \mathfrak{m} - 
I_1\Big(\frac{\mathfrak{m}}{s} \Big)s
- h_{\text{max}} = 
2\lambda \mathfrak{m} - 
I_1\Big(\frac{\mathfrak{m}}{s} \Big)s
- \Big[ 2\lambda r - I_1 \Big(\frac{r}{s} \Big)s\Big]
\leq 
-\Big(\lambda - \frac{\varepsilon} {100}\Big) (r - \mathfrak{m})
\]
With this inequality we get
our final bound on $\widetilde h(x,y_1, y_2)$:
\begin{equation*}
	\widetilde h(x,y_1, y_2) \leq
	h_{\text{max}} - 
	\Big(\lambda - \frac{\varepsilon} {100}\Big) (r - \mathfrak{m})
	- 0.9 \partial 
	- c \mathfrak{m} \alpha ^2 - c \mathfrak{m} (\gamma_1 + \gamma_2)^2
	-
	0.5 \varepsilon |y_1|
	- 0.12
	\mathfrak{m}\varepsilon  (\gamma _1 - \gamma _2).
\end{equation*}
With this bound we now derive the desired
inequality for $\mce$. 
Let $G$ be a shorthand for 
$$\big(1 + r  - \langle x + y _1, \theta _1 \rangle \big)
\big(1 + r  - \langle x + y _2, \theta _2 \rangle \big)
(|y_1|\vee|y_2| + 1 )^{2\d-2}.$$
Now
we apply the above bound on $\widetilde h(x,y_1, y_2)$ to \eqref{ausgiebig = extensive, thorough} and use
domination of sums with exponentially decaying elements 
by the first terms:
\begin{align}
	s \mce & \leq 
	\sum\limits _{(x,y_1, y_2) \in \mci _{360}} 
	\mkern-12mu
	e^{ \widetilde h(x,y_1, y_2)}
	{\mathrm{\mathbf{p}}}(x)
	G
	\\ 
	& \leq 
	\sum\limits _{(x,y_1, y_2) \in \mci _{360}} 
	\mkern-12mu
	\exp  \begin{multlined}[t]
		\Big\{h_{\text{max}} - 
		\Big(\lambda - \frac{\varepsilon} {100}\Big) (r - \mathfrak{m})
		- 0.9 \partial 
		- c \mathfrak{m} \alpha ^2 - c  \mathfrak{m} (\gamma_1 + \gamma_2)^2
		\\
		-
		0.5 \varepsilon |y_1|
		- 0.12
		\mathfrak{m}\varepsilon  (\gamma _1 - \gamma _2) \Big\}
		{\mathrm{\mathbf{p}}}(x) G\end{multlined}
	\\ 
	\notag
	& \leq 
	C \mkern-13mu \sum\limits _{\substack{(x,y_1, y_2) \in \mci _{360}:
			\\ \partial \leq 10 \d}} 
	\mkern-12mu
	\exp   \begin{multlined}[t]\Big\{h_{\text{max}} - 
		\Big(\lambda - \frac{\varepsilon} {100}\Big) (r - \mathfrak{m})
		- c \mathfrak{m}\alpha ^2 - c\mathfrak{m} (\gamma_1 + \gamma_2)^2
		\\
		\ \hspace{2.2cm} \mkern 37mu  - 
		0.5 \varepsilon |y_1|
		- 0.12
		\mathfrak{m}\varepsilon  (\gamma _1 - \gamma _2)  \Big\}
		{\mathrm{\mathbf{p}}}(x)
		G \end{multlined}
	\\ 
	\notag
	& \leq 
	C \mkern-17mu \sum\limits _{\substack{(x,y_1, y_2) \in \mci _{360}:
			\\ \partial \leq 10 \d, r - \mathfrak{m} \leq 10 \d }} 
	\mkern-27mu
	\exp  \begin{multlined}[t]\Big\{h_{\text{max}}  
	- c r \alpha ^2 - cr (\gamma_1 + \gamma_2)^2
	\\
    \
	 \hspace{2.2cm} \mkern 37mu -
	0.5 \varepsilon |y_1|
	- 0.12
	r\varepsilon  (\gamma _1 - \gamma _2)  \Big\}
	{\mathrm{\mathbf{p}}}(x)
	G\end{multlined}
	\\ 
	\label{jdn troesten = comfort smb}
	& \leq 
	C \mkern-17mu \sum\limits _{\substack{(x,y_1, y_2) \in \mci _{360}:
			\\ \partial \leq 10 \d, r - \mathfrak{m} \leq 10 \d, \\
			\gamma_1 = \gamma_2 = \gamma  }} 
	\mkern-12mu
	\exp  \Big\{h_{\text{max}}  
	- c r\alpha ^2 - 4 cr \gamma^2
	-
	0.5 \varepsilon |y_1| \Big\}
	{\mathrm{\mathbf{p}}}(x) G
	\\ 
	\notag
	& \leq 
	C \mkern-17mu \sum\limits _{\substack{(x,y_1, y_2) \in \mci _{360}:
			\\ \partial \leq 10 \d, r - \mathfrak{m} \leq 10 \d, \\
			\gamma_1 = \gamma_2 = \gamma, y_1 = 0_2  }} 
	\mkern-12mu
	\exp  \Big\{h_{\text{max}}  
	- c r\alpha ^2 - 4 cr \gamma^2
	\Big\}
	{\mathrm{\mathbf{p}}}(x) G
	\\
  &	= C e^{h_{\text{max}} - c r \alpha ^2  }
	\mkern-32mu
	\sum\limits _{\substack{(x,0_2, 0_2) \in \mci _{360}:
			\\  r - \mathfrak{m} \leq 10 \d,
			\gamma_1 = \gamma_2 = \gamma  }} 
	\mkern-22mu
	e^{ - 4 c r \gamma^2 }
	{\mathrm{\mathbf{p}}}(x) G.
\end{align}
Next we work on the last sum.
For $(x,0_2, 0_2) \in \mci _{360}$
with $r - \mathfrak{m} \leq 10 \d$
we have $|x| \geq r - 10$, and hence 
by \eqref{Laien = lay ppl}
\[
{\mathrm{\mathbf{p}}}(x) \leq \frac{C (1+y)}{s}.
\]
Consequently
\begin{align}
	\notag
	\sum\limits _{\substack{(x,0_2, 0_2) \in \mci _{360}:
			\\   r - \mathfrak{m} \leq 10 \d, \\
			\gamma_1 = \gamma_2 = \gamma  }} 
	\mkern-12mu
	e^{ - 4 c r \gamma^2 }
	{\mathrm{\mathbf{p}}}(x) G 
	& \leq 
	\frac{C (1+y)}{s}
	\sum\limits _{\substack{(x,0_2, 0_2) \in \mci _{360}:
			\\   r - \mathfrak{m} \leq 10 \d, \\
			\gamma_1 = \gamma_2 = \gamma  }} 
	\mkern-12mu
	e^{ - 4 c r \gamma^2 } \big(1 + r  - \langle x + y _1, \theta _1 \rangle \big)
	\big(1 + r  - \langle x + y _2, \theta _2 \rangle \big)
	\\ \notag
	& = 
	\frac{C (1+y)}{s}
	\sum\limits _{\substack{(x,0_2, 0_2) \in \mci _{360}:
			\\   r - \mathfrak{m} \leq 10 \d, \\
			\gamma_1 = \gamma_2 = \gamma  }} 
	\mkern-12mu
	e^{ - 4 c r \gamma^2 } \big(1 + r (1 - \cos(\gamma - \alpha))\big)\big(1 + r (1 - \cos( \gamma + \alpha ))  \big)
	\\ \notag
	& \leq  
	\frac{C (1+y)}{s}
	\sum\limits _{\substack{(x,0_2, 0_2) \in \mci _{360}:
			\\   r - \mathfrak{m} \leq 10 \d, \\
			\gamma_1 = \gamma_2 = \gamma \geq 0  }} 
	\mkern-12mu
	e^{ - 4 c r \gamma^2 } \big(1 + r (1 - \cos( \gamma + \alpha ))^2  \big)
	\\ 
	\label{der Ausschlag = rash}
	& \leq  
	\frac{C (1+y)}{s}
	\sum\limits _{\substack{(x,0_2, 0_2) \in \mci _{360}:
			\\   r - \mathfrak{m} \leq 10 \d, \\
			\gamma_1 = \gamma_2 = \gamma \geq 0  }} 
	\mkern-12mu
	e^{ - 4 c r \gamma^2 } \big(1 + 0.3 r ( \gamma + \alpha )^2  \big).
\end{align}
We estimate the last sum using a comparison with an integral and then the substitution ${ z = r^{1/2}\gamma}$.
We have
\begin{align}
	\notag
	\sum\limits _{\substack{(x,0_2, 0_2) \in \mci _{360}:
			\\   r - \mathfrak{m} \leq 10 \d, \\
			\gamma_1 = \gamma_2 = \gamma \geq 0  }} 
	\mkern-12mu
	e^{ - 4 c r \gamma^2 } \big(1 + 0.3 r ( \gamma + \alpha )^2  \big) 
	& \leq C r\int\limits _{0} ^{\frac{\pi}{360}}
	e^{ - 4 c r \gamma^2 } \big(1 + 0.3 r ( \gamma + \alpha )^2  \big) d\gamma
	\\[-15pt] \notag
	& \leq C r^{1/2} \int\limits _{0} ^{\infty}
	e^{ - 4 c z^2 } \big(1 + 0.3 r ( r^{-1/2}z + \alpha )^2  \big) dz
	\\
	& \leq C r^{1/2} \int\limits _{0} ^{\infty}
	e^{ - 4 c z^2 } \big(1 + 0.6 r \alpha ^2  + 0.6 z^2 \big) dz
	\\
	& \leq C r^{1/2} (1 +  r \alpha ^2 ).
\end{align}
Combining this with \eqref{jdn troesten = comfort smb}
and \eqref{der Ausschlag = rash}
and using that for some $ c_1 \in (0,c)$ and $C>0$
\[
q e^{- c q  } \leq C e^{- c_1 q  },
\ \  \ q \geq 0, \]
we get
\begin{equation}\label{naschen = nibble}
	\mce \leq C s^{-1} e^{h_{\text{max}} - c r \alpha ^2  }
	\frac{ (1+y)}{s}
	r^{1/2} (1 +  r \alpha ^2 )
	\leq  \frac{ C(1+y)}{s ^{3/2}} e^{h_{\text{max}} - c_{_1} r \alpha ^2  }.
\end{equation}
Recall that  $r = f^{t,y} _{s}$ and
\[
e^{h_{\text{max}}} = \exp \Big\{2 \lambda r - I_1\Big( \frac rs \Big)s \Big\}.
\]
Since 
$I_1\Big(\frac{\Psi(\lambda)}{\lambda} \Big) = \ln m $,
$I_1'\Big(\frac{\Psi(\lambda)}{\lambda} \Big) = \lambda $,
and uniformly in $t$, $t \geq s$, $\frac rs = \frac{f^{t,y} _{s}}{s} = \frac{\Psi(\lambda)}{\lambda} + o(s)  $,
we have
\begin{equation*}
	I_1\Big(\frac{f^{t,y} _{s}}{s} \Big) = I_1\Big(\frac{\Psi(\lambda)}{\lambda} \Big)
	+ \lambda \Big(\frac{f^{t,y} _{s}}{s} - \frac{\Psi(\lambda)}{\lambda} \Big) + O\Big(\frac{f^{t,y} _{s}}{s} - \frac{\Psi(\lambda)}{\lambda} \Big)^2 
\end{equation*}
and hence
\begin{equation}\label{die Faehre = ferry}
	I_1\Big(\frac{r}{s} \Big)s = 
	I_1\Big(\frac{f^{t,y} _{s}}{s} \Big)s =  s \ln m
	+ \lambda \Big(f^{t,y} _{s} - \frac{\Psi(\lambda)}{\lambda}s \Big) + o(1). 
\end{equation}
Writing down  
\begin{align*}
	\exp\Big\{ - \lambda \Big( f^{t,y} _{s} - \frac{\Psi(\lambda)}{\lambda}s \Big)\Big\}  =
	& \,
	\exp\Big\{ -\frac{\d-1}{2} \ln (s+1)
	+ \frac{3}{2} \ln \frac{t+1}{t-s +1} -
	{3} M -  y\Big\}
	\\
	= & \
	e^{-3M -  y}
	\frac{  (t+1) ^{3/2}}{(s+1)^{(\d-1)/2} (t-s+1) ^{3/2}}
\end{align*}
we obtain 
by \eqref{naschen = nibble}
and \eqref{die Faehre = ferry}
\begin{align}
	\mce \leq 
	\frac{ C(1+y)}{s ^{3/2}} & \, e^{2 \lambda r - c_{_1} r \alpha ^2} \exp\Big\{ -I_1\Big(\frac{r}{s} \Big)s \Big\}
	\\
	&\leq 
	\frac{ C(1+y)e^{ 2 \lambda r - c_{_1} r \alpha ^2  }}{s ^{3/2}} m^{-s}
	e^{-  y}
	\frac{  (t+1) ^{3/2}}{(s+1)^{(\d-1)/2} (t-s+1) ^{3/2}}
	\\
	&\leq 
	\frac{ C(1+y)e^{ 2 \lambda r - c_{_1} r \alpha ^2 -y  }}{m^2} 
	\times 
	\frac{  (t+1) ^{3/2}}{(s+1)^{(\d+2)/2} (t-s+1) ^{3/2}},
\end{align}
and \eqref{pull someone up for something} is proven.
\end{proof}

\begin{lem}\label{umbrage; take an umbrage new}
	Let $\theta, \theta ' \in \S ^{\d - 1}$, 
	and let $2 \alpha \in [0, \pi]$
	be the angle between $\theta$ and $\theta'$, 
	$\cos (2\alpha) = \langle \theta, \theta ' \rangle$
	There exists $A > 0$ such that for $t \geq 1$
	and $\theta, \theta ' \in \S ^{\d - 1}$
	satisfying $|\theta - \theta '| \geq A t^{-\frac 12}$
	\begin{equation}
		\E  [\#  \mathcal{A}^{t,y}_\theta \# \mathcal{A}^{t,y}_{\theta'}]
		\leq C (1+y)	e^{-y} \Big[(t+1)^{-\d+1} \alpha ^{-\d +2} + 
		(t+1)^{-\frac{\d-1}{2}} e^{-c \alpha ^2 t} \, \Big],
	\end{equation}
	for some $C,c> 0$.
\end{lem}
  
\begin{proof}
Take $A$ so large that 
\begin{equation}\label{indolent = lazy}
	\{x \in \R ^\d: |x| \leq f_t^{t,y}, \langle x,\theta \rangle \geq  f_t^{t,y} - 1 \} \cap\{x \in \R ^\d: |x| \leq f_t^{t,y},
	\langle x,\theta ' \rangle \geq  f_t^{t,y} - 1 \} = \varnothing 
\end{equation}
whenever $|\theta - \theta '| \geq A t^{-\frac 12}$, $t \in \N$.
By \eqref{indolent = lazy} a.s.
$\mathcal{A}^{t,y}_\theta \cap \mathcal{A}^{t,y}_{\theta'} = \varnothing $
if $|\theta - \theta '| \geq A t^{-\frac 12}$.
The many-to-two lemma now gives
\begin{align}  \label{trough}
	\E  [\# & \mathcal{A}^{t,y}_\theta \# \mathcal{A}^{t,y}_{\theta'}]   = \E [\# \mathcal{A}^{t,y}_\theta \cap \mathcal{A}^{t,y}_{\theta'}]
	\\
	\notag
	& + 
	m_2 \sum\limits _{s=1} ^{t-1} m ^{2t-s-2} 
	\P\Big\{|Q ^{\langle s \rangle} _k| \leq f^{t,y} _k, \langle Q ^{\langle s \rangle} _t, \theta \rangle
	\geq f^{t,y} _t - 1, |Q^{[s]}_k| \leq f^{t,y} _k, \langle Q^{[s]}_t, \theta' \rangle
	\geq f^{t,y} _t - 1, k = 0,...,t  \Big\}
	\\
	\notag
	= & \, m_2
	\sum\limits _{s=1} ^{t-1} m ^{2t-s-2} 
	\P\Big\{|Q ^{\langle s \rangle} _k| \leq f^{t,y} _k, \langle Q ^{\langle s \rangle} _t, \theta \rangle
	\geq f^{t,y} _t - 1, |Q^{[s]}_k| \leq f^{t,y} _k, \langle Q^{[s]}_t, \theta' \rangle
	\geq f^{t,y} _t - 1, k = 0,...,t  \Big\}. 
\end{align}
By the Markov property 
\begin{multline*}
	\P\Big\{|Q ^{\langle s \rangle}_k| \leq f^{t,y} _k, \langle Q ^{\langle s \rangle}_t, \theta \rangle
	\geq f^{t,y} _t - 1, |Q^{[s]}_k| \leq f^{t,y} _k, \langle Q^{[s]}_t, \theta '\rangle
	\geq f^{t,y} _t - 1, k = 0,...,t  \Big\}
	\\
	\leq  \E \Big[ \phi_{s} \big(\langle Q_s + \Delta _1, \theta \rangle\big) \phi_{s} \big(\langle Q_s + \Delta _2, \theta ' \rangle\big)\1 \{|Q_k| \leq f^{t,y} _k, k = 0,...,s,
	\\
	|Q_s+\Delta _1|\vee |Q_s+\Delta _2| \leq f^{t,y} _{s+1}  \} 
	\Big].
\end{multline*}
with
$\phi_s$  defined in \eqref{phi def}.
Hence
by \eqref{disparate} we get 
\begin{multline}\label{dovish = conciliatory}
	\P\Big\{|Q ^{\langle s \rangle}_k| \leq f^{t,y} _k, \langle Q ^{\langle s \rangle}_t, \theta \rangle
	\geq f^{t,y} _t - 1, |Q^{[s]}_k| \leq f^{t,y} _k, \langle Q^{[s]}_t, \theta' \rangle
	\geq f^{t,y} _t - 1, k = 0,...,t  \Big\}
	\\ 
	\begin{aligned}
		\leq \,  &
		C  m^{-2(t-s)} \frac{\Phi(\lambda) ^{-2s}
			e^{-2y}
		} {(t-s + 1)^{3}(t+1)^{{(\d-4)}}} 
		\times
		\E \bigg[
		\big(1 + f^{t,y} _{s+1}  - \langle Q_s + \Delta _1, \theta \rangle \big)
		\\
		&
		\times
		\big(1 + f^{t,y} _{s+1}  - \langle Q_s + \Delta _2, \theta ' \rangle \big) e^{\lambda \langle Q_s+ \Delta _1, \theta \rangle
			+	 \lambda \langle Q_s+ \Delta _2, \theta ' \rangle}
		\\
		& \times
		\1 \{|Q_k| \leq f^{t,y} _k, k = 0,...,s, |Q_s+\Delta _1|\vee |Q_s+\Delta _2| \leq f^{t,y} _{s+1}  \}    
		\bigg].
	\end{aligned}
\end{multline}
In the last expectation we recognize $\mce_\alpha$
from Lemma \ref{pauper = poor person}. 
Hence 
\begin{multline}
	\P\Big\{|Q ^{\langle s \rangle}_k| \leq f^{t,y} _k, \langle Q ^{\langle s \rangle}_t, \theta \rangle
	\geq f^{t,y} _t - 1, |Q^{[s]}_k| \leq f^{t,y} _k, \langle Q^{[s]}_t, \theta' \rangle
	\geq f^{t,y} _t - 1, k = 0,...,t  \Big\}
	\\ 
	\leq 
	C  m^{-2(t-s)} \frac{\Phi(\lambda) ^{-2s}
		e^{-2y}
	} {(t-s + 1)^{3}(t+1)^{{(\d-4)}}} 
	\times
	\frac{C e^{2 \lambda f^{t,y} _{s} -c \alpha ^2s -y }}{m^s}
	\times \frac{(t+1) ^{3/2} (1+y)}{(s+1)^{(\d+2)/2} (t-s+1) ^{3/2}}.
\end{multline}
Now \eqref{trough}  gives
\begin{align}
	\notag
	\E  [\#  \mathcal{A}^{t,y}_\theta \# \mathcal{A}^{t,y}_{\theta'}]  \leq & \, C \sum\limits _{s=1} ^{t-1} m ^{2t-s} 
	m^{-2(t-s)} \frac{\Phi(\lambda) ^{-2s}
		e^{-2y} (1+y)
	} {(t-s+1)^{3}(t+1)^{{(\d-4)}}} 
	\\ & 
	\times \underbrace{
		\frac{\Phi(\lambda) ^{2s}(s+1)^{\d-1}(t-s+1)^{3} e^{2y}}{(t+1)^{3}}
	} _{e^{2 \lambda f^{t,y} _{s} }}
	\times
	\frac{ e^{ - c \alpha ^2 s - y} (t+1) ^{3/2}}{m^s  (s+1)^{(\d+2)/2} (t-s+1) ^{3/2} }
	\notag
	\\
	\leq & \, C (1+y)	e^{-y} \sum\limits _{s=1} ^{t-1} \frac{ (s+1)^{\d/2 - 2}  
	} {(t+1)^{{\d- 5/2}}(t-s+1)^{3/2}}   e^{-c \alpha ^2 s}.
	\label{per diem = per day}
\end{align}
To bound the last sum
we split it into two parts. We have
\begin{align}
	\sum\limits _{ 1 \leq s \leq t/2} \frac{ (s+1)^{\d/2 - 2}  
	} {(t+1)^{{\d- 5/2}}(t-s+1)^{3/2}}   e^{-c \alpha ^2 s}
	\leq & \,  C (t+1)^{-\d+1} 
	\sum\limits _{ 1 \leq s \leq t/2}
	(s+1)^{\d/2 - 2}   e^{-c \alpha ^2 s}
	\notag 
	\\ \notag 
	\leq & \,
	C (t+1)^{-\d+1} 
	\int\limits _0 ^ \infty (s+1)^{\d/2 - 2}   e^{-c \alpha ^2 s} ds
	\\ \label{neidisch auf}
	\leq & \,
	C (t+1)^{-\d+1} \alpha ^{-\d +2}
\end{align}
and since $\sum\limits _{n \in \N}\frac{1}{n^{3/2}} < \infty$
\begin{align}
	\sum\limits _{ t/2 \leq s \leq t}  \frac{ (s+1)^{\d/2 - 2}  
	} {(t+1)^{{\d- 5/2}}(t-s+1)^{3/2}}   e^{-c \alpha ^2 s}
	\leq & \,  C (t+1)^{\d/2 - 2-\d+ 5/2} 
	e^{-c \alpha ^2 t/2}
	\sum\limits _{ t/2 \leq s \leq t}
	\frac{1}{(t-s+1)^{3/2}}
	\notag 
	\\ \label{nachdenken ueber}
	\leq & \,
	C (t+1)^{-\d/2 +1/2} e^{-c \alpha ^2 t/2}.
\end{align}
The statement of the lemma follows from \eqref{per diem = per day},
\eqref{neidisch auf}, and \eqref{nachdenken ueber}.
\end{proof}

Recall that 
\begin{equation*}
	f^{t,y}_t 
	= 
	\frac{\Psi (\lambda)}{\lambda}t
	+\frac{\d-4}{2\lambda} \ln (t+1) + \frac{3M}{\lambda} + \frac y\lambda.
\end{equation*}

\subsection{Proof of Proposition \ref{prop fartherst point lower bound}: the final part}

\begin{prop}\label{volcano plume}
	There exists $C>0$ such that
	for large $t \in \N$ and $y \in [1,t^{1/2}]$
	\begin{equation}
		\P \left\{ \exists u \in \T _t:
		|X_t(u)| \geq r_t+ y \right\}
		\geq \frac{(1+y) e^{-y}}{C},
	\end{equation}
	where we recall \eqref{rdef} for the definition of $r_t$.
\end{prop}
  
\begin{proof}
Let $A>0$ be the constant from the proof of Lemma \ref{umbrage; take an umbrage new}, so that \eqref{indolent = lazy} is satisfied.
Take a sequence $ \{\mathscr{K}_t, t \in \N \}$ of finite subsets of $\S ^{d-1}$ such that 
\begin{itemize}
	\item $\sup\limits _{t \in \N} \frac{\# \mathscr{K}_t}{t ^{(\d-1)/2}} < \infty$, 
	\  
	$\inf\limits _{t \in \N} \frac{\# \mathscr{K}_t}{t ^{(\d-1)/2}}  > 0 \ \  $,
	\item  For $\theta, \kappa \in \mathscr{K}_t$, $\theta \ne \kappa$, we have $|\theta - \kappa| \geq \frac{A}{\sqrt {t}}$.
\end{itemize} 
For a given $t \in \N$ a set $\mathscr{K}_t$
can be obtained for example as follows: 
choose the first point $x_1$ randomly uniformly on 
$\S ^{d-1}$, and then given $x_1, ..., x_n$, choose $x_{n+1}$ randomly uniformly from the set
$$
\mathfrak{S} = \Big\{ y \in \S ^{d-1}: |y - x_i| >  \frac{A}{\sqrt {t}}, i =1,...,n \Big\}
$$
as long as $\mathfrak{S}$ is non-empty; if $\mathfrak{S}$ is empty, then 
$\mathscr{K}_t = \{x_1, ..., x_n \}$.
Recall the definition of $\mathcal{A}^{t,y}_\theta$ in \eqref{der Teil} on Page \pageref{der Teil}.
Set 
$$\mathcal{N}_t = \# \{ u \in \T _t: \langle X_t(u), \kappa \rangle \geq f^{t,y}_t - 1 \text{ for some } \kappa \in \mathscr{K}_t,
|X_s(u)| \leq f^{t,y}_s, s \leq t \} =
\sum \limits _{\kappa \in \mathscr{K}_t}
\#\mathcal{A}^{t,y}_\kappa $$
and let $\theta \in \S ^{\d-1}$.
By Lemma \ref{hoodwink = deceive, trick2}
for large $t$
\begin{equation} \label{shank = betw knee and ankle}
	\E \big[\mathcal{N}_t \big] =
	|\mathscr{K}_t| \E \big[\# \mathcal{A}^{t,y}_{\theta}\big]
	\geq
	C t ^{(\d-1)/2} (1+y) e^{-y} t^{-(\d-1)/2} 
	= C (1+y) e^{-y}.
\end{equation}
Next we bound the second moment of $\mathcal{N}_t$,
\begin{equation}\label{toady = smn obsequious}
	\E \big[\mathcal{N}_t ^2\big] = 
	\sum \limits _{\kappa \in \mathscr{K}_t}
	\E\big[\big( \# \mathcal{A}^{t,y}_\kappa \big) ^2\big]
	+
	\sum \limits _{\kappa, \kappa ' \in \mathscr{K}_t}
	\E \big[\#\mathcal{A}^{t,y}_\kappa \#\mathcal{A}^{t,y}_\kappa\big].
\end{equation}
By Lemma \ref{fledgling}
\begin{multline}\label{take a dig at}
	\sum \limits _{\kappa \in \mathscr{K}_t}
	\E\big[\big( \# \mathcal{A}^{t,y}_\kappa \big) ^2\big]
	= |\mathscr{K}_t| \E\big[\big( \# \mathcal{A}^{t,y}_\theta \big) ^2\big] 
	\\
	\leq C t ^{(\d-1)/2} (1+y) e^{-y}  t^{-(\d-1)/2} =  C  (1+y) e^{-y},
\end{multline}
whereas by Lemma  \ref{umbrage; take an umbrage new}
for some $c>0$
\begin{align}\label{frigid = very cold}
	\sum \limits _{\kappa, \kappa ' \in \mathscr{K}_t, \, \kappa \ne \kappa '}
	\E \#\mathcal{A}^{t,y}_\kappa \#\mathcal{A}^{t,y}_\kappa
	\leq & \,  C (1+y) e ^{-y}
	\sum \limits _{\substack{\kappa, \kappa ' \in \mathscr{K}_t, \,\kappa \ne \kappa '
			\\ 2\alpha = \arccos(\langle \kappa, \kappa '\rangle)}}
	\bigg[\frac{1}{\alpha ^{\d-2}(t+1)^{\d - 1}}
	\notag
	+ (t+1)^{-\frac{\d-1}{2} } e^{-c \alpha ^2 t } \bigg]
	\notag
	\\
	\leq & \, 
	\frac{C (1+y) e ^{-y}}{(t+1)^{\d - 1}}
	\sum \limits _{\substack{\kappa, \kappa ' \in \mathscr{K}_t,\, \kappa \ne \kappa '
			\\ 2\alpha = \arccos(\langle \kappa, \kappa '\rangle)}} \frac{1}{\alpha ^{\d-2}}
	\\ 
	+ & \, 
	\frac{ C (1+y) e ^{-y}}{(t+1)^{\d/2 - 1/2 }}
	\sum \limits _{\substack{\kappa, \kappa ' \in \mathscr{K}_t,
			\, \kappa \ne \kappa '
			\\ 2\alpha = \arccos(\langle \kappa, \kappa '\rangle)}}
	e^{-c \alpha ^2 t}. \notag
\end{align}
Now we are going to bound the two sums on the right hand side 
of \eqref{frigid = very cold}. For some $C > 0$
for all $\kappa \in \mathscr{K}_t$
there are no more than $C n ^{\d-2}$ points
$\kappa '  \in \mathscr{K}_t$ with distance  $\frac{n-1}{\sqrt{t}} \leq |\kappa - \kappa'| \leq \frac{n}{\sqrt{t}}$ from $\kappa$. 
Since $\arccos(\langle \kappa, \kappa '\rangle) \geq  |\kappa - \kappa'|$ 
we have for the first sum on the right 
hand side of \eqref{frigid = very cold}
\begin{multline}\label{obsequious}
	\sum \limits _{\substack{\kappa, \kappa ' \in \mathscr{K}_t,
			\,	\kappa \ne \kappa '
			\\ 2\alpha = \arccos(\langle \kappa, \kappa '\rangle)}} \frac{1}{\alpha ^{\d-2}} \leq 
	C \sum \limits _{\substack{\kappa, \kappa ' \in \mathscr{K}_t},
		\,	\kappa \ne \kappa '} \frac{1}{|\kappa - \kappa'| ^{\d-2}}
	= C \sum \limits _{\substack{\kappa \in \mathscr{K}_t}}  \sum \limits _{\substack{ \kappa ' \in \mathscr{K}_t, \kappa' \ne \kappa}} \frac{1}{|\kappa - \kappa'| ^{\d-2}}
	\\
	\leq C t ^{(\d-1)/2}
	\sum \limits _{n \in \N: 1\leq n \leq 2\sqrt{t} + 1}
	\bigg(\frac{\sqrt{t}}{n} \bigg) ^{\d-2} n ^{\d-2}
	=  C t ^{(2\d-3)/2} 
	\sum \limits _{n \in \N: 1\leq n \leq 2\sqrt{t} + 1} 1 
	\leq C  t^{\d - 1}.
\end{multline}
To bound the second sum on the right hand side of \eqref{frigid = very cold}
note that 
$\arccos(\langle \kappa, \kappa '\rangle) \leq C^{-1} |\kappa - \kappa'|$ for some $C > 0$, and hence 
\begin{align}
	\sum \limits _{\substack{\kappa, \kappa ' \in \mathscr{K}_t,
			\,	\kappa \ne \kappa '
			\\ 2\alpha = \arccos(\langle \kappa, \kappa '\rangle)}}
	e^{-c \alpha ^2 t} \notag 
	= & \sum \limits _{\substack{\kappa \in \mathscr{K}_t}}  \sum \limits _{\substack{ \kappa ' \in \mathscr{K}_t, \kappa' \ne \kappa,
			\\ 2\alpha = \arccos(\langle \kappa, \kappa '\rangle)}}
	e^{-c \alpha ^2 t}  \notag
	\\
	\leq & \,  C  t ^{(\d-1)/2}
	\sum \limits _{n \in \N: 1\leq n \leq 2\sqrt{t} + 1}
	\exp\Big\{-c  \Big( \frac{n}{\sqrt{t}}\Big)^2  t \Big\}
	\times n^{\d-2}
	\label{flummox = perplex, bewilder}
	\\
	= & \,   C  t ^{(\d-1)/2}
	\sum \limits _{n \in \N: 1\leq n \leq 2\sqrt{t} + 1}
	n^{\d+2}e ^{-cn} 
	\notag
	\\
	\leq & \,  C  t ^{(\d-1)/2}.
	\notag
\end{align}
Combining \eqref{frigid = very cold}, \eqref{obsequious},
and \eqref{flummox = perplex, bewilder}
we get 
\begin{equation}
	\sum \limits _{\kappa, \kappa ' \in \mathscr{K}_t}
	\E \big[\#\mathcal{A}^{t,y}_\kappa \#\mathcal{A}^{t,y}_\kappa\big]
	\leq C (1+y) e ^{-y}.
\end{equation}
Hence by \eqref{toady = smn obsequious} and \eqref{take a dig at} 
\begin{equation}
	\E \big[\mathcal{N}_t ^2 \big]\leq C (1+y) e ^{-y}.
\end{equation}
By the Cauchy--Schwarz inequality
$\E\big[\mathcal{N}_t \big] ^2 \leq \P \{ \mathcal{N}_t \geq 1\} \E \big[\mathcal{N}_t ^2\big]$
and therefore by \eqref{shank = betw knee and ankle}
for large $t$
\begin{equation}
	\P \{ \mathcal{N}_t \geq 1   \}
	\geq \frac{\E\big[\mathcal{N}_t \big] ^2}{\E \big[\mathcal{N}_t ^2\big]}
	\geq \frac{(1+y) e^{-y}}{C}.
\end{equation}
It remains to recall that 
$f_t ^{t,y} = \frac{\Psi (\lambda)}{\lambda}t +  \frac{\d-4}{2\lambda} \ln (t+1) + \frac {3M}{\lambda}  
+ \frac y \lambda$, and hence
for $t \in \N$
\begin{equation}
	\P \left\{ \exists u \in \T _t:
	|X_t(u)| \geq r_t + y \right\}
	\geq  \P \{ \mathcal{N}_t \geq 1   \}.
\end{equation}
\end{proof}

\begin{proof}
[Proof of Proposition \ref{prop fartherst point lower bound}]
\label{Proof prop fartherst point lower bound}
Let $\eta _t = \{ X_t(u), u \in \T_t\}$
be the set of spatial locations of the particles
at time $t \in \N$.
Further, denote by $X^R _t$ the spatial location of the farthest point from the origin at time $t$, so that $|X^R _t| = R_t$;
if the branching random walk has died out by the time $t$, we set $X^R _t = 0 _\d$. Should there be multiple farthest points, one of them is chosen uniformly at random, so that the direction 
of $X^R _t$  given $|X^R _t| >0$ is uniformly distributed on the unit sphere. 
Let us introduce an independent copy
$(\widetilde \T, \{\widetilde X_t(u), t \in \Z_+, u \in  \widetilde\T \})$
of the BRW $(\T, \{X_t(u), t \in \Z_+, u \in \T \})$
defined on a different probability space with the probability measure $\widetilde  \P $;
we write $\widetilde R_t, \widetilde X^ R_t$ for the corresponding quantities.
Recall  that we follow the standard convention that the product over the elements of the empty set equals $1$.
Conditioning on $\eta _h$, $h \in \N$,
we find for large $t$
\begin{align*}
	\P \{R_{t+h} \leq & \, r_t  | \eta _h   \} = \E \Big[  \prod\limits _{x \in \eta _h}  
	\widetilde \P\{\forall u \in \widetilde \T, |u | = t: |\widetilde X_t(u) + x| \leq r_t   \} \Big| \eta _h \Big]
	\\
	\leq & \, \E \Big[ \prod\limits _{x \in \eta _h}  
	\widetilde \P\big\{ |\widetilde X^R _t + x| \leq r_t  \big\} \Big| \eta _h \Big]
	\leq  \E \Big[ \prod\limits _{x \in \eta _h}  
	\big[1-
	\widetilde \P\big\{ |\widetilde X^R _t + x| \geq r_t  \big\}  
	\big] \Big| \eta _h \Big]
	\\
	\leq &  \, \E \Big[ \prod\limits _{x \in \eta _h}  
	\big[1- 
	\frac 12
	\widetilde \P\big\{ |\widetilde R_t| \geq r_t  \big\}  
	\big] \Big| \eta _h \Big]
	= \E \Big[1- 
	\frac 12
	\widetilde \P\big\{ |\widetilde R_t| \geq r_t  \big\}  
	\Big| \eta _h \Big] ^{|\eta _h|}
	\leq  \E
	\Big[ (1- 
	a
	) ^{|\eta _h|} \Big| \eta _h \Big],
\end{align*}
where $a : = \frac 13\liminf\limits_{t \in \N}\P\big\{ |R_t| \geq r_t \big\} > 0$ by Proposition \ref{volcano plume}.
Recall that
$\mathcal{S}$ is the event that 
the BRW survives. Since $\mcs \subset \{|\eta _h| \geq 1\}$ we obtain 
by conditioning on $\eta _h $
\begin{align*}
	\P\Big( \{R_{t+h} \leq & \,  r_t \} \cap \mcs   \Big)
	\leq  \P\Big( \{R_{t+h} \leq   r_t \} \cap \{ |\eta _h| \geq 1\}   \Big)
	\\
	= &  \, \E\left[\E \big[ \1 \{R_{t+h} \leq   r_t \} \1 \{ |\eta _h| \geq 1\} \big| \eta _h \big]\right]
	= \E \big[\1 \{ |\eta _h| \geq 1\} \E \big[  \1 \{R_{t+h} \leq   r_t \}  \big| \eta _h \big] \big]
	\\
	= &  \,
	\E \Big[\1 \{ |\eta _h| \geq 1\} \P \big\{R_{t+h} \leq   r_t  \big| \eta _h \big\} \Big] \leq 
	\E \big[\1 \{ |\eta _h| \geq 1\}  (1- 
	a
	) ^{|\eta _h|} \big].
\end{align*} 
The process $ \{|\eta _h|, h \in \Z_+  \}$
is a supercritical Galton-Watson process started from $1$.
Therefore for any  $n \in \N$, $\P \big\{|\eta _h| \in [1,n] \big\}\xrightarrow{h \to \infty} 0$, 
and consequently
\begin{equation}\label{take a red eye}
	\rho_h : = 
	\E \big[\1 \{ |\eta _h| \geq 1\}  (1- 
	a
	) ^{|\eta _h|} \big] \xrightarrow{h \to \infty} 0.
\end{equation}
Since $r_{t+h} - r_t \leq c h$
for some $c > 0$,
we get for large $t \in \N$
for all $h \in \N $
$$\P \big\{R_{t+h} \leq  r_{t+h} - c h | \mathcal{S} \big\} \leq
\frac{\P\big( \{R_{t+h} \leq   r_t \} \cap \mcs   \big)}{\P(\mcs)} \leq \frac{\rho_h}{\P(\mcs)}, $$
hence for large $s \in \N$
and $h \in \N$ satisfying $h \leq \frac s 2$,
$$\P \Big\{R_{s} \leq  r_s - c h \big| \mathcal{S} \Big\} \leq 
\frac{\rho_h}{\P(\mcs)}. $$
The statement 
\begin{equation*}
	\sup\limits_{t \geq 0} \P\big\{R_t - r_t  \leq -y  \big| \mathcal{S} \big\} \to 0  \  \mbox { as } \ y \to \infty\, 
\end{equation*}
now follows from \eqref{take a red eye}.
\end{proof}

\section{Ballot theorem with moving barrier}\label{sec ballot}
In this section, $(S_n, n \in \Z_+)$ is a one-dimensional centered non-lattice non-degenerate
random walk
with finite $(3 + \varepsilon)$-th moment, i.e. $0 < \E \big[|S_1|^{3 + \varepsilon}\big] < \infty$ for some 
$\varepsilon > 0$. We always assume $S_0=0$.
We   prove 
Theorem \ref{thm homing device = guidance system}
and Theorem \ref{thm beinhalten = include} in this section.
For $a \geq 0$ denote  $$T_{-a} = \min\{n \in \N: S_n < -a\}.$$
The rough proof idea of Theorem \ref{thm homing device = guidance system} is to show that the moment when $(S_n, n \in \Z_+)$ crosses  the barriers $f(n)-a$ or $-f(n)-a$ is unlikely to  drastically differ from  $T_{-a}$. The proof of Theorem \ref{thm beinhalten = include} follows a `gluing' argument with the trajectory of the random walk being split into three parts. This argument is well known (\cite{Ballot_unpublished}). Over the next couple of pages we collect preliminary results involving the crossing times and overshoots of the random walk. After that we proceed with  a few auxiliary lemmas. The proof of Theorem \ref{thm homing device = guidance system}
then follows starting on Page
\pageref{Fuelle an, von}. 
The proof of Theorem \ref{thm beinhalten = include} is located on Page \pageref{Fuelle = abundance, wealth}.
\\
  
The following lemma from analysis is given without proof.
\begin{lem}\label{whittle}
	For an increasing  $f: \N \to \R_+$, 
	the series 
	$$
	\sum\limits _{n \in \N } \frac{f(n)}{n^{3/2}}
	$$ converges if and only if
	$$
	\sum\limits _{m \in \N } \frac{f(2^m)}{2^{m/2}} < \infty.
	$$
\end{lem}
  
We now collect a few results from the literature 
on random walks.

\begin{lem}[{{\cite[Theorem A]{Kozlov76} or \cite[Lemma 3.3]{PP95}}}]
	\label{chow = food}
	There exists $C> 1$ such that
	for $a \in [0, n^{1/2} ]$
	\begin{equation}
		\frac{(a+1)}{C\sqrt{n}}
		\leq
		\P\{  T_{-a} > n \} \leq  \frac{C(a+1)}{\sqrt{n}}.
	\end{equation}
\end{lem}
  
The following statement is a corollary 
to \cite[Theorem 1]{Rog62}, 
and is basically 
\cite[Theorem 1]{Ballot_unpublished}; see also \cite{Stone67}.

\begin{lem}
	\label{countenance = facial expression}
	There exists $c > 0$ such that
	for all $x \in \R$ and $n \in \N$
	\begin{equation*}
		\P \{ x \leq S_n \leq x+1  \}
		\leq \frac{c}{\sqrt{n}}.
	\end{equation*}
\end{lem}
  
Let $\{(F _n, H_n), n \in \Z_+  \}$
be the strict ascending ladder variables 
for the random walk $(S_n, n \in \Z_+)$, that is,
$F _0 = H_0 = 0$, 
and for $n \in \N$

\begin{equation*}
	F_n =  \inf\{k \in \N: 
	S_{F _{n-1} + k} >  H _{n-1} \}, 
	\ \ \ 
	H _n =  S_{F _n} - H_{n-1}.
\end{equation*}
  
The sequence of the ladder heights $\{H _n, n \in \N \}$ is a sequence of strictly positive i.i.d random variables
(more information on ladder variables can be found 
in \cite[Chapter 12 and elsewhere]{FellerVol2}, \cite[Chapter 8]{Chung}).
The following lemma relates moments 
of $ S_1 \vee 0$ and $H_1$.

\begin{lem}
	[{\cite[Theorem, (i)]{Ladder_heights_moments}}]
	\label{verdauen = digest}
	Let $\phi: \R_+ \to \R_+$ be an increasing, differentiable function with $\lim\limits _{x \to \infty} \phi (x) = \infty$ and for every $K > 0$
	$$ \limsup\limits _{x \to \infty} \frac{\phi(x+K)}{\phi(x)} < \infty.$$
	Set $\Phi (y) = \int\limits _0 ^y \phi(x)dt $, $y  \geq 0$.
	If $\E \big[ \Phi ( S_1 \vee 0 ) \big] < \infty$, 
	then $\E \big[ \phi ( H_1   ) \big] < \infty$.
\end{lem}
  
For $b \geq 0$
set $T ^+ _b = \min\{n \in \N: S_n \geq b\}$.
We will need bounds involving 
the   `overshoot'
random variables $ S_{T ^+ _b} - b$, $b\geq 0$.
It follows from the construction of the ladder heights
that the overshoots of the ladder heights 
$\{H _n, n \in \N \}$
coincide with those of $(S_n, n \in \Z_+)$:
\begin{equation}\label{have someone eating out of your hand}
	S_{T ^+ _b} - b = 
	\min\{H_1 + \dots + H _n - b: n \in \N, H_1 + \dots + H _n \geq b \}, 
	\ \ \ b \geq 0.
\end{equation}
Let a random variable $L$ have the size-biased distribution of $H_1$, 
that is, the respective distribution functions satisfy
\begin{equation}\label{eitel = vain}
	F_L(dx) = \frac{x F _{H_1} (dx)}{\E[H_1]}
\end{equation}
The following lemma is a direct corollary 
to \cite[Proposition 4.3]{Ineq_for_overshoot94}.
\begin{lem}\label{minnow = a small fish}
	For $p \geq 0$ and $b \geq 0$
	\begin{equation}
		\E \Big[ \big( S_{T ^+ _b} - b\big) ^p \Big]
		\leq \E \big[ L ^p \big] + \big( \E [ L ] \big) ^p.
	\end{equation}
\end{lem}

\subsection{Proof of Theorem \ref{thm homing device = guidance system}}

The next lemma shows that the `overshoot'
random variables $ S_{T ^+ _b} - b$
have a finite $(1+\varepsilon)$-th moment uniformly in $b$.
\begin{lem} \label{zuneigung = affection} 
	For some $C > 0$ for all $b \geq 0$ 
	\begin{equation}
		\E \Big[ \big( S_{T ^+ _b} - b\big) ^{1+\varepsilon} \Big] \leq C. 
	\end{equation}
\end{lem}
  
\begin{proof}
Since $\E \big[|S_1|^{3 + \varepsilon}\big] < \infty $
we can apply
Lemma \ref{verdauen = digest}
with $\phi(x) = x^{2 + \varepsilon}$
to find 
$$
\E \big[ H_1 ^{2 + \varepsilon} \big] < \infty.
$$
Hence by   \eqref{eitel = vain}
$$
\E \big[ L ^{1 + \varepsilon} \big]
= \int\limits _{x \geq 0} x^{1 + \varepsilon} F_L(dx)
= \frac{1}{\E[H_1]} \int\limits _{x \geq 0} x^{2 + \varepsilon} F_{H_1}(dx) = \frac{\E \big[ H_1 ^{2 + \varepsilon} \big]}{\E[H_1]} < \infty.
$$
Applying
Lemma \ref{minnow = a small fish} 
with $p = 1+\varepsilon$ gives the desired result.
\end{proof}

	For a function $g: \Z_+ \to \R$ and $k,n \in \Z_+$, $k \leq n$, 
	denote by $A(g;k,n)$ the event $\{ S_i \geq g(i) \text{ for all } k \leq i \leq n \}$, and let $A(g;n)$ be a shorthand for $A(g;0,n)$, \label{letzten Endes = ultimately} as in \cite{PP95}.
	We also write $A(-a; n)$ for $\{T_{ -a} > n\}$.
	The following Lemma is a uniform version of a limit law
		for a random walk conditioned to stay positive, see \cite{CC08}.
	\begin{lem} \label{der Zettel}
		For some $c> 0$ for $n$ large enough
		and $a \in [0,n^{1/2}]$
		\begin{equation}
			\P \big\{  A(-a; n), 
			S_{ n } \in  [\sqrt{3n}, 2\sqrt{3n}] \big\}
			\geq 
			\frac{ c (a+1)}{\sqrt{n}}.
		\end{equation}
	\end{lem}
	  
	\begin{proof}
	By Lemma \ref{chow = food} for some $c_1> 0$
	\begin{equation}\label{prowl copied2}
		\PP{  A(-a; n)  }\geq
		\frac{ c_1(a+1)}{\sqrt{n}}.
	\end{equation}
	Recall that
	for $b \geq 0$,
	$T ^+ _b = \min\{n \in \N: S_n \geq b\}$.
	Set $ \mathfrak{q} = \frac{\sqrt{3n}}{3}$.
	By the central limit theorem for some $c_2> 0$
	for large  $n \in \N$
	\begin{equation}\label{nachsagen = repeat copied2}
		\PP{ S_{ n} \geq 4\mathfrak{q}  } \geq c_2
	\end{equation}
	and hence also
	\begin{equation}\label{nachsagen = repeat2 copied2}
		\PP{ T ^+_{  4\mathfrak{q}} \leq n  } \geq c_2
	\end{equation}
	The events
	$  A(- a; n)$ and $\{ T ^+_{4\mathfrak{q}} \leq n\}$
	are increasing events determined by the trajectory of the 
	random walk $(S_k)_{1\leq k \leq n}$. Therefore they are positively correlated
	by the FKG inequality (see e.g. \cite[Proposition 2.7]{legrand2024fkginequalitiesstochasticprocesses}), and by \eqref{prowl copied2}
	and \eqref{nachsagen = repeat2 copied2} we get for 
	${ B : = A( - a; n)\cap \big\{ T ^+_{ 4\mathfrak{q}} \leq  n \big\}}$
	for large $n$
	
	\begin{equation}\label{zuversictlich = optimistic, confident copied2}
		\P (B)  \geq \frac{ c_1c_2 (a+1) }{ \sqrt{n} }.
	\end{equation}
	Set $  B _1 = \big\{ T ^+_{ 4\mathfrak{q}} < T  _{ -a},
	T ^+_{ 4\mathfrak{q}} \leq n     \big\}$.
	Then $B \subset  B _1$
	and hence by \eqref{zuversictlich = optimistic, confident copied2} 
	\begin{equation}\label{Zaehler = numerator}
		\P(B_1)  \geq \frac{ c_1c_2 (a+1) }{ \sqrt{n} }.
	\end{equation}
	Lemma \ref{zuneigung = affection}
	implies that given $ B_1$ the overshoot 
	$S_{T ^+_{4\mathfrak{q}}} - 4\mathfrak{q}$ is typically  small
	compared to $\sqrt{n}$: 
	there is $C> 0$ such that
	\begin{multline}\label{umsichtig = prudent copied2}
		\P \Big\{ S_{T ^+_{4\mathfrak{q}}} - 4\mathfrak{q} \geq \mathfrak{q} \Big|  B_1 \Big\}
		\leq \frac{\P \big\{ S_{  T ^+_{ 4\mathfrak{q}}} - 4\mathfrak{q} \geq \mathfrak{q}  \big\} }{\P (  B_1)}
		\leq \frac{\E\big[  (S_{  T ^+_{ 4\mathfrak{q}}} - 4\mathfrak{q})^{1+\varepsilon}\big] }{\P (  B_1)  \mathfrak{q} ^{1+\varepsilon}}
		\leq C n^{-\varepsilon/2}, \ \ \ n \in \N,
	\end{multline}
	and hence for large 
	$n$
	\begin{equation}\label{Nenner = denominator}
		\P \Big\{ S_{T ^+_{4\mathfrak{q}}} - 4\mathfrak{q} \leq \mathfrak{q} \Big|  B_1 \Big\}
		\geq  \frac 23.
	\end{equation}
	Next define the event
	\begin{equation}
		B_2 =  \big\{  T ^+_{ 4\mathfrak{q}} \leq n, |S_k - S_{ T ^+_{4\mathfrak{q}} }| \leq 
		\mathfrak{q}
		\text{ for }   T ^+_{ 4\mathfrak{q}} + 1 \leq k \leq n \big\} 
	\end{equation} 
	For large $n$, the conditional probabilities of $B_2$, given 
	${  B_1 \cap \{  T ^+_{ 4\mathfrak{q}}  = i\}}$,
	$1 \leq i \leq n$, are bounded away from $0$: 
	for some $c_{19} > 0$ 
	\begin{equation}	\label{harken back to copied2}
		\inf _{i = 1,...,n}\P \big\{ B_2 | B_1,   T ^+_{ 4\mathfrak{q}}  = i\big\} \geq c_{19}.	
	\end{equation}
	Note that 
	$$B_1 \cap \big\{ S_{T ^+_{4\mathfrak{q}}} - 4\mathfrak{q} \leq \mathfrak{q}\big\} \cap B_2 \subset   \big\{  A(-a; n), 
	S_{ n } \in  [\sqrt{3n}, 2\sqrt{3n}] \big\}.$$
	Given $B_1 \cap \big \{ 
	T ^+_{ 4\mathfrak{q}}  = i \big\} $, the events 
	$B_2$ and $\big\{ S_{T ^+_{4\mathfrak{q}}} - 4\mathfrak{q} \leq \mathfrak{q} \big\}$
	are conditionally independent.
	In combination with 
	\eqref{Nenner = denominator}
	and \eqref{harken back to copied2}
	this gives 
	for large $n$ 
	\begin{align}
		\P\big\{  A(-a; & \, n), 
		S_{ n } \in  [\sqrt{3n}, 2\sqrt{3n}] \big\}
		\geq   \P\big\{ B_2, S_{T ^+_{4\mathfrak{q}}} - 4\mathfrak{q} 
		\leq  \frac{\sqrt{n}}{3},  B_1 \big\}
		\\
		= & \,
		\sum\limits _{i=1} ^{n} 
		\P\big\{ B_2, S_{T ^+_{4\mathfrak{q}}} - 4\mathfrak{q} 
		\leq  \mathfrak{q} \big| B_1,  T ^+_{ 4\mathfrak{q}}  = i \big\} 
		\P\big\{  T ^+_{ 4\mathfrak{q}}  = i, B_1  \big\} 
		\\ \notag
		= & \,
		\sum\limits _{i=1} ^{n} 
		\P\big\{ B_2 \big| B_1,  T ^+_{ 4\mathfrak{q}}  = i \big\} 
		\P\big\{  S_{T ^+_{4\mathfrak{q}}} - 4\mathfrak{q} 
		\leq  \mathfrak{q} \big| B_1,  T ^+_{ 4\mathfrak{q}}  = i \big\}  
		\P\big\{  T ^+_{ 4\mathfrak{q}}  = i, B_1  \big\} 
		\\
		\geq & \,
		\sum\limits _{i=1} ^{n} 
		c_{19}
		\P\big\{  S_{T ^+_{4\mathfrak{q}}} - 4\mathfrak{q} 
		\leq  \mathfrak{q} \big| B_1,  T ^+_{ 4\mathfrak{q}}  = i \big\} 
		\P\big\{  T ^+_{ 4\mathfrak{q}}  = i, B_1  \big\} 
		\\
		\label{der Bias}
		= & \,
		c_{19} 
		\P\big\{  S_{T ^+_{4\mathfrak{q}}} - 4\mathfrak{q} 
		\leq  \mathfrak{q}, B_1 \big\} 
		\geq \frac{2c_{19}}{3}  \P (B_1).
	\end{align}
	The statement of the lemma now follows 
	from \eqref{der Bias} and  \eqref{Zaehler = numerator}.
	\end{proof}
	
	 The following lemma  is a corollary
	to \cite[Lemma 3.9]{Mallein_interfaces}. It 
	is used in the proof of Theorem \ref{thm homing device = guidance system}.
	
	\begin{lem}\label{entnehmen}
		For some $c > 0$
		for $a, b \in [0, \sqrt{n}]$
		\begin{equation}\label{verlegen auf = reschedule}
			\P \{ A(-a; n), S_n  \in [-a +b,-a+b+1)    \} \geq  \frac{c(a+1) (b+1)}{n^{3/2}}.
		\end{equation}
	\end{lem}

	The following auxiliary lemma is used in the proof of  Theorem \ref{thm homing device = guidance system}.
	\begin{lem} \label{perusal = examine sth}
		Let $f$ be as in Theorem \ref{thm beinhalten = include}. For some $C > 0$
		for $a \in [0, \sqrt{n}]$
		for large $n \in \N$ and $N \geq n$, $N \in \N$
		\begin{equation}\label{in sich gehen}
			\P\{A(0;n), A(-f; N) \} \leq \frac{C}{a+1} \P \{A(-a;n), A(-f-a; N)  \}.
		\end{equation}
	\end{lem}
	  
	\begin{proof} 
	We first show that for some 
	$c > 0$
	\begin{equation}\label{solder = payat}
		\P\{A(0;n), A(-f; N), S_n \in [0, 2\sqrt{n}] \}
		\geq c \P\{A(0;n), A(-f; N) \}, \ \ \ n \in \N.
	\end{equation}
	Without loss of generality and to simplify notation we assume in this proof that the variance $\E [S_1 ^2] = 1$.
	Later  we will use a limit theorem
	for a random walk 
	conditioned to stay non-negative (\cite{Igl74}, or e.g. \cite{LLT_RW_cond}): for $k \in \Z_+$
	\begin{equation}\label{Optionen abwaegen}
		\lim\limits _{n \to \infty}\P\big\{S_n \in [kn^{1/2}, (k+1) n^{1/2}] | A(0;n)  \big\} = \int\limits _{k}^{k+1} x e^{-x^2/2}dx.
	\end{equation}
	For $n,k \in \N $
	consider a stopping time $ \tau ^{(k)} = \min\{ i \in \N:
	S_i \geq kn^{1/2} \text{ or }  S_i \leq -n^{1/2} \}$. 
	By the optional stopping theorem 
	$\E \big[S_{\tau ^{(k)}}\big]  =0$, 
	and with Lemma \ref{zuneigung = affection} 
	we can control the overshoot 
	above $k n^{1/2}$, 
	that is, 
	\begin{equation}\label{schielen = squint, peer}
		\E\big[ (S_{\tau ^{(k)}} - kn^{1/2} ) 
		\1 \{S_{\tau ^{(k)}} \geq kn^{1/2} \}\big]
		\leq \E \big[(S_{T ^+ _ {kn^{1/2}}} - kn^{1/2} )\big]  \leq C,
	\end{equation}
	as well as  the undershoot below  $-n^{1/2}$,
	\begin{equation}\label{das Blei = lead}
		\E\big[ |S_{\tau ^{(k)}} + n^{1/2} | \1 \{S_{\tau ^{(k)}} \leq  -n^{1/2} \}\big]
		\leq \E \big[(-S)_{T^+_ {n^{1/2}}} - n^{1/2} \big] \leq C.
	\end{equation}
	In combination with 
	\begin{multline*}
		0 = \E \big[S_{\tau ^{(k)}}\big]  = \E \big[(S_{\tau ^{(k)}} - kn^{1/2} ) 
		\1 \{S_{\tau ^{(k)}} \geq kn^{1/2} \} \big]
		+  kn^{1/2} \P \{S_{\tau ^{(k)}} \geq kn^{1/2} \}
		\\
		- \E\big[ |S_{\tau ^{(k)}} + n^{1/2} | \1 \{S_{\tau ^{(k)}} \leq  -n^{1/2} \}\big]
		- n^{1/2} \P \{S_{\tau ^{(k)}} \leq -n^{1/2} \} 
	\end{multline*} 
	\eqref{schielen = squint, peer} 
	and \eqref{das Blei = lead}
	imply 
	\begin{equation}\label{Einhalt gebieten = arrest, stem, prevent}
		\P\{S_{\tau ^{(k)}} \geq  kn^{1/2} \}
		\geq \frac{c}{k},  \ \ \  k \in \N,
	\end{equation}
	for some $c > 0$.
	Now define $ \sigma  ^{(k)}
	= 
	\min\{ i \in \N, i > n:
	S_i \geq kn^{1/2} \text{ or }  S_i \leq 0 \}$.
	By the strong Markov property (or simply because 
	of the independence of the increments of the random walk)
	and \eqref{Einhalt gebieten = arrest, stem, prevent}
	we get 
	\begin{equation}\label{Verdacht nahe legen}
		\P \big\{ S_{\sigma  ^{(k)}} \geq k n^{1/2} 
		|  A(0; n),S_n \in [n^{1/2}, 2 n^{1/2}] \big\}
		\geq \frac{c}{k},  \ \ \ k \in \N.
	\end{equation}
	To put it into words, given $S_n \in [n^{1/2}, 2 n^{1/2}]$,
	the probability for $(S_i)$
	to go above $k n^{1/2} $
	before falling below $0$ is at least $\frac{c}{k}$.
		\begin{figure}[t!]
		\vspace{-2.1cm}
		\centering
		\includegraphics[scale=0.65]{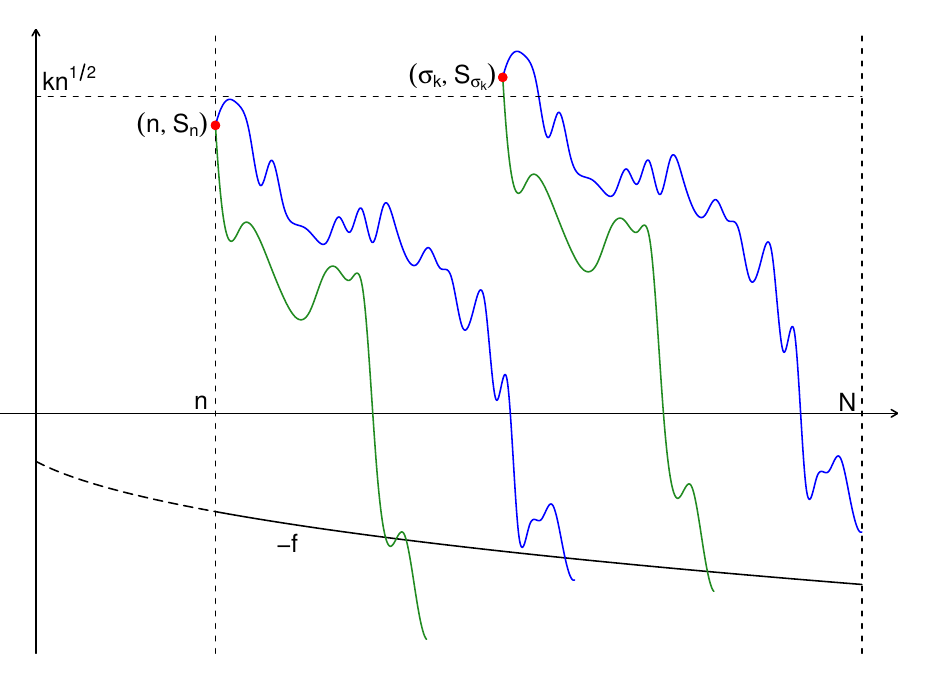}
		\vspace{-0.4cm}
		\caption{ \small
			An illustration to \eqref{hervorrufen = cause, trigger}. 
  Given $\{ A(0; n),S_n \in [n^{1/2}, 2 n^{1/2}],
  S_{\sigma  ^{(k)}} \geq k n^{1/2} \}$, the trajectory started at $(\sigma_k, S_{\sigma_k})$ has to come down a long way to cross the barrier $-f$. Every trajectory that manages that would also cross $-f$ if it started at $(n, S_n)$ - an example is given by the green trajectory. The reverse is not true however, as demonstrated by the blue trajectory. 
		}
		\label{f_for_proof3}
		\vspace{-0.31cm}
	\end{figure}
	Note that 
	\begin{multline}\label{hervorrufen = cause, trigger}
		\P \Big\{  A(-f; n+1, N) \big| A(0; n),S_n \in [n^{1/2}, 2 n^{1/2}],
		S_{\sigma  ^{(k)}} \geq k n^{1/2} \Big\}
		\\
		\geq 
		\P \Big\{  A(-f; n+1, N) \big| A(0; n),S_n \in [(k-1)n^{1/2}, k n^{1/2}]
		\Big\}
	\end{multline}
	because given $\{ A(0; n),S_n \in [n^{1/2}, 2 n^{1/2}],
	S_{\sigma  ^{(k)}} \geq k n^{1/2} \}$
	the random walk has to go further down (from $S_{\sigma  ^{(k)}} $
	to below the barrier $k \mapsto -f(k)$)
	in a fewer steps in order for $A(-f; n+1, N)$
	not to occur, as compared to when conditioning on 
	$ \{A(0; n),S_n \in [(k-1)n^{1/2}, k n^{1/2}] \}$.
	An illustration can be found in Figure \ref{f_for_proof3}.
	It follows from \eqref{Verdacht nahe legen} 
	and \eqref{hervorrufen = cause, trigger}
	that for large $n \in \N$
	\begin{align}
		\P \Big\{  A(-f; n+1, N)  \big| & A(0; n),S_n \in [n^{1/2}, 2 n^{1/2}] \Big\} \notag 
		\\
		\notag 
		\geq \, &
		\P \Big\{  A(-f; n+1, N) \cap \{ S_{\sigma  ^{(k)}} \geq k n^{1/2}  \} \big| A(0; n),S_n \in [n^{1/2}, 2 n^{1/2}] \Big\}
		\\
		\notag 
		= \, &
		\begin{multlined}[t]
			\P \Big\{  A(-f; n+1, N) \big| A(0; n),S_n \in [n^{1/2}, 2 n^{1/2}],
			S_{\sigma  ^{(k)}} \geq k n^{1/2} \Big\} 
			\\
			\times
			\P \big\{ S_{\sigma  ^{(k)}} \geq k n^{1/2} 
			|  A(0; n),S_n \in [n^{1/2}, 2 n^{1/2}] \big\}
		\end{multlined}
		\\
		\geq \, & \frac ck 
		\P \Big\{  A(-f; n+1, N) \big| A(0; n),S_n \in [(k-1)n^{1/2}, k n^{1/2}]
		\Big\}.
		\label{fortschreiten = proceed, progress}
	\end{align}
	Now we can finally  
	obtain 
	\eqref{solder = payat}
	by 
	conditioning, decomposing according to the values
	of $S_n$, 
	and applying
	\eqref{Optionen abwaegen} 
	and \eqref{fortschreiten = proceed, progress}:
	for large $n \in \N$
	
	\begin{equation*}
		\frac
		{\P\{A(0;n), A(-f; N) \} }
		{\P\{A(0;n), A(-f; N), S_n \in [n ^{1/2},  2n ^{1/2}] \}} 
		\hspace{7.0cm}
	\end{equation*}
	\vspace{-0.2cm}
	\begin{align*}
		& = 
		\frac
		{\sum\limits _{k \in \N }\P\{A(0;n), A(-f; N),
			S_n \in [(k-1)n ^{1/2}, kn ^{1/2}] \} }
		{\P\{A(0;n), A(-f; N), S_n \in [n ^{1/2}, 2n ^{1/2}] \}}
		\\
		& = 
		\begin{aligned}[t]
			\sum\limits _{k \in \N } \Bigg[ &\frac
			{\P\big\{ A(-f; n+1, N) \big| S_n \in [(k-1)n ^{1/2}, kn ^{1/2}], 
				A(0;n) \big\} }
			{\P\big\{ A(-f; n+1, N) \big| S_n \in [n ^{1/2}, 2n ^{1/2}], 
				A(0;n) \big\} }
			\\
			\times &
			\frac
			{ \P \big\{ S_n \in [(k-1)n ^{1/2}, kn ^{1/2}]
				\big|A(0;n)  \big\}  \P \big\{ A(0;n) \big\}}
			{ \P \big\{ S_n \in [n ^{1/2}, 2n ^{1/2}]
				\big|A(0;n)  \big\}  \P \big\{ A(0;n) \big\} }
			\Bigg]
		\end{aligned}
		\\
		&
		\overset{\substack{
				\eqref{Optionen abwaegen}
				,
				\eqref{fortschreiten = proceed, progress}}}{\leq } 
		C 
		\sum\limits _{k \in \N } k \times k e^{-k^2/2},
	\end{align*}
	and \eqref{solder = payat} follows. \\
	  
	We now proceed with the proof of \eqref{in sich gehen}.
	The idea is to
	decompose over the values of $S_n$ and prove an inequality similar to \eqref{in sich gehen} but with the restriction  $S_n \in [m, m+1)$ (see \eqref{lithe = thin, supple,young, healthy}), and then apply \eqref{solder = payat}.
	By the aforementioned limit theorem
	for a conditional random walk
	\cite{LLT_RW_cond}
	and Lemma \ref{entnehmen}
	for $m \in [0, 2 \sqrt{n}]$ there exists 
	$C > 0$ such that 
	for all $a \in [0, \sqrt{n}]$
	and $n \in \N $
	\begin{align*}
		\P \{ A(0; n), S_n  \in [m,m+1)    \} \leq & \frac{C (m+1)}{n^{3/2}},
		\\
		\P \{ A(-a; n), S_n  \in [m,m+1)    \} \geq & \frac{(a+1) (m+1)}{Cn^{3/2}}.
	\end{align*}
	\nopagebreak
	  
	(Note that we can replace $\sqrt{n}$ by $2 \sqrt{n}$
	in Lemma \ref{entnehmen}
	for example by rescaling the random walk.)
	Consequently for some $c > 0 $
	\begin{equation}\label{jaded = fatigued}
		\P \{ A(0; n), S_n  \in [m,m+1)    \} \leq 
		\frac{c}{a+1}\P \{ A(-a; n), S_n  \in [m,m+1)    \}, 
		\ \ \ a \in [0, \sqrt{n}], n \in \N.
	\end{equation}
	Next we note that
	for some $C>0$
	for large $n$
	and $a \geq 0$,  $m \in [0, 2 \sqrt{n}]$
	a.s.\ on $ \{ S_n \in [m, m+1) \}$
	\begin{multline} \label{equinox}
		\P\big\{ A(-f;n+1, N) | S_n   \big\} 
		\leq \P \{S_i + m + 1 \geq -f(i), i = n+1,...,N  \}
		\\
		\leq C \P\big\{ A(-f-a;n+1, N) | S_n   \big\}.
	\end{multline}
	\nopagebreak
	(the probability  in the middle of \eqref{equinox} equals $\P \{ A(-f;n+1, N) | S_n = m+1   \} $)
	and hence by conditional independence for  $ m \in [0, 2 \sqrt{n}]$
	\begin{align}
		\P\{A(0;n), & \, A(-f;n+1,  N), S_n \in [m, m+1) \} 
		\\
		= & \,
		\E\Big[  \P\big\{A(0;n), A(-f;n+1, N) | S_n  \big\}
		\1\{ S_n  \in [m,m+1) \} \Big]
		\\
		= & \,
		\E\Big[  \P\big\{A(0;n) | S_n  \big\}
		\P\big\{ A(-f;n+1, N) | S_n   \big\} \1\{ S_n  \in [m,m+1) \} 
		\Big]
		\\
		\leq  & \,
		\P \{S_i + m + 1 \geq -f(i), i = n+1,...,N  \}
		\E\Big[  \P\big\{A(0;n) | S_n  \big\}
		\1\{ S_n  \in [m,m+1) \} 
		\Big]
		\\
		\leq  & \,
		\P \{S_i + m + 1 \geq -f(i), i = n+1,...,N  \}
		\P\big\{A(0;n),  S_n  \in [m,m+1) \big\}.
	\end{align}
	We continue by applying \eqref{jaded = fatigued} 
	\begin{align}
		\notag
		\P \{S_i + m + 1 \geq & \, -f(i), i = n+1,...,N  \}
		\P\big\{A(0;n),  S_n  \in [m,m+1) \big\}
		\\ \notag
		\leq  & \,
		\frac{c}{a+1} \P \{S_i + m + 1 \geq  -f(i), i = n+1,...,N  \}
		\P\big\{A(-a;n),  S_n  \in [m,m+1) \big\}
		\\ \notag
		=  & \,
		\frac{c}{a+1}  \P \{S_i + m + 1 \geq -f(i), i = n+1,...,N  \}
		\\  & \times
		\E\Big[  \P\big\{A(-a;n) | S_n  \big\}
		\1\{ S_n  \in [m,m+1) \} 
		\Big]
		\\ \notag
		\overset{\eqref{equinox}}{\leq}  & \,\frac{C}{a+1} 
		\E\Big[  \P\big\{A(-a;n) | S_n  \big\}
		\P\big\{ A(-f - a;n+1, N) | S_n   \big\} \1\{ S_n  \in [m,m+1) \} 
		\Big]
		\\ \notag
		=  & \,
		\frac{C}{a+1} 
		\P\{A(-a;n),   A(-f -a ;n+1,  N), S_n \in [m, m+1) \}. 
	\end{align}
	Thus for large $n$
	and $a \in [0, \sqrt{n}]$, $m \in [0, 2 \sqrt{n}]$
	\begin{multline}\label{lithe = thin, supple,young, healthy}
		\P\{A(0;n),  A(-f;n+1,  N), S_n \in [m, m+1) \} 
		\\
		\leq 
		\frac{C}{a+1} 
		\P\{A(-a;n),   A(-f-a;n+1,  N), S_n \in [m, m+1) \}.
	\end{multline}
	To simplify notation in the remaining part of the proof we treat $\sqrt{n}$
	as if it was an integer.
	By \eqref{solder = payat} and \eqref{lithe = thin, supple,young, healthy}
	\begin{align*}
		\P\big\{A(0;n), & \,  A(-f; n+1, N)  \big\}
		\leq  C  \P\{A(0;n), A(-f; N), S_n \in [0, 2\sqrt{n}] \}
		\\
		& \ =   C \sum\limits _{m =0 } ^{ 2\sqrt{n}  - 1}
		\P\{A(0;n), A(-f;n+1,  N), S_n \in [m, m+1) \}
		\\
		& \overset{\eqref{lithe = thin, supple,young, healthy}}{\leq}  \frac{C}{a+1} \sum\limits _{m =0 } ^{ 2\sqrt{n}  - 1}
		\P\{A(-a;n), A(-f-a;n+1,  N), S_n \in [m, m+1) \} 
		\\ 
		&\  =   \frac{C}{a+1}
		\P\{A(-a;n), A(-f-a;n+1,  N), 
		S_n \in [0, 2\sqrt{n}] \}.
	\end{align*}
	Since  $$\{A(-a;n), A(-f-a;n+1,  N), 
	S_n \in [0, 2\sqrt{n}] \} \subset \{A(-a;n), A(-f-a; N)  \},$$
	\eqref{in sich gehen} follows.
	\end{proof}

	\begin{proof}[
Proof of Theorem \ref{thm homing device = guidance system}]
	\label{Fuelle an, von} 
	Having carried out the preparatory work
	we now closely follow
	the proof of Theorem 3.2 in \cite{PP95}.
	We start with \eqref{a hill to die on}. Take $m \in \N$ and $N \geq 2 ^{m-1}$. Set
	$$\kappa _m = \inf \{k \in \N: k \geq 2^{m-1} + 1,  S_k \leq f(2^m) - a \}.$$
	Applying the strong Markov property, the independence of the increments of the random walk, and  Lemma \ref{chow = food} 
	we find 
	\begin{align}
		\P \big\{  T_{-a} > & 4N,
		\kappa _m \leq 2^m \big\}
		\notag
		\\ \notag
		\leq  & \, \P  \Big\{  
		T_{-a} >  2^{m-1}, \kappa _m \leq 2^m,
		S_j - S_{\kappa _m } \geq  -f(2^m), j = \kappa_m +1, ..., 4N
		\Big\}
		\\ \notag
		=  & \, \P  \Big\{    
		S_j - S_{\kappa _m } \geq  -f(2^m), j = \kappa_m +1, ..., 4N
		\Big|T_{-a} >  2^{m-1}, \kappa _m \leq 2^m
		\Big\}
		\\ &  \times  \P \{ T_{-a} >  2^{m-1}, \kappa _m \leq 2^m\}
		\\ \notag
		\leq   & \, \P  \Big\{    
		S_j \geq  -f(2^m), j = 1, ...,2N
		\Big\} \P \{ T_{-a} >  2^{m-1}\}
		\\ \notag
		= & \,
		\P \{  T_{- f(2^m)} > 2N   \}
		\P \{  T_{-a} > 2^{m-1} \}
		\\ \notag
		\leq & \, c_1 \frac{f(2^m)+1}{N ^{1/2}} \times  \frac{a+1}{2^{m/2}} .
	\end{align}
	Hence for $N \geq 2^{m-1}$
	\begin{equation}\label{oven mitts = oven gloves}
		\P \Big\{ \kappa_m \leq 2^m \big| T_{-a} > 4N \Big\}
		\leq  c_1 \frac{(a+1)(f(2^m) + 1)}{2^{m/2} N ^{1/2}} 
		\times
		\frac{N ^{1/2}}{c_2(a+1)} = c_3
		\frac{(f(2^m) + 1)}{2^{m/2} } 
	\end{equation}
	By Lemma \ref{whittle}, 
	$\sum _{m \in \N } {f(2^m)}{2^{-m/2}} < \infty$.
	Choose $m_f \in \N$ satisfying
	\begin{equation*}
		c_3 \sum\limits _{m = m_f }^\infty \frac{f(2^m)+1}{2^{m/2}}
		< \frac 12.
	\end{equation*}
	Note that for $k \in \N$, $2^{\lceil \ln _2 k \rceil}$
	is the smallest power of $2$ greater or equal than $k$.
	By \eqref{oven mitts = oven gloves} for 
	$M > m_f$ we have
	\begin{align*}
		\P \Big\{\exists k \in \N, 2^{m_f-1}< k \leq & \, 2^M: 
		S_k \leq f(2^{\lceil \ln _2 k \rceil}) - a \big|T_{-a}> 2^{M+1} \Big\}
		\\
		\leq & \, \sum\limits _{m = m_f} ^M 
		\P \Big\{\exists k \in \N, 2^{m-1}< k \leq  2^m: 
		S_k \leq f(2^m) - a \big| T_{-a}> 2^{M+1} \Big\}
		\\
		= & \, \sum\limits _{m = m_f} ^M 
		\P \Big\{ \kappa _m \leq  2^m| T_{-a}> 2^{M+1} \Big\}
		\\
		\leq &  \, c_3 \sum\limits _{m = m_f} ^M 
		\frac{(f(2^m) + 1)}{2^{m/2} }  \leq \frac 12.
	\end{align*}
	Hence
	\begin{equation*}
		\P \Big\{
		S_k \geq f(2^{\lceil \ln _2 k \rceil}) - a, 2^{m_f-1}< k \leq  2^M  \big|T_{-a}> 2^{M+1} \Big\} \geq \frac 12
	\end{equation*}
	and 
	\begin{equation*}
		\P \big\{ S_k \geq f(2^{\lceil \ln _2 k \rceil}) - a, 2^{m_f-1}< k \leq  2^M \big\}\geq
		\frac { \P\{T_{-a}> 2^{M+1} \} }{2}
		\geq c  (a+1) 2^{-M/2}.
	\end{equation*}
	Since $f(2^{\lceil \ln _2 k \rceil}) \geq f(k)$,
	$k \in \N$, 
	we get for $n \in \N$, $n>  2^{m_f-1}$
	\begin{multline*}
		\P \big\{ S_k \geq f(k) - a, 2^{m_f-1}< k \leq  n \big\}
		\geq
		\P \big\{ S_k \geq f(2^{\lceil \ln _2 k \rceil}) - a, 2^{m_f-1}< k \leq  2^{\lceil \ln _2 n \rceil} \big\}
		\\
		\geq 
		c  (a+1) \frac{1}{\sqrt{2^{\lceil \ln _2 n \rceil}}}
		\geq c  (a+1) \frac{1}{\sqrt{n}},
	\end{multline*}
	and \eqref{a hill to die on} follows.\\
	  
	The proof of \eqref{being a trooper}  too goes almost exactly 
	as in \cite{PP95}.
	Recall that $A(g;n)$ was defined on Page \pageref{letzten Endes = ultimately}. 
	The key part in the proof of \eqref{being a trooper} is the inequality for some $c > 0$
	(cf. \cite[(3.12)]{PP95})
	\begin{equation}\label{mellow - soft, smooth}
		\P \Big\{ T_{-a} < 2n \big| T_{-a} \geq  n ,  A(-f-a; N) \Big\} \leq c  \frac{f(3n)}{\sqrt{n}}
	\end{equation}
	which holds for large $n \in \N$
	and $N \geq 4n$, $ N \in \N$.
	Assuming \eqref{mellow - soft, smooth} we find
	for large $M \in \N$
	and $m_0 \leq M $
	\begin{align}
		\P \big\{ T_{-a} \geq & \, 2^{M+1}  |   A(-f-a; 2^{M+2}) \big\} 
		\notag 
		\\ 
		= & \,
		\P \big\{ T_{-a} \geq   2^{m_0}  |  A(-f-a; 2^{M+2}) \big\}  \notag
		\\
		& \,
		\times 
		\prod\limits _{m= m_0} ^M 
		\P \big\{ T_{-a} \geq   2^{m +1}  | T_{-a} \geq   2^{m },   A(-f-a; 2^{M+2})  \big\}  \notag
		\\
		\geq & \, \tilde c_{m_0}
		\prod\limits _{m= m_0} ^M 
		\Big(1 - c_6 \frac{f(3\cdot 2^m)}{2^{m/2}}\Big),
		\label{hackneyed = overused}
	\end{align}
	where $\tilde c_{m_0} =  \P \big\{ T_{-a} \geq   2^{m_0}  |   A(-f-a; 2^{M+2}) \big\} $
	and $m_0$ is chosen in such a way that every term
	on the right hand side of 
	\eqref{hackneyed = overused} is positive.
	The convergence \eqref{fib = unimportant lie}
	and Lemma \ref{whittle} imply 
	$$
	\prod\limits _{m= m_0} ^\infty 
	\Big(1 - c_6 \frac{f(3\cdot 2^m)}{2^{m/2}}\Big)
	> 0,
	$$
	and hence by \eqref{hackneyed = overused} for some $c_8>0$
	\begin{equation*}
		\P \big\{ T_{-a} \geq   2^{M+1}  |  A(-f-a; 2^{M+2}) \big\} \geq c_8.
	\end{equation*}
	Consequently $\P \big\{ T_{-a} \geq   2^{M+1} \big\}
	\geq c_8  \P \big\{   A(-f-a; 2^{M+2}) \big\} $ and by Lemma \ref{chow = food}
	\begin{equation*}
		\P \big\{   A(-f-a; 2^{M+2}) \big\}\leq c_8 ^{-1 }\P \big\{ T_{-a} \geq   2^{M+1} \big\}
		\leq \frac{c_9 (a+1) }{2^{M/2 }}
	\end{equation*}
	for large $M$. The inequality \eqref{being a trooper} follows.\\
	  
	It remains to show \eqref{mellow - soft, smooth}. To this end we 
	condition on $T_{-a} = k \in [n, 2n)$
	and $y \in (0, f(k)]$ to find 
	\begin{multline}\label{nach bestem Wissen}
		\P \big\{ n \leq T_{-a} < 2n, A(-f-a; N)  \big\}
		\\ \leq 
		\P \big\{ n \leq T_{-a} < 2n  \big\}
		\sup\limits_{\substack{n\leq k < 2n, \\
				0 < y \leq f(k) }} \P \big\{  
		A(-f-a; k,N)
		| S_{k}  = -a -y    \big\}
	\end{multline}
	We now use the inequality obtained
	with  the cutting and pasting argument   in \cite{PP95}
	on Pages 118-121 (it is the second inequality from below on Page 121):
	\begin{equation*}
		\sup\limits_{\substack{n\leq k < 2n, \\
				0 < y \leq f(k) }} \P \big\{   A(-f ; k,N) | S_{k}  =   -y    \big\}
		\leq c_{18} f(3n) \P \{ A(0;n), A(-f; N)  \},
	\end{equation*}
	or equivalently
	\begin{equation*}
		\sup\limits_{\substack{n\leq k < 2n, \\
				0 < y \leq f(k) }} \P \big\{   A(-f-a; k,N) | S_{k}  = -a -y    \big\}
		\leq c_{18} f(3n) \P \{ A(0;n), A(-f; N)  \}.
	\end{equation*}
	Therefore by
	\eqref{nach bestem Wissen}
	and
	Lemma \ref{perusal = examine sth} 
	\begin{multline}
		\P \Big\{ T_{-a} < 2n \big| T_{-a} \geq  n ,  A(-f-a;N)  \Big\} 
		= \frac{\P \Big\{ n \leq T_{-a} < 2n,  A(-f-a;N)  \Big\} }{\P \Big\{  T_{-a} \geq  n ,  A(-f-a;N)  \Big\}} 
		\\
		\leq 
		\frac{c_{18} f(3n) \P \big\{  T_{-a} \geq n \big\}
			\P \{ A(0;n), A(-f; N)  \}}{\P \Big\{ A(-a;n) ,  A(-f-a;k,N)  \Big\}} 
		\leq C f(3n) n^{-1/2},
	\end{multline}
	and \eqref{mellow - soft, smooth} follows.
	\end{proof}
	
	\subsection{Proof of Theorem \ref{thm beinhalten = include}}
	We need one more auxiliary lemma before we can proceed to the proof of the theorem.
	\begin{lem} \label{verdrehen = twist, distort}
		There is  $c > 0$ such that
		for large enough $n \in \N$ 
		\begin{equation}\label{sich zufriedengeben mit = settle for}
			\PP{ S_k \geq f(k) - a \text{ for } 1 \leq k \leq n, 
				S_{ n } \in  [\sqrt{3n}, 2\sqrt{3n}] }\geq
			\frac{ c(a+1)}{\sqrt{n}}.
		\end{equation}	 
	\end{lem}
	  
	\begin{proof}
	By 
	Theorem \ref{thm homing device = guidance system}
	there exists $n_f \in \N$
	and $c_1 > 0$ such that 
	for $n \geq n_f$ and $a \geq 0 $
	\begin{equation}\label{prowl copied}
		\PP{ S_k \geq f(k) - a \text{ for } n_f \leq k \leq n }\geq
		\frac{c_1(a+1)}{\sqrt{n}}.
	\end{equation}
	From here on this proof mirrors 
	the proof of Lemma \ref{der Zettel} - the only difference is that the level $-a$ is replaced by $f-a$. 
	Set again  $ \mathfrak{q} = \frac{\sqrt{3n}}{3}$.
	The events
	$  A(f- a; n)$ and $\{ T ^+_{4\mathfrak{q}} \leq n\}$
	are increasing events determined by the trajectory of the 
	random walk $(S_k)$. Therefore they are positively correlated
	by the FKG inequality, and by \eqref{prowl copied2}
	and \eqref{nachsagen = repeat2 copied2} we get for 
	${  B : = A( f- a; n)\cap \big\{ T ^+_{ 4\mathfrak{q}} \leq  n \big\}}$
	for large $n$
	
	\begin{equation}\label{zuversictlich = optimistic, confident copied}
		\P (B)  \geq \frac{ c (a+1) }{ \sqrt{n} }.
	\end{equation}
	Set $  B _1 = \big\{  T ^+_{ 4\mathfrak{q}} \leq n, S_k \geq f(k) - a, 1 \leq k \leq T ^+_{ 4\mathfrak{q}} \big\}$.
	Then $B \subset  B _1$
	and hence by \eqref{zuversictlich = optimistic, confident copied2} 
	\begin{equation}\label{Zaehler = numerator copied}
		\P(B_1)  \geq \frac{   (a+1) }{ \sqrt{n} }.
	\end{equation}
	Lemma \ref{zuneigung = affection}
	helps us to control the overshoot 
	$S_{T ^+_{4\mathfrak{q}}} - 4\mathfrak{q}$: for large 
	$n$
	\begin{equation}\label{Nenner = denominator copied}
		\P \Big\{ S_{T ^+_{4\mathfrak{q}}} - 4\mathfrak{q} \leq \mathfrak{q} \Big|  B_1 \Big\}
		\geq  \frac 23, \ \ \ n \in \N.
	\end{equation}
	Next define the event
	\begin{equation}
		B_2 =  \Big\{  T ^+_{ 4\mathfrak{q}} \leq n, |S_k - S_{ T ^+_{4\mathfrak{q}} }| \leq 
		\mathfrak{q}
		\text{ for }   T ^+_{ 4\mathfrak{q}} + 1 \leq k \leq n \Big\} 
	\end{equation} 
	For large $n$ the conditional probabilities of this event, given ${  B_1 \cap \{  T ^+_{ 4\mathfrak{q}}  = i\}}$,
	$1 \leq i \leq n$, are bounded away from $0$:
	for some $c_{20} > 0$ 
	\begin{equation}	\label{harken back to copied}
		\inf _{i = 1,...,n}\P \big\{ B_2 | B_1,   T ^+_{ 4\mathfrak{q}}  = i\big\} \geq c_{20}.	
	\end{equation}
	Note that again
	$$B_1 \cap \Big\{ S_{T ^+_{4\mathfrak{q}}} - 4\mathfrak{q} \leq \mathfrak{q}  \Big\} \cap B_2 \subset   \big\{  A(f-a; n), 
	S_{ n } \in  [\sqrt{3n}, 2\sqrt{3n}] \big\}.$$
	Given $B_1 \cap \big \{ 
	T ^+_{ 4\mathfrak{q}}  = i \big\} $, the events 
	$B_2$ and $\big\{ S_{T ^+_{4\mathfrak{q}}} - 4\mathfrak{q} \leq \mathfrak{q} \Big|  B_1 \big\}$
	are conditionally independent.
	Since 
	\begin{equation*}
		\P\big\{  A(f-a; \, n), 
		S_{ n } \in  [\sqrt{3n}, 2\sqrt{3n}] \big\}
		\geq   \P\big\{ B_2, S_{T ^+_{4\mathfrak{q}}} - 4\mathfrak{q} 
		\leq  \frac{\sqrt{n}}{3},  B_1 \big\},
	\end{equation*}
	the same transformations verbatim as in \eqref{der Bias} 
	give for some $c>0$ for large $n$
	\begin{equation*}
		\P\big\{  A(f-a; \, n), 
		S_{ n } \in  [\sqrt{3n}, 2\sqrt{3n}] \big\}
		\geq  c \P(B_1).
	\end{equation*}
	The statement of the lemma follows 
	from this inequality and  \eqref{Zaehler = numerator copied}.
	\end{proof}
	\begin{proof}
	[Proof of Theorem \ref{thm beinhalten = include}] 
	\label{Fuelle = abundance, wealth} 
	The convergence of the series in \eqref{fib = unimportant lie}
	and the monotonicity of $f$ imply that 
	${f(n) < \frac{\sqrt{n}}{10} }$ for all but finitely many $n \in \N $.
	Without loss of generality we assume 
	that $f(n) < \frac{\sqrt{n}}{10} $ for all $n \in \N $. 
	To simplify notation and without loss of generality
	we  consider only $n \in \N$
	divisible by $3$.
	We start with the lower bound in \eqref{verschleiern = disguise}. \\
	  
	By Lemma \ref{verdrehen = twist, distort}
	for some 
	$c > 0$
	for large $n \in \N$ 
	\begin{equation}\label{E1}
		\PP{ S_k \geq f(k) - a \text{ for } 1 \leq k \leq n/3, 
			S_{ n/3 } \in  [\sqrt{n}, 2\sqrt{n}] }\geq
		\frac{ c(a+1)}{\sqrt{n}}.
	\end{equation}
	and
	\begin{equation}\label{E2}
		\PP{\widetilde S_k \geq f(k) - b \text{ for } 1 \leq k \leq n/3, 
			\widetilde S_{ n/3 } \in  [\sqrt{n}, 2\sqrt{n}] }\geq
		\frac{ c(b+1)}{\sqrt{n}},
	\end{equation}
	where $c_5> 0$.
	The remaining part of the proof of
	the lower bound in \eqref{verschleiern = disguise}
	mirrors 
	the proof of Lemma \ref{entnehmen}, 
	except that here we deal with moving barriers.
	To ensure that 
	$S _n \in [-a + f_n(n)+ b , -a+ f_n(n) + b + 1 ]$
	as in the event in \eqref{verschleiern = disguise} 
	we adjust the middle part of the random walk $S_1,...,S_n$. 
	If $S_{n/3} = p$ and $\widetilde S_{n/3} = q$, $p,q \in [\sqrt{n}, 2\sqrt{n}]$,
	then $S _n \in [-a+f_n(n)+ b , -a+ f_n(n) + b + 1 ]$ 
	if and only if
	$${p + S_{2n/3} - S_{n/3} - q \in [-a+f_n(n)+ b ,-a+ f_n(n) + b + 1 ] }, $$
	or 
	$$
	S_{2n/3} - S_{n/3} \in [q - p -a+ f_n(n)+ b ,q - p -a+ f_n(n) + b + 1 ].
	$$
	Since $|f_n(n)| \leq f(n) < \frac{\sqrt{n}}{10}$,
	by the local limit theorem for some $c  >0 $
	uniformly in $p,q \in [\sqrt{n}, 2\sqrt{n}]$
	we have 
	\begin{equation}\label{glean = gather laboriously}
		\PP{ S_{2n/3} - S_{n/3} \in [q - p -a+ f_n(n)+ b ,q - p -a+ f_n(n) + b + 1 ]  } \geq \frac{c }{\sqrt{n}}.
	\end{equation}
	Since 
	$q - p + f_n(n)+ b \geq 0.8\sqrt{n} - p -a$, by the invariance principle for the conditioned sum
	of independent random variables 
	\cite[Theorem 4 and the corollary on page 568]{Lig68}
	for some $c  > 0$
	uniformly in $p,q \in [\sqrt{n}, 2\sqrt{n}]$
	for large $n \in \N$
	\begin{multline}
		\P \Big\{   
		S_{k} - S_{n/3} \geq \frac{\sqrt{n}}{2} - p -a, 
		n/3 \leq	k  \leq 2n/3 
		\\
		\big| S_{2n/3} - S_{n/3} \in [q - p -a + f_n(n)+ b ,q - p  -a +f_n(n) + b + 1 ]  \Big\} \geq c.
	\end{multline}
	Combining this with \eqref{glean = gather laboriously}
	gives 
	\begin{multline}\label{E3}
		\P\Big\{ S_{2n/3} - S_{n/3} \in [q - p -a + f_n(n)+ b ,q - p -a + f_n(n) + b + 1 ], \\
		S_{k} - S_{n/3} \geq \frac{\sqrt{n}}{2} - p -a, 
		n/3 \leq	k  \leq 2n/3   \Big\} \geq \frac{c }{\sqrt{n}}.
	\end{multline}
	Denote by $E_1$, $E_2$, and  $E_3(p,q)$
	the events in the left hand sides
	of \eqref{E1}, \eqref{E2},
	and \eqref{E3}, respectively.
	The events $E_1$, $E_2$, and  $E_3(p,q)$
	are independent 
	as they are defined in terms of
	disjoint parts of the increments of the random walk.
	By construction and definition of $f$, 
	\begin{multline}\label{zugunsten + gen}
		\P\{S_k \geq -a + f _n(k), k \leq n, S _n \in [ - a + f_n(n)+ b , -a  + f_n(n) + b + 1 ]  \}
		\\
		\geq 
		\P(E_1 \cap E_2 \cap E_3 (S_{n/3}, \widetilde S_{n/3} ))
		=  \P(E_1) \P( E_2)  \E \big[\P (E_3(S_{n/3}, \widetilde S_{n/3} ) | E_1, E_2)\big]
		\\
		\geq \frac{c(a+1)(b+1)}{n^{3/2}}.
	\end{multline}
	Now let us turn to the upper bound in \eqref{verschleiern = disguise}.
	Define the events
	\begin{gather*}
		D_1 = \{   
		S_k \geq - f(k) - a \text{ for } 1 \leq k \leq n/3 
		\},
		\\
		D_2 = \{   
		\widetilde S_k \geq - f(k) - b - 1 \text{ for } 1 \leq k \leq n/3
		\},
	\end{gather*}
	and for $p, q \in \R$
	\begin{equation}
		D_3 (p,q) = \{ 
		S_{2n/3} - S_{n/3} \in [q - p -a + f_n(n)+ b ,q - p-a + f_n(n) + b + 1 ]
		\}.
	\end{equation}
	Note that 
	if $S _n \in [-a + f_n(n)+ b , -a+ f_n(n) + b + 1 ]$
	then for   $2n/3\leq k \leq n$
	\begin{multline}
		S_k \geq -a + f _n(k) \Rightarrow 
		S_k - S_n \geq -a + f _n(k) - S_n 
		\\
		\Rightarrow \widetilde S_{n-k} \geq 
		-a + f _n(k) -  (-a+ f_n(n) + b + 1)
		\Rightarrow  \widetilde S_{n-k} \geq 
		+ f _n(k) -   f_n(n) - b - 1
		\\
		\Rightarrow  \widetilde S_{n-k} \geq 
		-   f(n-k) - b - 1.
	\end{multline}
	Therefore
	by Theorem \ref{thm homing device = guidance system} 
	and Lemma \ref{countenance = facial expression}
	\begin{multline*}
		\P\{S_k \geq -a + f _n(k), k \leq n, S _n \in [-a + f_n(n)+ b , -a+ f_n(n) + b + 1 ]  \}
		\\
		\leq 
		\P(D_1 \cap D_2 \cap D_3(S_{n/3}, \widetilde S_{n/3} ))
		=  \P(D_1) \P( D_2)  \E \big[\P (D_3(S_{n/3}, \widetilde S_{n/3} ) | D_1, D_2 )\big]
		\\
		\leq \frac{C(a+1)(b+1)}{n^{3/2}}.
	\end{multline*}
	\end{proof}

	The next Lemma is used
	on Page \pageref{pageref Lemma deign}
	in the proof of Lemma \ref{hoodwink = deceive, trick2}.
	\begin{lem}\label{deign}
		In the setting of Theorem \ref{thm beinhalten = include}
		denote by $B$ the event in the middle part 
		of  
		\eqref{verschleiern = disguise} and let $s\in \{1,2,...,n\}$.
		Then a.s.
		\begin{equation}\label{plaid}
			\P \big( B \big| S  _s - S _{s-1} \big)
			\leq C \big(  |S _s - S _{s-1}|^4+1 \big)
			\P (B ).
		\end{equation}
	\end{lem}
	\begin{proof}
		Consider the trajectory of the random walk $\{{S}_k, 1\leq k\leq n\}$
		without the $s$-th step: 
		\begin{equation}
			\breve{S}_k = \begin{cases}
				S_k, & k \leq s-1,
				\\
				S_{k+1} - (S_s - S_{s-1}), &  s\leq k \leq n-1. 
			\end{cases}
		\end{equation}	
		Define also $\breve f_{n-1}:\{0,1,...,n-1\} \to \R $ by
		\begin{equation}
			\breve{f}_{n-1} (k) = \begin{cases}
				{f}_{n} (k), & k \leq s-1,
				\\
				{f}_{n} (k+1), &  s\leq k \leq n-1. 
			\end{cases}
		\end{equation}	
		Note that the trajectory of $\{\breve{S}_k, 1\leq k\leq n-1\}$ is obtained from $\{{S}_k, 1\leq k\leq n\}$ by removing the $s$-th step. We have
		\begin{multline}
			B \subset \breve{B} :=
			\Big\{ 
			\breve{S}_k \geq \breve{f}_{n-1} (k) - a - |S_s-S_{s-1}|, 
			\breve{S}_{n-1} \in 
			\big[f_n(n) -a - (S_s-S_{s-1})+ b , 
			\\ f_n(n)-a - (S_s-S_{s-1}) + b + 1 \big] 
			\Big\}.
		\end{multline}
		Since $S_s-S_{s-1}$ and $\{\breve{S}_k, 1\leq k\leq n-1\}$ are independent,
		we see by 
		Theorem \ref{thm beinhalten = include} that
		a.s. on $\{S_s-S_{s-1}  \leq n ^{0.4} \}$
		\begin{equation}
			\P \big( B \big| S  _s - S _{s-1} \big) \leq
			\frac{C (a+ |S_s-S_{s-1}| + 1) (b + 2 |S_s-S_{s-1}| + 1)}{n^{3/2} }
			\leq C \P (B) [(S_s-S_{s-1})^2 + 1].
		\end{equation}
		On the other hand, 
		by 
		Theorem \ref{thm beinhalten = include}
		there exists $\widetilde C>0$  such that
		for all $n \in \N$
		$$\widetilde C( n ^{0.4} + 1) P(B) \geq 1.$$ 
		For such $C = \widetilde C$ inequality  \eqref{plaid} trivially holds
		a.s. on $\{S_s-S_{s-1}  > n ^{0.4} \}$.
	\end{proof}

\begin{appendix}
\section*{}
We used the following classical result on the distribution of 
the projection of a random point on the unit sphere. It is attributed to E. Borel, Introduction géométrique à quelques théories physiques. Gauthier-Villars, Paris, 1914 in \cite{Spruill07}. We also refer to \cite{thesisYin-Ting}, Lemma 3.3.1., for a proof which is available online.

\begin{lem}\label{gully = ravine}
	Let $\mathcal{H} = (\mathcal{H}_1, \mathcal{H}_2, ..., \mathcal{H}_\d)$ be a uniformly chosen random point on
	the unit sphere $\S ^{\d-1}$.
	Then the density of $\mathcal{H}_1$
	is given by 
	\begin{equation}\label{onesie}
		f _{ \mathcal{H}_1 } (x)=
		\begin{cases} \frac{1}{B(\frac 12, \frac{\d-1}{2})} 
			(1 - x^2) ^{\frac{\d-3}{2}} & x \in [-1, 1]\\
			0 & \mbox{ otherwise}.
		\end{cases}
	\end{equation}
Here, $B(\cdot, \cdot)$ is the beta function given by
\begin{equation}
B(\alpha, \beta): = \int\limits_0^1 x^{\alpha -1}(1-x)^{\beta - 1} dx
\end{equation}
for $\alpha, \beta > 0$.
\end{lem}

   Let $F$ be the cumulative distribution function of $|Q_1|$.
   The pair $(\X_1, |Q_1|)$ has the same distribution as
   $(\mathcal{H}_1|Q_1|, |Q_1|).$ 
By Lemma \ref{gully = ravine}
\begin{align}
	\P\{\X_1 \in A, & \, |Q_1| \in \mathcal{Q}  \}
	\notag
	\\
	\label{take a rain check=decline}
	 = & \, \frac{1}{B\big(\frac 12, \frac{\d-1}{2} \big)} 
	\int\limits _{y \in A, r \in \mathcal{Q}}
	\frac {\1 \{ |y| \leq r \} \1 \{ r > 0  \} }{r}
	\Big(1 - \frac{y^2}{r^2}\Big)^{\frac{\d-3}{2}}dy F(dr)
	\\
 &	+ \1\{0 \in A, 0 \in \mathcal{Q} \} F(\{0\}).
 \notag
\end{align}
Hence for a Borel set $\mathcal{Q} \subset \R $
and $a \geq 0  $  
\begin{equation}\label{cond distr}
	\P \big\{ |Q_1| \in \mathcal{Q}  | \X_1 = - a \big\}
	=
	\P \big\{ |Q_1| \in \mathcal{Q}  | \X_1 = a \big\}
	=  \frac{\int\limits _{\mathcal{Q} \cap [a,\infty)}   r^{-1}(1 - \frac{a^2}{r^2})^{\frac{\d-3}{2}} F(dr)  }
	{\int\limits _{a} ^{\infty}  r^{-1}(1 - \frac{a^2}{r^2})^{\frac{\d-3}{2}} F(dr)  }.
\end{equation}

\begin{proof}
[Proof of Lemma \ref{lem many-to-two}]
\label{proof of many-to-two} For $u,v \in \T$ denote by 
$u \wedge v$ their most recent common ancestor 
in $\T$; $u \wedge u = u$. Grouping the terms of  the sum by 
the generation of the most recent common ancestor 
we obtain 
\begin{align}
	\sum\limits _{ |v| = |u| = n}
	f(X(v_0), & \,  ..., X(v_n), X(u_0), ..., X(u_n) ) \notag
	\\  \label{upsell = persuade to buy more expensive}
	= & \,
	\sum\limits_{k=0}^n \sum\limits _{\substack{ |v| = |u| = n: \\ |u \wedge v| = k}}
	f(X(v_0),  ..., X(v_n), X(u_0), ..., X(u_n) ). 
\end{align}
For the term with $k = n$ in \eqref{upsell = persuade to buy more expensive} we have $u = v$ and by the many-to-one lemma 
\begin{equation}
	\E\Big[ \sum\limits _{\substack{  |u| = n}}
	f(X(u_0),  ..., X(u_n), X(u_0),  ..., X(u_n) )\Big] = 
	m^n\E \Big[f(Q_0,  ..., Q_n, Q_0 ,  ..., Q_n)\Big].
\end{equation}
Consider now the terms with $k < n$
in the sum on the right hand side 
of \eqref{upsell = persuade to buy more expensive}.
For 
$a,b  \in \T $ we write $a \prec b$ 
if $a$ is an ancestor of $b$ and $a \ne b$;
we write $a \preccurlyeq b$
if either $a \prec b$  or $a = b$.
Set $w = u \wedge v$
and for $a \in \T$ denote by $\mathcal{C} _a$ the set of children of $a$.
We have 
\begin{multline}\label{impertinent = rude}
	\sum\limits _{\substack{ |v| = |u| = n: \\ |u \wedge v| = k}}
	f(X(v_0),  ..., X(v_n), X(u_0), ..., X(u_n) )
	\\
	= 
	\sum\limits _{\substack{ |w| =  k}}
	\sum\limits _{\substack{ a,b \in  \mathcal{C}_w: \\
			a \ne b}}
	\sum\limits _{\substack{|u| = n: \\ a \preccurlyeq u}}
	\sum\limits _{\substack{|v| = n: \\ b \preccurlyeq v}}
	f(X(v_0),  ..., X(v_n), X(u_0), ..., X(u_n) ),
\end{multline}
where we follow the usual convention that the sum 
over the empty set is set to zero.
Define $g_k: (\R ^\d)^k \to \R$ by 
\begin{equation}\label{sycophant}
	g_k(x_0,...,x_k) = \E \Big[\begin{aligned}[t] f(x_0,...,x_k, x_k + \Delta _1, x_k + \Delta _1+Q'_1,  ..., x_k + \Delta _1+Q'_{n-k-1},
		\\
		x_0,...,x_k, x_k + \Delta _2, x_k + \Delta _2+Q''_1,  ..., x_k + \Delta _2+Q''_{n-k-1})\Big].
	\end{aligned}
\end{equation}
Denote by $\mathcal{F} _t$ the 
$\sigma$-algebra generated by the the spatial positions
of the particles
of the BRW $(\T, \{X(v) \mid v \in \T \})$ up to time $t$:
$$
\mathcal{F} _t = \sigma \{ X(v): |v| \leq t \},
$$
and by $\emptyset $ the root of $\T$. Also,
denote by $N_t$ the number of particles in generation $t $, i.e. $N_t = \# \T_t$.
For any bounded measurable $h$ and $k, \ell \in \N$,
$k < \ell$ 
and  $w \in \T $ with $|w| = k$
we find, almost surely,
\begin{align}
& \E 	\Big[ \sum\limits _{\substack{ a,b \in  \mathcal{C}_w: \\
		a \ne b}} 
\sum\limits _{\substack{|u| = \ell: \\ a \preccurlyeq u}}
\sum\limits _{\substack{|v| = \ell: \\ b \preccurlyeq v}}
h(X(v_0),   ..., X(v_\ell), X(u_0), ..., X(u_\ell) ) \Big| 
\mathcal{F}_k  \Big]
\notag 
\\
\notag 
= & \,
m^{2 \ell -2}
\E \bigg[\sum\limits _{\substack{ a,b \in  \mathcal{C}_w: \\
		a \ne b}}
h\Big(X(v_0), ..., X(w),  X(a), X(a) + Q_1, X(a)+Q_2, ..., X(a)+Q_{\ell - k -1},
\\ 
& \hspace{3.1cm} X(v_0), ..., X(w), X(b), X(b) + Q'_1, X(b)+Q'_2, ..., X(b)+Q'_{\ell - k -1} \Big) \Big| 
\mathcal{F}_k  \bigg]
\notag 
\\
= & \,  
m^{2 \ell -2} {m_2} \E\Big[
h(X(v_0), ..., X(w), X(w)+  \Delta _1, X(w) + \Delta _1 + Q'_1,  ..., X(w) +\Delta _1+Q'_{\ell - k -1}, 
\\
\label{sleuth = detective}
& \hspace{2.1cm}X(v_0), ..., X(w), X(w)+ \Delta _2 + Q''_1, ..., X(w)+ \Delta _2+Q''_{\ell - k -1} ) \Big| 
\mathcal{F}_k  \Big].
\end{align}
In the last step in \eqref{sleuth = detective}
the distribution of $\Delta$ given in  \eqref{morose, sulky, sullen}
was used.
Conditioning on $\mathcal{F}_k$
we get by 
\eqref{sleuth = detective}
and
the definition of $g_k$  in \eqref{sycophant}
\begin{multline}\label{prevaricate = speak or act in an evasive way}
\E\Big[\sum\limits _{\substack{ |w| =  k}}
\sum\limits _{\substack{ a,b  \in \mathcal{C}_w: \\
		a \ne b}}
\sum\limits _{\substack{|u| = n: \\ a \preccurlyeq u}}
\sum\limits _{\substack{|v| = n: \\ b \preccurlyeq v}}
f(X(v_0),  ..., X(v_n), X(u_0), ..., X(u_n) ) \Big]
\\
\begin{aligned}[t]
	= & \,m^{2 n - 2k -2} m_2 \E\Big[\sum\limits _{\substack{ |w| =  k}} 
	g_k( X(v_0), X(v_1), ..., X(v_k))\Big]
	\\
	= & \, m^{2 n - k -2} m_2 \E\Big[ g_k( Q_0, Q_1, ..., Q_k)\Big].
\end{aligned} 
\end{multline}
By construction 
$$
\E \big[g_k( Q_0, Q_1, ..., Q_k)\big] = \E \big[f(Q_0 ^{\langle k \rangle}, Q_1 ^{\langle k \rangle}, ..., Q_n ^{\langle k \rangle}, Q^{[k]} _0  ,  Q^{[k]} _1 , ..., Q^{[k]} _n )\big],
$$
therefore \eqref{Skrei} follows from \eqref{prevaricate = speak or act in an evasive way} by summing over $k$.
\end{proof}

\end{appendix}

\begin{acks}[Acknowledgments]
 We thank Bastien Mallein for discussions and for pointing us to Lemma \ref{lem many-to-two}. We are grateful to two anonymous referees for reading very carefully and suggesting many improvements.
\end{acks}

\bibliographystyle{imsart-nameyear} 
\bibliography{Sinus}       

\end{document}